\documentclass[11pt]{article}

\usepackage{fancyhdr}
\usepackage{amsfonts}
\usepackage{amsmath}
\usepackage{amssymb}
\usepackage{amsthm}
\usepackage{tikz}

\usepackage{amscd}
\usepackage{amstext}

\usepackage{oldgerm} 

\usepackage{fancyhdr}
\usepackage{amsfonts}
\usepackage{amsmath}
\usepackage{amssymb}
\usepackage{amsthm}

\usepackage{amscd}
\usepackage{amstext}
\usepackage{color}

\usepackage{graphics}
\usepackage{graphicx}

\let\svthefootnote\thefootnote
\newcommand\freefootnote[1]{%
	\let\thefootnote\relax%
	\footnotetext{#1}%
	\let\thefootnote\svthefootnote%
}

\newlength{\fixboxwidth}
\newcommand{\com}{\mathbb{C}}

\newcommand{\cs}{{\mathcal S}}

\newcommand{\cf}{{\mathcal F}}
\newcommand{\cfi}{{\cf}^{-1}}

\newcommand{\supp}{{\rm supp \, }}
\newcommand{\sgn}{{\rm sgn \, }}

\newcommand{\be}{\begin{equation}}
\newcommand{\ee}{\end{equation}}
\newcommand{\beq}{\begin{eqnarray}}
\newcommand{\beqq}{\begin{eqnarray*}}
\newcommand{\eeq}{\end{eqnarray}}
\newcommand{\eeqq}{\end{eqnarray*}}

\newtheorem{satz}{Theorem}
\newtheorem{cor}{Corollary}
\newtheorem{rem}{Remark}
\newtheorem{defi}{Definition}
\newtheorem{lem}{Lemma}
 
\newtheorem{prop}{Proposition}

\begin{document}


\title{Quarklet Characterizations for bivariate Bessel-Potential Spaces on the Unit Square via Tensor Products}

\author{Marc Hovemann$^{\ast}$}

\date{\today}

\maketitle

\freefootnote{
	\hspace{-21.5pt}
	This work has been supported by Deutsche Forschungsgemeinschaft (DFG), grant HO 7444/1-1 with project number 528343051.

\noindent	
${}^\ast$Philipps-Universit{\"a}t Marburg, Fachbereich Mathematik und Informatik,   Hans-Meerwein Str. 6,  Lahnberge, 35043 Marburg, Germany, Email: \texttt{hovemann@mathematik.uni-marburg.de}.	
	}


\vspace{0.2 cm }


\vspace{0.2 cm }

\textbf{Abstract.}  In this paper we deduce new characterizations for bivariate Bessel-Potential spaces defined on the unit square via B-spline quarklets. For that purpose in a first step we use univariate boundary adapted quarklets to describe univariate Bessel-Potential spaces on intervals. To obtain the bivariate characterizations a recent result of Hansen and Sickel is applied. It yields that each bivariate Bessel-Potential space on a square can be written as an intersection of function spaces which have a tensor product structure. Hence our main result is a characterization of bivariate Bessel-Potential spaces on squares in terms of quarklets that are tensor products of univariate quarklets on intervals.

\section{Introduction}

Many problems in science, engineering and finance can be modeled as partial differential equations. Often a closed form of the unknown solution is not available, so that numerical schemes for its constructive approximation are indispensable tools. A very popular approach is the finite element method (FEM). The classical $h$-FEM is based on a space refinement of the domain under consideration. Alternatively it is also possible to increase the polynomial degree of the ansatz functions, which is known as the $p$-method. A combination of both, so-called $hp$-FEM, also can be applied. To handle large-scale problems it is essential to employ adaptive strategies to increase the overall efficiency. Then the goal is to end up with a satisfying approximation in a reasonable time. Especially for adaptive $h$-FEM there is a huge amount of literature. Let us refer to \cite{Ci}, \cite{Ha}, \cite{NSV09}, \cite{Sch} and \cite{Verf}. Nowadays the convergence analysis of adaptive $p$- and $hp$-methods is more and more in the focus of research. These schemes converge very fast and often show exponential convergence. However, when it comes to the theoretical analysis and to rigorous convergence proofs only a few results have been found recently. Some state of the art results concerning the convergence of $hp$-adaptive strategies are \cite{bib:BPS13,bib:BD11,bib:DV20,bib:DH07} and \cite{bib:CNSV14,bib:CNSV17}, which also include optimality results. 

Another approach is to use wavelets. The advantages of wavelets are their very strong analytical properties that can be used to attain adaptive schemes that are guaranteed to converge with the optimal convergence order of the best $N$-term wavelet approximation. We refer to \cite{CaDaDe} and \cite{bib:Ste09}. Adaptive wavelet schemes are essentially space refinement methods and can therefore be classified as $h$-methods. Then a very natural question is whether $hp$-versions of adaptive wavelet methods can be designed. At this point the concept of \emph{quarklets} comes into play. These polynomially enriched wavelets have been introduced in the last decade, see \cite{DaKRaa}. They are constructed out of biorthogonal compactly supported Cohen-Daubechies-Feauveau spline wavelets, where the primal generator is a cardinal B-spline. For the theory of such biorthogonal wavelets we refer to Section 6.A in \cite{CoDau}. Roughly speaking the quarklets are linear combinations of translated cardinal B-splines that are multiplied with some monomials. A precise definition can be found below in Definition \ref{def_quarklet}. The properties of the quarklets have been studied in \cite{DaFKRaa,DaHoRaVo,DaKRaa,DaRaaS,HoDa,HoKoRaVo,SiDiss}. In particular they can be used to design schemes that resemble $hp$-versions of adaptive wavelet methods. It is our long term goal to design an adaptive scheme for solving partial differential equations based on quarklets that can be proven to converge and which realizes the optimal (exponential) rate. This paper can be seen as a further step in this direction.

One very interesting property of quarklets is that they can be used to characterize certain function spaces. So they have been utilized to describe univariate Sobolev, Besov and Triebel-Lizorkin spaces on the real line, see \cite{SiDiss} and \cite{HoDa}. Such results are important for numerical applications, in particular since usually the regularity of the solution in Bessel-Potential or Triebel-Lizorkin spaces describes the approximation order that can be achieved by adaptive numerical schemes. So nowadays it is well-known that certain Besov spaces are closely related to approximation spaces for N-term wavelet approximation, see \cite{DeJaPo} and Chapter 7.1 in \cite{CoSchn}. However, it has been observed in \cite{DaHaSchSi}, that at least for some PDEs sharper results can be obtained, when studying the regularity of the solution in the context of Bessel-Potential and Triebel-Lizorkin spaces. When we deal with the multivariate setting up to now only $L_{2}(\Omega)$ and the Sobolev spaces $H^{s}(\Omega) = H^{s}_{2}(\Omega)$ with $s \geq 0$ have been described via multivariate quarklets, at which $\Omega \subset \mathbb{R}^d $ is either a unit cube or a domain that can be decomposed into cubes. Here we refer to Section 3 in \cite{DaFKRaa} for the details. In other words so far only function spaces that are Hilbert spaces have been characterized by multivariate quarklets. Function spaces that are more difficult have not been investigated up to now. In this paper it is our main goal to overcome this hurdle. Indeed, in our main result Theorem \ref{THM_d2_Hsr2_main1} we will see, that it is possible to describe bivariate Bessel-Potential spaces $H^{s}_{r}((0,1)^2)$ with $ 1 < r < \infty  $ in terms of bivariate tensor quarklets. In connection with that also a new equivalent norm will be obtained. In order to prove this result we require two main tools which are also of interest on their own. On the one hand we derive a new characterization of univariate Bessel-Potential spaces $H^{s}_{r}((0,1))$ on bounded intervals in terms of quarklets, see Theorem \ref{mainresult1_Hsr_d=1}. For that purpose we deal with special boundary adapted quarklets which have been introduced in \cite{DaFKRaa}. On the other hand we use a recent result of Hansen and Sickel, see \cite{HaSi2022}. It tells us that bivariate Bessel-Potential spaces defined on squares can be written as an intersection of function spaces that are tensor products of univariate spaces, see also Proposition \ref{prop_tensor_Hsr_d2} below. A combination of both auxiliary results then yields our main theorem, namely the characterization of bivariate Bessel-Potential spaces via bivariate quarklets. Indeed, the bivariate quarklets we use for this purpose are tensor products of univariate boundary adapted quarklets, see Theorem \ref{THM_d2_Hsr2_main1} for the details. Due to the particular method we utilized for their construction, they have a highly anisotropic structure. In general, our tensor product quarklets can be interpreted as a wavelet version of sparse grids including polynomial enrichment. Therefore in the long run they most likely can be applied to design numerical approximation techniques with dimension-independent convergence rates.  

This paper is organized as follows. In Section \ref{sec_def} we recall the definition of Bessel-Potential spaces. Moreover, we explain the construction of univariate quarklets on the real line and on bounded intervals. In Section \ref{sec_bessel_d1} we prove a characterization in terms of quarklets for univariate Bessel-Potential spaces on bounded intervals. In Section \ref{sec_tensor} some basic facts concerning tensor products are collected. Furthermore, we recall a recent result of Hansen and Sickel which tells us that bivariate Bessel-Potential spaces on squares can be written as an intersection of function spaces that have a tensor product structure, see \cite{HaSi2022}. In Section \ref{sec_bessel_multi} we prove that bivariate Bessel-Potential spaces defined on unit squares can be characterized in terms of bivariate quarklets with tensor product structure. In connection with that we obtain the main result of this paper, see Theorem \ref{THM_d2_Hsr2_main1}. 

First of all we introduce some notation. As usual $\mathbb{N}$ denotes the natural numbers, $\mathbb{N}_0$ the natural numbers including $0$, $\mathbb{Z}$ the integers and $\mathbb{R}$ the real numbers. The number $d \in \mathbb{N}$ always refers to the dimension. Often we have $d = 1$ or $ d=2$. All functions in this paper are assumed to be complex-valued, that means we consider functions $f:~ \mathbb{R}^d \to \com$. Let $\mathcal{S}(\mathbb{R}^d)$ be the collection of all Schwartz functions on $\mathbb{R}^d$ endowed with the usual topology and denote by $\mathcal{S}'(\mathbb{R}^d)$ its topological dual, namely the space of all bounded linear functionals on $\mathcal{S}(\mathbb{R}^d)$ endowed with the weak $\ast$-topology. The symbol $\cf$ refers to  the Fourier transform, $\cfi$ to its inverse transform, both defined on $\cs'(\mathbb{R}^d)$. Almost all function spaces which we consider in this paper are subspaces of $\cs'(\mathbb{R}^d)$, namely spaces of equivalence classes with respect to almost everywhere equality. However, if such an equivalence class  contains a continuous representative, then usually we work with this representative and call also the equivalence class a continuous function. By $C^\infty_0(\mathbb{R}^d) = D(\mathbb{R}^d) $ we mean the set of all infinitely often differentiable functions on $\mathbb{R}^d$ with compact support. For $ 0 < r \leq \infty $ by $ L_{r}(\mathbb{R}^d)  $ we denote the usual Lebesgue spaces. Given a function $ f \in  L_{r}(\mathbb{R}^d) $ we use the symbol $ \Vert f \vert L_{r}(\mathbb{R}^d) \Vert   $ for the associated quasi-norm. When we have $ f, g \in L_{2}(\mathbb{R}^d)  $ we use the abbreviation
\begin{align*}
\left\langle f , g     \right\rangle_{L_{2}(\mathbb{R}^d)} = \int_{\mathbb{R}^d} \overline{f(x)} g(x) dx .
\end{align*}
For $ j, j' \in \mathbb{Z}   $ the symbol  $ \delta_{j , j'}   $ refers to the Kronecker delta. The symbols  $C, C_1, c, c_{1} \ldots $ denote  positive constants that depend only on the fixed parameters and probably on auxiliary functions. Unless otherwise stated their values may vary from line to line. In connection with that by $ A \lesssim B   $ we mean $ A \leq C B   $. Moreover, $  A \approx B   $ stands for $ A \lesssim B  $ and $ B \lesssim A   $. For $ j \in \mathbb{Z}  $ and $ k \in \mathbb{Z}  $ we define the dyadic intervals $ I_{j,k} $ via $ I_{j,k} := 2^{-j} ( [0,1) + k   )$. By $ \chi_{j,k}  $ we denote the characteristic function of an interval $ I_{j,k} $. Moreover, for $ j \in \mathbb{N}  $ and $ k \in \mathbb{Z}  $ we define 
\[
\tilde{\chi}_{j,k} := \left\{ \begin{array}{lll}
 \chi_{j,k}  & \quad & \mbox{for} \qquad k \in \{ 0, 1, \ldots , 2^{j} - 1 \} ;
\\  
\chi_{j,0} & \quad & \mbox{for}\qquad k < 0  ;
\\
\chi_{j,2^{j}-1} & \quad & \mbox{for}\qquad k > 2^{j}-1  .
\\
\end{array}
\right.
\]

\section{Basic Definitions: Function Spaces and Quarklets}\label{sec_def}

\subsection{Bessel-Potential Spaces on $\mathbb{R}^d$ and on Domains}

It is the main goal of this paper to prove quarklet characterizations for bivariate Bessel-Potential spaces defined on $ (0,1)^2  $. For that purpose in a first step we recall the definition of Bessel-Potential spaces on $ \mathbb{R}^d  $.

\begin{defi}\label{def_bes_pot}
\textbf{Bessel-Potential Spaces.} 

\noindent
Let $ s \in \mathbb{R}  $ and $ 1 < r < \infty  $. Then the Bessel-Potential space $ H^{s}_{r}(\mathbb{R}^{d})   $ is the collection of all $  f \in \mathcal{S}'(\mathbb{R}^{d})  $, such that $ \mathcal{F}^{-1} [(1 + \vert \xi \vert^{2} )^{s/2} \mathcal{F} f(\xi)](\cdot)    $ is a regular distribution and
\begin{align*}
\Vert f \vert H^{s}_{r}(\mathbb{R}^{d}) \Vert := \Big \Vert  \mathcal{F}^{-1} [(1 + \vert \xi \vert^{2} )^{\frac{s}{2}} \mathcal{F} f(\xi)](\cdot)   \Big \vert  L_{r}(\mathbb{R}^{d})  \Big \Vert  < \infty .
\end{align*}
\end{defi} 
\noindent 
The Bessel-Potential spaces are generalizations of the Sobolev spaces $W^{m}_{r}(\mathbb{R}^{d})$. So for $ 1 < r < \infty   $ and $ m \in \mathbb{N}   $ we find $ H^{m}_{r}(\mathbb{R}^{d}) =  W^{m}_{r}(\mathbb{R}^{d})   $ in the sense of equivalent norms, see the theorem in Chapter 2.5.6. in \cite{Tr83}. Therefore the spaces $ H^{s}_{r}(\mathbb{R}^{d})  $ sometimes are called fractional Sobolev spaces. Moreover, for $ 1 < r < \infty $ and $ s \in \mathbb{R} $ the Bessel-Potential spaces coincide with the Triebel-Lizorkin spaces $F^{s}_{r,2}(\mathbb{R}^d)$, that means we have 
\begin{equation}\label{eq_bes_eq_Tri}
H^{s}_{r}(\mathbb{R}^{d}) = F^{s}_{r,2}(\mathbb{R}^d)
\end{equation}
with equivalent norms, see Chapter 2.5.6 in \cite{Tr83}. For more details concerning Triebel-Lizorkin spaces we refer to \cite{Tr83}, \cite{Tr92}, \cite{Tr06} and \cite{Tr20} at least. It is possible to define Bessel-Potential spaces (and also Sobolev and Triebel-Lizorkin spaces) on domains $\Omega \subset \mathbb{R}^d$. Here domain means open connected set. For the definition we follow \cite{Tr06}, see Definition 1.95. 

\begin{defi}\label{def_bespot_dom}
\textbf{Bessel-Potential Spaces on Domains.} 

\noindent
Let $ s \in \mathbb{R}  $ and $ 1 < r < \infty  $. Let $\Omega \subset \mathbb{R}^d $ be a domain.
Then  $H^{s}_{r}(\Omega)$ is defined as the collection of all $f \in D' (\Omega)$, such that
there exists a distribution $g \in H^{s}_{r}( \mathbb{R}^{d} )$ satisfying
\begin{align*}
f (\varphi) = g (\varphi) \qquad \mbox{for all} \qquad \varphi \in D(\Omega) \, .
\end{align*}
Here $\varphi \in D(\Omega)$ is extended by zero on $ \mathbb{R}^{d} \setminus \Omega$.
We put
\begin{align*}
\| \, f\, \vert H^{s}_{r}(\Omega) \| := \inf \Big\{\| \, g\, \vert H^{s}_{r}( \mathbb{R}^{d} ) \|: \quad g_{|_\Omega} =f  \Big\} \, .
\end{align*}
\end{defi}
\noindent
In what follows almost all domains that we investigate are either the unit interval $(0,1)$ if $d = 1$ or the unit cube $(0,1)^d$ if $d \in \mathbb{N} $. It is possible to describe Bessel-Potential spaces which are defined on $(0,1)^d$ in terms of differences of higher order. For that purpose we use the following notation. Let $ f : \mathbb{R}^d \rightarrow \mathbb{C}$ be a function and let $ x, h \in \mathbb{R}^d $. Then the difference of the first order is defined as
\begin{align*}
\Delta_{h}^{1}f (x) := f ( x + h ) - f (x).
\end{align*}
For $ N \in \mathbb{N}  $ with $  N \geq 2  $ the difference of order $ N $ is given by
\begin{align*}
\Delta_{h}^{N}f (x) := \left  (\Delta_{h} ^1 \left ( \Delta_{h} ^{N-1}f \right  )\right ) (x) .
\end{align*}
Now we are prepared to formulate the following characterization. 

\begin{prop}\label{satz_diff1}
Let $\Omega = (0,1)^d$ with $d \in \mathbb{N}$. Let $1 < r < \infty$, $1 \leq v \leq \infty$, $N \in \mathbb{N}$ and 
\begin{align*}
d \max \Big ( 0 , \frac{1}{r} - \frac{1}{v} , \frac{1}{2} - \frac{1}{v} \Big ) < s < N .
\end{align*}
Then $H^{s}_{r}(\Omega)$ is the collection of all $f \in L_{\max(r,v)}(\Omega)$ such that 
\begin{align*}
\Vert f \vert L_{r}(\Omega) \Vert + \Big \Vert \Big ( \int_{0}^{1} t^{-2s}   \Big ( t^{-d} \int_{h \in V^{N}(x,t)} \vert ( \Delta^{N}_{h}f )(x) \vert^{v} dh \Big )^{\frac{2}{v}}   \frac{dt}{t} \Big )^{\frac{1}{2}} \Big \vert L_{r}(\Omega) \Big \Vert < \infty
\end{align*}
in the sense of equivalent norms. Here we use the abbreviation \\
 $V^{N}(x,t) := \{ h \in \mathbb{R}^d : |h| < t \; \mbox{and} \; x + \tau h \in \Omega \; \mbox{for} \; 0 \leq \tau \leq N   \} $.
\end{prop}

\begin{proof}
This result is a consequence of Theorem 1.118 in \cite{Tr06} and \eqref{eq_bes_eq_Tri}. Theorem 1.118 is formulated for bounded Lipschitz domains. Of course the unit cube $(0,1)^d$ is a bounded Lipschitz domain. For a definition concerning Lipschitz domains we refer to Stein, see \cite[VI.3.2]{Stein}.
\end{proof}

\subsection{B-Splines, Quarks and Quarklets on the Real Line}\label{sec_inner_quark}

In this paper it is our main goal to describe Bessel-Potential spaces $H^{s}_{r}((0,1)^2)$ in terms of bivariate quarklets. For that purpose we intend to use tensor quarklets resulting out of univariate quarklets. To this end in the following section we recall the definition of univariate quarklets for the shift-invariant setting on $\mathbb{R}$. We follow \cite{DaKRaa}. In a first step we recall the definition of cardinal B-splines. The first order cardinal B-spline $  N_{1}  $ is just the characteristic function of the interval $ [0,1)  $, namely $ N_{1} := \chi_{[0,1)} $. Higher order cardinal B-splines of order $ m \in \mathbb{N}  $ with $ m \geq 2  $ are defined by induction using the convolution $ \ast $. We have 
\begin{align*}
N_{m} := N_{m - 1} \ast N_{1} := \int_{0}^{1} N_{m-1}( \cdot - t ) dt.
\end{align*}
The cardinal B-splines possess some very nice properties. They are collected in the following lemma, see Chapter 5.2 and 5.3 in \cite{DeLo} as well as \cite{Ch1} for the proof.

\begin{lem}\label{Bspl_elem}
Let $ m \in \mathbb{N}   $ and $ x \in \mathbb{R}  $. Then for the cardinal B-splines we observe the following elementary properties. 

\begin{itemize}
\item[(i)] We can write $  N_{m}(x) = \frac{1}{(m-1)!} \sum_{k = 0}^{m}(-1)^{k} { m \choose k } (x - k)_{+}^{m-1}   $.

\item[(ii)] For $ m \geq 2  $ we have $ N_{m}(x) = \frac{x}{m-1} N_{m-1}(x) + \frac{m-x}{m-1} N_{m-1} (x-1)   $.

\item[(iii)] For $ m \geq 3  $ for the derivatives we observe $  N_{m}'(x) = N_{m-1}(x) - N_{m-1}(x-1)  $.

\item[(iv)] The B-splines are compactly supported with $ \supp N_{m} = [0,m]  $.

\item[(v)] We have $ \sum_{k \in \mathbb{Z} } N_{m}(x - k) = 1   $.

\end{itemize}
\end{lem}
When we think on Bessel-Potential spaces it is interesting to know under what conditions on the parameters the cardinal B-splines belong to the spaces $ H^{s}_{r}(\mathbb{R})   $. There is the following useful observation.

\begin{lem}\label{spline_inF}
Let $ s \in \mathbb{R}   $ and $ 1 < r < \infty  $. Let $ m \in \mathbb{N}  $. Then we have
\begin{align*}
N_{m} \in H^{s}_{r}(\mathbb{R}) \qquad \qquad \mbox{if and only if} \qquad \qquad s < m - 1 + \frac{1}{r} .
\end{align*}
\end{lem}

\begin{proof}
This result can be found in \cite{RS}, see Lemma 3 in Chapter 2.3.1. We have to use the relation \eqref{eq_bes_eq_Tri}.
\end{proof}

In what follows for fixed $ m \in \mathbb{N} $ we will work with the symmetrized cardinal B-spline $ \varphi(x)  :=   N_{m}  ( x + \lfloor \frac{m}{2} \rfloor   )    $. We observe $ \supp \varphi = [ - \lfloor \frac{m}{2} \rfloor    ,  \lceil \frac{m}{2} \rceil    ]   $. The symmetrized cardinal B-spline shows up in the following definition of the so-called quarks.

\begin{defi}\label{Bquark}
Let $ m \in \mathbb{N} $ and $ p \in \mathbb{N}_{0}  $. Then the p-th cardinal B-spline quark $ \varphi_{p}  $ is defined by
\begin{equation}
\varphi_{p}(x)  := \Big ( \frac{x}{\lceil \frac{m}{2} \rceil } \Big )^{p}  N_{m} \Big ( x + \lfloor \frac{m}{2} \rfloor  \Big ) .
\end{equation}
\end{defi}
\noindent
The quarks will be very important for us in order to define the quarklets. Their properties have been studied in \cite{DaKRaa}. It is shown in \cite{CoDau} by Cohen, Daubechies and Feauveau that for a given $ \tilde{m} \in \mathbb{N}  $ with $ \tilde{m} \geq m    $ and $  m + \tilde{m} \in 2 \mathbb{N}   $ there exists a compactly supported spline wavelet $ \psi  $ (sometimes called CDF-wavelet) with
\begin{equation}\label{def_CDF_wav}
\psi = \sum_{k \in \mathbb{Z}} b_{k}  \varphi ( 2 \cdot - k   )
\end{equation}
with expansion coefficients $ b_{k} \in \mathbb{R}   $. Only finitely many of them are not zero. Moreover $ \psi  $ has $ \tilde{m}   $ vanishing moments and the  system
\begin{align*}
\Big \{ \varphi (  \cdot - k )  \ : \ k \in \mathbb{Z}  \Big  \} \cup \Big \{ 2^{\frac{j}{2}} \psi (2^{j} \cdot - k) \ : \ j \in \mathbb{N}_{0} \ , \ k \in \mathbb{Z}  \Big \}
\end{align*}
is a Riesz basis for $ L_{2}(\mathbb{R})   $. To construct such a $ \psi   $ we have to work with a compactly supported dual generator $ \tilde{\varphi}   $ associated to the primal generator $ \varphi $ that fulfills
\begin{equation}\label{biorto1}
\left \langle  \varphi , \tilde{\varphi} (\cdot - k)   \right\rangle_{L_{2}(\mathbb{R})} = \delta_{0,k} , \qquad k \in \mathbb{Z}.
\end{equation}
Connected with that there is another compactly supported biorthogonal wavelet $ \tilde{\psi} \in L_{2}(\mathbb{R})  $ with 
\begin{equation}\label{def_biort_wav1}
\tilde{\psi} = \sum_{k \in \mathbb{Z}}  \tilde{b}_{k} \tilde{\varphi} ( 2 \cdot - k   ).
\end{equation}
Here only finitely many of the $ \tilde{b}_{k} \in \mathbb{R}   $  are not zero. Moreover $ \tilde{\psi}  $ has $ m \in \mathbb{N}  $ vanishing moments and the system
\begin{align*}
\Big \{ \tilde{\varphi} (  \cdot - k )  \ : \ k \in \mathbb{Z}  \Big  \} \cup \Big \{ 2^{\frac{j}{2}} \tilde{\psi} (2^{j} \cdot - k) \ : \ j \in \mathbb{N}_{0} \ , \ k \in \mathbb{Z}  \Big \}
\end{align*}
is a Riesz basis for $ L_{2}(\mathbb{R})   $. For $ j \in \mathbb{N}_{0}   $ and $ k \in \mathbb{Z}   $ let us write
\begin{equation}\label{def_biort_wav2}
\psi_{j,k} := 2^{\frac{j}{2}} \psi ( 2^{j} \cdot - k ) \qquad \qquad \mbox{and} \qquad \qquad \tilde{\psi}_{j,k} := 2^{\frac{j}{2}} \tilde{\psi} ( 2^{j} \cdot - k ) .
\end{equation}
Moreover, for $ k \in \mathbb{Z}   $ we put $ \psi_{-1,k} := \varphi (  \cdot - k )   $ and $ \tilde{\psi}_{-1,k} := \tilde{\varphi} (  \cdot - k )   $. Then we observe
\begin{equation}\label{biorto2}
\langle \psi_{j,k} , \tilde{\psi}_{j',k'}     \rangle_{L_{2}(\mathbb{R})} = \delta_{j , j'} \delta_{k , k'} , \qquad j, j' \in \mathbb{N}_{0}, \quad k, k' \in \mathbb{Z} .
\end{equation}
For each $ f \in L_{2}(\mathbb{R})  $ we have
\begin{align}
f & = \sum_{k \in \mathbb{Z}} \langle f , \tilde{\psi}_{-1,k}     \rangle_{L_{2}(\mathbb{R})} \psi_{-1,k}  + \sum_{j \in \mathbb{N}_{0} ,k \in \mathbb{Z}} \langle f , \tilde{\psi}_{j,k}     \rangle_{L_{2}(\mathbb{R})} \psi_{j,k} \nonumber \\
& = \sum_{k \in \mathbb{Z}} \langle f , \psi_{-1,k}     \rangle_{L_{2}(\mathbb{R})} \tilde{\psi}_{-1,k}  +  \sum_{j \in \mathbb{N}_{0}, k \in \mathbb{Z}} \langle f , \psi_{j,k}     \rangle_{L_{2}(\mathbb{R})} \tilde{\psi}_{j,k}    \label{def_biort_wav3}
\end{align}
with convergence in $ L_{2}(\mathbb{R})  $. For details and proofs concerning the above construction we refer to \cite{CoDau}, see especially Section 6.A. Now we can use the CDF-wavelets $ \psi $ to define the quarklets.  
\begin{defi}\label{def_quarklet}
Let $ p \in \mathbb{N}_{0}  $. Then the p-th quarklet $ \psi_{p} $ is defined by 
\begin{equation}
\psi_{p} := \sum_{k \in \mathbb{Z}} b_{k} \varphi_{p}(2 \cdot - k).
\end{equation}
Here the $ b_{k}  $ are the same as in \eqref{def_CDF_wav}. Furthermore, for $ j \in \mathbb{N}_{0}   $ and $ k \in \mathbb{Z}  $ we define
\begin{equation}
\psi_{p,j,k} := 2^{\frac{j}{2}} \psi_{p}(2^{j} \cdot - k) \qquad \qquad \mbox{and} \qquad \qquad \psi_{p,-1,k} := \varphi_{p}( \cdot - k) .
\end{equation}
\end{defi}

\begin{rem}\label{rem_hist_of_qua}
The quarklets given in Definition \ref{def_quarklet} have been introduced around 2017 in \cite{DaKRaa}. Later their properties have been studied in detail in \cite{DaFKRaa}, \cite{DaRaaS} and \cite{HoKoRaVo}. A systematic and comprehensive treatise can be found in \cite{SiDiss}. When defining the quarklets the main focus is on numerical applications. They are specially tailored for adaptive approximation of functions and the numerical treatment of PDEs with very good convergence properties. For that purpose let us refer to \cite{DaHoRaVo} and \cite{DaFKRaa}. For the definition of the quarklets B-spline wavelets are used for the following reasons. B-splines possess optimal smoothness properties compared to their support size. Moreover, explicit formulas exist which make point evaluations quite simple. This issue is important for the construction of suitable quadrature formulas, that are necessary for any numerical scheme for the treatment of PDEs. Of course, in principle quarkonial decompositions also can be provided using other wavelets such as orthonormal Daubechies wavelets. However, Daubechies wavelets are not symmetric which sometimes is disadvantageous. Moreover, for these wavelets no explicit formulas exist which makes point evaluations much more difficult. Using B-spline pre-wavelets would also be a possible choice when constructing quarkonial decompositions. However, the biorthogonal approach we used in Definition \ref{def_quarklet} has the advantage that the lenghts of all filters involved in the associated decomposition and reconstruction schemes are finite, which is usually not the case in the pre-wavelet setting. Of course, due to the polynomial enrichment, the quarklet dictionary is highly redundant. For the construction of adaptive wavelet $hp$-methods, this fact cannot be avoided. At the first glance, this might look as a disadvantage, but it seems to be clear that this is not the case for the following reason. Our long-term goal is the development of adaptive numerical schemes based on quarklets. The art of adaptivity is to find a sparse expansion of an unknown object, namely the solution of a PDE. Now if we work with a very rich dictionary, then the chance to find such a sparse expansion is much higher compared to the basis case where the expansion is unique. From this point of view, redundancy is very helpful. Indeed, in \cite{DaHoRaVo} some numerical experiments showed, that our quarklets can be used for adaptive $hp$-tree approximation of functions with inverse-exponential convergence rates. Moreover, a rigorous proof that certain model singularities showing up in the solution theory of elliptic PDEs can be approximated via quarklets with inverse-exponential rates can be found in \cite{DaRaaS}.  
\end{rem}

\subsection{Quarks and Quarklets on the Interval}\label{Subsec_bound_qua}

When we deal with function spaces defined on bounded intervals such as $I := (0,1) \subset \mathbb{R}$ we require special boundary adapted quarks and quarklets. Their construction is explained in \cite{SiDiss} and \cite{DaFKRaa} and will be summarized in the following section. The foundation of the construction is given by a wavelet basis designed by Primbs, see \cite{Pri1}. In a first step we recall the definition of the so-called Schoenberg B-splines. Let $ m, \tilde{m} \in \mathbb{N}_{0}   $ with $ \tilde{m} \geq m \geq 2   $ and $ m + \tilde{m} \in 2 \mathbb{N}   $. Let $ j_{0} \in \mathbb{N}  $ be a fixed number that depends on $m$ and $\tilde{m}$ and is sufficiently large, see Chapter 4.4 in \cite{Pri1} for further explanations and a precise definition. For $ j \in \mathbb{N} $ with $ j \geq j_{0}  $ let
\begin{equation}\label{index_set_delta0}
\Delta_{j} = \{ -m+1, \ldots , 2^{j} - 1   \} .
\end{equation}
We define the knots 
\[
t^{j}_{k}:= \left\{ \begin{array}{lll}
0    & \quad & \mbox{for} \qquad k = -m+1, \ldots , 0 ;
\\  
2^{-j} k & \quad & \mbox{for}\qquad k = 1 , \ldots , 2^{j} - 1 ;
\\
1 & \quad & \mbox{for}\qquad k = 2^{j} , \ldots , 2^{j} + m - 1 .
\\
\end{array}
\right.
\]
Here the boundary knots have multiplicity $m$. Now the Schoenberg B-splines $ B^{m}_{j,k}  $ are defined by
\begin{equation}\label{eq_def_SchoenbergBs}
B^{m}_{j,k}(x) := ( t^{j}_{k+m} - t^{j}_{k} ) ( \cdot - x )^{m-1}_{+} [ t^{j}_{k}, \ldots , t^{j}_{k+m} ] ,  \qquad k \in \Delta_{j} , x \in I .
\end{equation}
Here the symbol $ ( \cdot - x )^{m-1}_{+} [ t^{j}_{k}, \ldots , t^{j}_{k+m} ]    $ stands for the $m-$th divided difference of the function $  ( \cdot - x )^{m-1}_{+}   $. The generating functions of the Primbs basis are defined by
\begin{equation}\label{eq_qua_gen_jjj}
\varphi_{j,k} := 2^{\frac{j}{2}} B^{m}_{j,k}, \qquad k \in \Delta_{j}.
\end{equation}
The Schoenberg B-splines are generalizations of the cardinal B-splines $ N_{m} $ and have some nice properties. They are collected in the following lemma, see \cite{Pri1} for more details. 

\begin{lem}\label{lem_schoenberg_elementary}
Let $ m \in \mathbb{N}   $, $ j \in \mathbb{N}  $ with $ j \geq j_{0}  $, $ k \in \Delta_{j}  $ and $ x \in I $. Then for the Schoenberg B-splines we observe the following elementary properties. 

\begin{itemize}
\item[(i)] We can write $ B^{m}_{j,k}(x)  = N_{m}(2^{j}x-k)   $ for $ k = 0 , \ldots , 2^{j} - m  $.

\item[(ii)] We have $ B^{m}_{j,k}(x) = B^{m}_{j,2^{j}-m-k}(1-x)   $ for $ k = -m+1 , \ldots , 2^{j} - 1 $.

\item[(iii)] It is $ \supp B^{m}_{j,k} =  [ t^{j}_{k}, t^{j}_{k+m} ]      $.

\item[(iv)] We have $ \sum_{k \in  \Delta_{j}} B^{m}_{j,k}(x)  = 1   $.

\end{itemize}
\end{lem}
\noindent
When dealing with quarklets it is our long-term goal to design adaptive PDE solvers with very good convergence properties that are based on quarklets. Since in most of the cases PDEs come along with some boundary conditions, we have to incorporate these boundary conditions in both the provided quarklet system and the underlying function spaces. For that purpose we have to introduce some additional notation. Let $H^{s}_{r}(I)$ be a Bessel-Potential space with $ 1 < r < \infty   $ and $0 < s < m - 1 + 1/r $. Let $\varepsilon > 0$ be arbitrary small. Then we put
\begin{equation}\label{def_bound_cond_sigma1}
\sigma := ( \sigma^{l} , \sigma^{r} ) \in \{ 0 ,  \ldots  , \lfloor  s + 1 - \frac{1}{r} - \varepsilon  \rfloor   \}^{2}   .
\end{equation} 
Here $  \sigma  $ denotes the order of the boundary conditions. For example $\sigma = ( 0, 1 )$ stands for free boundary conditions at $x=0$ and homogeneous Dirichlet boundary conditions at $x = 1$. We can define Bessel-Potential spaces with incorporated boundary conditions $\sigma$ by
\begin{equation}
(H^{s}_{r}(I))_{\sigma} := \{ f \in    H^{s}_{r}(I) : f(0) = \ldots = f^{(\sigma^{l}-1)}(0) = 0 = f(1) = \ldots = f^{(\sigma^{r}-1)}(1) \} ,
\end{equation}
whereby the equations are meaningful due to Theorem 3.3.1 in \cite{SiTri}. The spaces $ (H^{s}_{r}(I))_{\sigma}   $ can be equipped with the usual norm. That means we use the norm inspired by Definition \ref{def_bespot_dom}, namely
\begin{equation}\label{def_norm_boundcond1}
\Vert f \vert  (H^{s}_{r}(I))_{\sigma}   \Vert := \inf \Big\{\| \, g\, \vert H^{s}_{r}( \mathbb{R} ) \|: \quad g_{|_I} =f \; , \mbox{$g$ fulfills conditions $\sigma$ at 0 and 1}  \Big\} \, .
\end{equation}
With respect to boundary conditions the index set given in \eqref{index_set_delta0} can be generalized by
\begin{equation}
\Delta_{j, \sigma} := \{ -m+1 + \sgn \sigma^{l} , \ldots , 2^{j} - 1 - \sgn \sigma^{r}   \} .
\end{equation}
Here we take into account that in case of active boundary conditions the leftmost and/or the rightmost function in \eqref{eq_qua_gen_jjj} is unused. Recall that the Primbs basis is a biorthogonal wavelet basis. Therefore a dual multiresolution analysis with dual generators $\tilde{\varphi}_{j,k}$ is necessary for the construction. If the generators are represented as column vectors $\Phi_{j} := \{ \varphi_{j,k} : k \in  \Delta_{j, \sigma}  \} $ and $\tilde{\Phi}_{j} := \{ \tilde{\varphi}_{j,k'} : k' \in  \Delta_{j, \sigma}  \} $, they fulfill the duality relation
\begin{align*}
\langle \Phi_{j} , \tilde{\Phi}_{j} \rangle := ( \langle \varphi_{j,k} , \tilde{\varphi}_{j,k'} \rangle_{L_{2}(I)})_{k , k' \in \Delta_{j, \sigma}  } = Id_{|\Delta_{j, \sigma}  |}.
\end{align*}
Notice that we have
\begin{equation}\label{eq_vert_delta_jsig}
|\Delta_{j, \sigma}  | =  2^{j} - 1 - \sgn \sigma^{r} + 1 + m - 1 - \sgn \sigma^{l} =  2^{j} - 1 + m - \sgn \sigma^{r}  - \sgn \sigma^{l} .  
\end{equation}
In the construction of the Primbs basis it is possible to choose different boundary conditions on the primal and the dual side, see Chapter 4.7 in \cite{Pri1}. We will need primal wavelets with free and zero boundary conditions. The corresponding dual wavelets always should be of free boundary type, since this allows us to construct primal wavelets with $\tilde{m}$ vanishing moments. For the construction of the wavelets the following index set is defined:
\[
\nabla_{j, \sigma}:= \left\{ \begin{array}{lll}
\{ 0, 1, \ldots , 2^{j} - 1 \}   & \quad & \mbox{for} \qquad j \geq j_{0} ;
\\  
\Delta_{j , \sigma} & \quad & \mbox{for}\qquad j = j_{0} - 1 .
\\
\end{array}
\right.
\]
To construct the Primbs wavelets we require a stable completion, see \cite{Pri2010} for the details. That means the construction of suitable matrices $M_{j,1}^{\sigma}$, $\tilde{M}_{j,1}^{\sigma}$, that contain the two-scale coefficients of the wavelet column vectors $\Psi = \{ \psi_{j,k}^{\sigma} : k \in  \nabla_{j, \sigma}  \}$. We have
\begin{equation}\label{eq_primbs_wav1}
\Psi_{j} := ( M_{j,1}^{\sigma} )^{T} \Phi_{j+1}, \qquad j \geq j_{0} , 
\end{equation}
with $( M_{j,1}^{\sigma} )^{T} := (b_{k,l}^{j, \sigma})_{k \in \nabla_{j, \sigma}, l \in \Delta_{j+1, \sigma}} \in \mathbb{R}^{|\nabla_{j, \sigma}| \times |\Delta_{j+1, \sigma}|} $. Similar relations hold for $\tilde{\Psi}_{j}$. Then the duality relations $\langle \Psi_{j} , \tilde{\Phi}_{j} \rangle = 0$, $\langle \Phi_{j} , \tilde{\Psi}_{j} \rangle = 0$ and $\langle \Psi_{j} , \tilde{\Psi}_{j} \rangle = Id_{|\nabla_{j, \sigma}|}$ are fulfilled. Now we are able to define the Schoenberg B-spline quarks. 

\begin{defi}\label{def_bound_quark}
Let $ m \in \mathbb{N}   $, $ j \in \mathbb{N}  $ with $ j \geq j_{0}  $ and $  p \in \mathbb{N}_{0}    $. Then the p-th Schoenberg B-spline quark $  \varphi_{p,j,k}  $ is defined by
\[
\varphi_{p,j,k} := \left\{ \begin{array}{lll}
\Big ( \frac{2^{j} \cdot}{k + m}  \Big )^{p} \varphi_{j,k}  & \quad & \mbox{for} \qquad k = -m+1, \ldots , -1 ;
\\  
\Big ( \frac{2^{j} \cdot - k - \lfloor \frac{m}{2} \rfloor }{\lceil \frac{m}{2} \rceil}  \Big )^{p} \varphi_{j,k}  & \quad & \mbox{for}\qquad k = 0 , \ldots , 2^{j} - m ;
\\
\varphi_{p,j,2^{j}-m-k}(1 - \cdot ) & \quad & \mbox{for}\qquad k = 2^{j} - m +1 , \ldots , 2^{j} - 1 .
\\
\end{array}
\right.
\]
\end{defi}
\noindent
Notice that the inner Schoenberg B-spline quarks are translated copies of the cardinal B-spline quarks, see Definition \ref{Bquark}. For $ m \in \mathbb{N}   $ there are $ m - 1   $ left boundary quarks and $ m-1$ right boundary quarks. Now let us turn to the construction of the quarklets. For the inner quarklets we can use a similar method as described in Subsection \ref{sec_inner_quark} above. Let $  p \in \mathbb{N}_{0}   $ and $ j \in \mathbb{N}   $ with $ j \geq j_{0}  $. Let $ k \in  \nabla_{j, \sigma}  $ with $ m - 1 \leq k \leq 2^{j} - m    $. We use the inner wavelets constructed by Primbs, see \cite{Pri1} and \cite{Pri2010}. They can be written as
\begin{equation}\label{wav_prim_bjk}
\psi_{j,k}^{\sigma} := \sum_{l \in \Delta_{j+1, \sigma}} b^{j , \sigma}_{k,l} \varphi_{j+1,l} ,
\end{equation}
see \eqref{eq_primbs_wav1}. Now to construct the inner quarklets for $ m - 1 \leq k \leq 2^{j} - m    $ and $ l \in \Delta_{j+1, \sigma}   $ we use the numbers $ b^{j , \sigma}_{k,l} $ from \eqref{wav_prim_bjk} (see also \eqref{eq_primbs_wav1}) without any changes. That means we put $ b^{p, j , \sigma}_{k,l} := b^{j , \sigma}_{k,l}    $ and define our inner quarklets by
\begin{equation}\label{def_inner_quark}
\psi_{p,j,k}^{\sigma} := \sum_{l \in \Delta_{j+1, \sigma}}  b^{p, j , \sigma}_{k,l} \varphi_{p,j+1,l} .
\end{equation} 
If the inner Primbs wavelets have $ \tilde{m}  $ vanishing moments, also the inner quarklets from \eqref{def_inner_quark} have $ \tilde{m}  $ vanishing moments. This result can be found in \cite{DaKRaa}, see Lemma 2. Next let us construct the boundary quarklets. It will be very important that also they have vanishing moments. Therefore in general for $ p > 0 $ we can not use the boundary wavelets constructed by Primbs, see \cite{Pri1}. A simple counterexample to illustrate that here we maybe do not have the desired vanishing moments can be found in \cite{SiDiss}, see page 67. Hence our strategy to construct the boundary quarklets is to fix that they have $\tilde{m}$ vanishing moments and to obtain a system of linear equations from that determination. Let us deal with the left boundary quarklets. That means we work with $ k = 0, 1, \ldots, m-2    $.  We assume that each left boundary quarklet consists of  $ \tilde{m} + 1 $ quarks, which are either left boundary or inner quarks, see Definition \ref{def_bound_quark}. Furthermore, the $k$-th quarklet representation should begin at the leftmost but $k$-th quark with respect to boundary conditions. This leads to a $ \tilde{m} \times (\tilde{m} + 1)      $ linear system of equations. It can be written as
\begin{equation}\label{eq_bou_quarklet_b}
\sum_{l = -m+1+ \sgn \sigma^{l} + k}^{-m+1+ \sgn \sigma^{l} + k + \tilde{m}} b^{p,j,\sigma}_{k,l} \int_{\mathbb{R}} x^{q} \varphi_{p,j+1,l}(x) dx = 0 , \qquad q = 0, 1, \ldots , \tilde{m} - 1 .
\end{equation} 
That means we have $ \tilde{m} $ linear equations, one for each moment condition. Each of them consists of $ \tilde{m} + 1 $ summands, one for each building quark. The resulting coefficient matrix is of size $ \tilde{m} \times (\tilde{m} + 1)      $ and has a nontrivial kernel, see Chapter 4.3 in \cite{SiDiss}. Consequently we can find nontrivial solutions for \eqref{eq_bou_quarklet_b} and therefore we are able to construct boundary quarklets with vanishing moments. There is the following definition. 

\begin{defi}\label{def_bound_quarklet}
Let $  k = 0, 1, \ldots, m-2   $ and $ \tilde{m} \in \mathbb{N}   $ with $ \tilde{m} \geq m   $. Let $ j \in \mathbb{N}   $ with $ j \geq j_{0}  $ and $ p \in \mathbb{N}_{0}   $. If the vector $ \textbf{b}^{p,j,\sigma}_{k} = (  b^{p,j,\sigma}_{k,-m+1+ \sgn \sigma^{l} + k}  , \ldots , b^{p,j,\sigma}_{k,-m+1+ \sgn \sigma^{l} + k + \tilde{m}}     ) \in \mathbb{R}^{\tilde{m} + 1}    $ with $ \textbf{b}^{p,j,\sigma}_{k} \not = 0   $ is a solution for \eqref{eq_bou_quarklet_b}, then we define the $k$-th left boundary quarklet by 
\begin{equation}\label{eq_kleftboundaryquar}
\psi^{\sigma}_{p,j,k} := \sum_{l = -m+1+ \sgn \sigma^{l} + k}^{-m+1+ \sgn \sigma^{l} + k + \tilde{m}} b^{p,j,\sigma}_{k,l}  \varphi_{p,j+1,l} .
\end{equation}
\end{defi}
\noindent
Here the parameter $k$ refers to the fact that we have $ m - 1  $ left boundary quarklets and $ m - 1  $ right boundary quarklets.
\begin{rem}
The boundary quarklets given in Definition \ref{def_bound_quarklet} have been constructed in \cite{DaFKRaa}, see Section 2.4. In \cite{DaFKRaa} they have been used to obtain frames for $L_{2}(I)$ and the Sobolev spaces $(H^{s}(I))_{\sigma} := (H^{s}_{2}(I))_{\sigma}$ with $ 0 \leq s < m - \frac{1}{2} $, see Theorem 2.7 and Theorem 2.9.
\end{rem}
\noindent
Later it will be convenient for us to use a uniform notation that incorporates both quarks and quarklets simultaneously. For that purpose for $p \in \mathbb{N}_{0}$ and $k \in \nabla_{j_{0}-1 , \sigma}$ we write 
\begin{equation}\label{eq_def_qua_as_qualet}
\psi^{\sigma}_{p,j_{0}-1,k} := \varphi_{p,j_{0},k} .
\end{equation}
Here the functions $\varphi_{p,j_{0},k}$ are the Schoenberg B-spline quarks given in Definition \ref{def_bound_quark}.

\section{Quarklet Characterizations for Bessel-Potential Spaces on Bounded Intervals}\label{sec_bessel_d1}

In this section it is our main goal to deduce quarklet characterizations for univariate Bessel-Potential spaces on bounded intervals. Thereto we use the boundary adapted quarklets introduced in Section \ref{Subsec_bound_qua} above. In connection with that we require some additional notation. We define the index set for the whole quarklet system by
\begin{equation}\label{index_quarkl_full1}
\nabla_{\sigma} :=  \{  (p,j,k) : p,j \in \mathbb{N}_{0}, j \geq j_{0}-1, k \in \nabla_{j, \sigma}     \} .
\end{equation}
It contains the Primbs basis index set, namely
\begin{equation}\label{index_Primbs1_full1}
\nabla_{\sigma}^{P} :=  \{  (0,j,k) : j \in \mathbb{N}_{0}, j \geq j_{0}-1, k \in \nabla_{j, \sigma}     \} .
\end{equation}
Of course we have $ \nabla_{\sigma}^{P} \subset     \nabla_{\sigma}$. Corresponding to $ \nabla_{\sigma}  $ the whole quarklet system itself is given by
\begin{equation}\label{def_quarklsy_ful1}
\Psi_{\sigma} := \{  \psi^{\sigma}_{p,j,k} : (p,j,k) \in   \nabla_{\sigma}    \} .
\end{equation}
Notice that for $j = j_{0}-1$ the system $\Psi_{\sigma}$ contains the Schoenberg B-spline quarks, see Definition \ref{def_bound_quark} and \eqref{eq_def_qua_as_qualet}. For $j \geq j_{0}$ it consists of the boundary adapted quarklets explained in Section \ref{Subsec_bound_qua}. Thereby for $k \in \nabla_{j, \sigma} $ with $ m - 1 \leq k \leq 2^{j} - m    $ it refers to the inner quarklets, see \eqref{def_inner_quark}. Otherwise for $k \in \{ 0, 1, \ldots , m - 2   \}$ or $k \in \{ 2^{j} - m + 1 , \ldots , 2^{j} - 1   \}$ the system $ \Psi_{\sigma} $ consists of left boundary quarklets or right boundary quarklets, respectively. For that we refer to Definition \ref{def_bound_quarklet}.

\subsection{Upper Estimates}

In what follows we prove quarklet characterizations for univariate Bessel-Potential spaces on intervals. Since the proof is rather technical, we subdivide it into several substeps. First of all the subsequent Proposition provides an upper estimate for the norm $ \|\, f \, |H^{s}_{r}(I)\|  $ in terms of a norm describing the quarklet sequence spaces.

\begin{prop}\label{res_Fspq_m-1}
Let $  1 < r < \infty $ and $ m \in \mathbb{N}  $ with $ m \geq 2  $. Let
\begin{equation}\label{res_Fspq_m-1_conds}
0 < s  < m - 1
\end{equation}
and $\delta > 1$. Let arbitrary boundary conditions $ \sigma $ as defined in \eqref{def_bound_cond_sigma1} be given. Let $ f \in L_{r}(I)  $ be, such that there exists a representation 
\begin{equation}\label{eq-inf3}
f = \sum_{(p,j,k) \in \nabla_{\sigma}}   c_{p,j,k} \psi_{p,j,k}^{\sigma}
\end{equation}
with convergence in $ \mathcal{S}'(I)  $, whereby
\begin{equation}\label{eq_later_introduced11}
\Big\| \Big[   \sum_{(p,j,k) \in \nabla_{\sigma}} (p+1)^{4m + 2   \delta }  2^{2 j s} 2^{j}    | c_{p,j,k}|^{2}  | \tilde{\chi}_{j,k}(x) |   \Big]^{\frac{1}{2}}\Big|  L_{r} (I) \Big\| < \infty 
\end{equation}
and $f$ fulfills the boundary conditions given by $\sigma$. Then there exists a constant $ C > 0  $ independent of $ f $, such that 
\begin{align*}
 \|\, f \, |(H^{s}_{r}(I))_{\sigma}\| \leq C \inf_{\eqref{eq-inf3}} \Big\| \Big[   \sum_{(p,j,k) \in \nabla_{\sigma}} (p+1)^{4m + 2   \delta }  2^{2 j s} 2^{j}    | c_{p,j,k}|^{2}  | \tilde{\chi}_{j,k}(x) |   \Big]^{\frac{1}{2}}\Big|  L_{r} (I) \Big\| .
\end{align*}
The infimum is taken over all sequences $ \{ c_{p,j,k}  \}_{(p,j,k) \in \nabla_{\sigma}} \subset \mathbb{C}  $ such that \eqref{eq-inf3} is fulfilled.
\end{prop}

\begin{proof}
Since the proof is rather technical, we split it into several substeps.

\textit{Step 1. Some preparations.} For the proof let $ f \in L_{r}(I)  $ be a function of the form
\begin{equation}\label{m-1_step1eq0}
f = \sum_{p \geq 0} f_{p} =  \sum_{p \geq 0} \sum_{j \geq j_{0} - 1} f_{p,j}  = \sum_{p \geq 0} \sum_{j \geq j_{0} - 1} \sum_{k \in \nabla_{j, \sigma} }   c_{p,j,k} \psi_{p,j,k}^{\sigma} ,
\end{equation}
such that \eqref{eq_later_introduced11} holds and the boundary conditions $ \sigma $ are fulfilled. For such functions we want to prove the estimate
\begin{equation}\label{m-1_step1eq1}
\|\, f \, |(H^s_{r}(I))_{\sigma}\| \lesssim \Big\| \Big[ \sum_{(p,j,k) \in \nabla_{\sigma} }   ( p + 1 )^{4 m + 2 \delta}  2^{2js} 2^{j}    |c_{p,j,k}|^{2}  \tilde{\chi}_{j,k}(x)   \Big]^{\frac{1}{2}}\Big|  L_{r} (I) \Big\| .
\end{equation}
Since $f$ fulfills the boundary conditions $ \sigma  $, we observe  $   \|\, f \, |(H^s_{r}(I))_{\sigma}\| = \|\, f \, |H^s_{r}(I) \| $, see Definition \ref{def_bespot_dom} and formula \eqref{def_norm_boundcond1}. Let $ t \in \mathbb{R} $ be such that $ s < t < m - 1 $. Because of \eqref{res_Fspq_m-1_conds} and $ f \in L_{r}(I)  = L_{\max(r,1)}(I) $ we can use Proposition \ref{satz_diff1} with $v=1$ to get
\begin{equation}\label{m-1_step1eq2}
\|\, f \, |(H^s_{r}(I))_{\sigma}\|   \lesssim \Vert f \vert L_{r}(I) \Vert + \Big \Vert \Big ( \int_{0}^{1} t^{-2s}   \Big ( t^{-1} \int_{h \in V^{N}(x,t)} \vert ( \Delta^{N}_{h}f )(x) \vert dh \Big )^{2}   \frac{dt}{t} \Big )^{\frac{1}{2}} \Big \vert L_{r}(I) \Big \Vert .
\end{equation}
Here we have $ N \in \mathbb{N}  $ with $ m - 1 \geq N > t > s  $, such that $N$ is as small as possible. In what follows we will investigate both terms in \eqref{m-1_step1eq2} separately. 

\textit{Step 2. An estimate concerning the Lebesgue-norm.}
Let us start with the term $ \|\, f \, |L_r (I)\|  $. By using the definition of the function $ f $, see \eqref{m-1_step1eq0}, we observe
\begin{align*}
\|\, f \, |L_r (I)\| & \leq \Big \|\, \sum_{p \geq 0}^{} \sum_{j \geq j_{0}-1} \sum_{k \in \nabla_{j, \sigma}} |c_{p,j,k}|  |\psi_{p,j,k}^{\sigma}| \, \Big |L_r (I) \Big \| .
\end{align*}
Notice that for $j \geq j_{0}$ by \eqref{def_inner_quark} and  \eqref{eq_kleftboundaryquar} it follows that both inner and boundary quarklets can be written in terms of Schoenberg B-spline quarks. For $j = j_{0}-1$ by \eqref{eq_def_qua_as_qualet} we find $  \psi_{p,j_{0}-1,k}^{\sigma} = \varphi_{p,j_{0},k} $ which also leads to the Schoenberg B-spline quarks. Consequently a combination of Definition \ref{def_bound_quark}, equation \eqref{eq_qua_gen_jjj} and Lemma \ref{lem_schoenberg_elementary} yields
\begin{align*}
\|\, f \, |L_r (I)\|  \leq C_{1} \Big \|\, \sum_{p \geq 0}^{} \sum_{j \geq j_{0}-1} \sum_{k \in \nabla_{j, \sigma}} 2^{\frac{j}{2}} |c_{p,j,k}|  | \tilde{\chi}_{j,k}(\cdot) | \, \Big |L_r (I) \Big \| .
\end{align*}
Here $ C_{1}  $ depends on $m$ and $\tilde{m}$, but is independent of $f$. We used that the Schoenberg B-splines are compactly supported and that only finitely many of the numbers $b^{p,j,\sigma}_{k,l} $ are not zero, see \eqref{def_inner_quark} and \eqref{eq_kleftboundaryquar}. Moreover, we utilized the boundedness of all involved functions. Now let $ \delta > 1  $ be fixed. Then the Cauchy-Schwarz inequality yields
\begin{align*}
& \|\, f \, |L_r (I)\| \\
& \qquad  \lesssim  \Big \|\, \sum_{p \geq 0}^{} \sum_{j \geq j_{0}-1} \sum_{k \in \nabla_{j, \sigma}} (p+1)^{\delta} (p+1)^{- \delta} 2^{js} 2^{-js} 2^{\frac{j}{2}} |c_{p,j,k}|  | \tilde{\chi}_{j,k}(\cdot) | \, \Big |L_r(I) \Big \| \\
& \qquad \lesssim \Big \|\,  \Big ( \sum_{p \geq 0}^{} \sum_{j \geq j_{0}-1} (p+1)^{- 2 \delta} 2^{-2js} \Big )^{\frac{1}{2}} \\
& \qquad \qquad \qquad  \times \Big ( \sum_{p \geq 0}^{} \sum_{j \geq j_{0}-1} (p+1)^{ 2 \delta} 2^{2js}  2^{j} \Big [   \sum_{k \in \nabla_{j, \sigma}}  |c_{p,j,k}|  | \tilde{\chi}_{j,k}(\cdot) | \Big ]^{2} \Big )^{\frac{1}{2}} \, \Big |L_r (I) \Big \| \\
& \qquad \lesssim  \Big \|\,  \Big ( \sum_{p \geq 0}^{} \sum_{j \geq j_{0}-1} (p+1)^{ 2 \delta} 2^{2js}  2^{j} \Big [   \sum_{k \in \nabla_{j, \sigma}}  |c_{p,j,k}|  | \tilde{\chi}_{j,k}(\cdot) | \Big ]^{2} \Big )^{\frac{1}{2}}   \, \Big |L_r (I) \Big \| .
\end{align*}
Here the first sum converges since $ s > 0  $ and $ 2 \delta > 1  $. Using the Cauchy-Schwarz inequality again, we find
\begin{align*}
& \|\, f \, |L_r (I)\| \\
& \qquad \lesssim  \Big \|\,  \Big ( \sum_{p \geq 0} \sum_{j \geq j_{0}-1} (p+1)^{2 \delta} 2^{2js}  2^{j} \Big [   \sum_{k \in \nabla_{j, \sigma}}  |c_{p,j,k}|  | \tilde{\chi}_{j,k}(\cdot) |^{\frac{1}{2}} | \tilde{\chi}_{j,k}(\cdot) |^{\frac{1}{2}} \Big ]^{2} \Big )^{\frac{1}{2}} \, \Big |L_r (I) \Big \| \\
& \qquad \lesssim  \Big \|\,  \Big (  \sum_{p \geq 0} \sum_{j \geq j_{0}-1} \sum_{k \in \nabla_{j, \sigma}} (p+1)^{ 2 \delta} 2^{2js}  2^{j}   |c_{p,j,k}|^{2}  | \tilde{\chi}_{j,k}(\cdot) |  \Big (   \sum_{k' \in \nabla_{j, \sigma}}     | \tilde{\chi}_{j,k'}(\cdot) | \Big )  \Big )^{\frac{1}{2}} \, \Big |L_r (I) \Big \| \\
& \qquad \lesssim   \Big \|\, \Big (  \sum_{p \geq 0} \sum_{j \geq j_{0}-1} \sum_{k \in \nabla_{j, \sigma}} (p+1)^{2 \delta}   2^{2js}  2^{j}   |c_{p,j,k}|^{2}  | \tilde{\chi}_{j,k}(\cdot) |  \Big )^{\frac{1}{2}} \, \Big |L_r (I) \Big \| .
\end{align*}
This is what we want to have, see \eqref{m-1_step1eq1}. Hence this step of the proof is complete.

\textit{Step 3. Estimates for differences of functions decomposed into quarklets.}
Now we have to deal with the second term in \eqref{m-1_step1eq2}. Since this issue is rather intricate, we split the calculations into several substeps.
 
\textit{Substep 3.1. Discover truncated quarklet systems.}
At first it is our objective to transform the second term in \eqref{m-1_step1eq2}, such that we can work with truncated quarklet systems. For that purpose we transform the integral concerning $t$ into a sum. We observe
\begin{align*}
& \Big \Vert \Big ( \int_{0}^{1} t^{-2s}   \Big ( t^{-1} \int_{h \in V^{N}(x,t)} \vert ( \Delta^{N}_{h}f )(x) \vert dh \Big )^{2}   \frac{dt}{t} \Big )^{\frac{1}{2}} \Big \vert L_{r}(I) \Big \Vert  \\
& \qquad = \Big \Vert \Big ( \sum_{i = 0}^{\infty} \int_{2^{-i-1}}^{2^{-i}} t^{-2s}   \Big ( t^{-1} \int_{h \in V^{N}(x,t)} \vert ( \Delta^{N}_{h}f )(x) \vert dh \Big )^{2}   \frac{dt}{t} \Big )^{\frac{1}{2}} \Big \vert L_{r}(I) \Big \Vert  \\
& \qquad \lesssim \Big \Vert \Big ( \sum_{i = 0}^{\infty} \int_{2^{-i-1}}^{2^{-i}} 2^{2i(s+1)} 2^{i}   \Big ( \int_{h \in V^{N}(x,t)} \vert ( \Delta^{N}_{h}f )(x) \vert dh \Big )^{2}   dt \Big )^{\frac{1}{2}} \Big \vert L_{r}(I) \Big \Vert  .
\end{align*}
Notice that for $ 2^{-i-1} \leq t \leq 2^{-i}$ we have
\begin{align*}
V^{N}(x,t) & = \{ h \in \mathbb{R} : |h| < t \; \mbox{and} \; x + \tau h \in I \; \mbox{for} \; 0 \leq \tau \leq N   \} \\
& \subseteq \{ h \in \mathbb{R} : |h| < 2^{-i} \; \mbox{and} \; x + \tau h \in I \; \mbox{for} \; 0 \leq \tau \leq N   \} \\
& = V^{N}(x,2^{-i}) .
\end{align*}
Consequently we get
\begin{align*}
& \Big \Vert \Big ( \int_{0}^{1} t^{-2s}   \Big ( t^{-1} \int_{h \in V^{N}(x,t)} \vert ( \Delta^{N}_{h}f )(x) \vert dh \Big )^{2}   \frac{dt}{t} \Big )^{\frac{1}{2}} \Big \vert L_{r}(I) \Big \Vert  \\
& \qquad \lesssim \Big \Vert \Big ( \sum_{i = 0}^{\infty} \int_{2^{-i-1}}^{2^{-i}} 2^{2i(s+1)} 2^{i}   \Big ( \int_{h \in V^{N}(x,2^{-i})} \vert ( \Delta^{N}_{h}f )(x) \vert dh \Big )^{2}   dt \Big )^{\frac{1}{2}} \Big \vert L_{r}(I) \Big \Vert  \\
& \qquad \lesssim \Big \Vert \Big ( \sum_{i = 0}^{\infty}  2^{2i(s+1)}   \Big ( \int_{h \in V^{N}(x,2^{-i})} \vert ( \Delta^{N}_{h}f )(x) \vert dh \Big )^{2}   \Big )^{\frac{1}{2}} \Big \vert L_{r}(I) \Big \Vert  .
\end{align*}
Next we use \eqref{m-1_step1eq0}. Then the Minkowski inequality yields
\begin{align*}
& \Big \Vert \Big ( \int_{0}^{1} t^{-2s}   \Big ( t^{-1} \int_{h \in V^{N}(x,t)} \vert ( \Delta^{N}_{h}f )(x) \vert dh \Big )^{2}   \frac{dt}{t} \Big )^{\frac{1}{2}} \Big \vert L_{r}(I) \Big \Vert  \\
& \qquad \lesssim \Big \Vert \Big ( \sum_{i = 0}^{\infty}  2^{2i(s+1)}   \Big ( \sum_{p \geq 0} \int_{h \in V^{N}(x,2^{-i})} \vert ( \Delta^{N}_{h}f_{p} )(x) \vert dh \Big )^{2}   \Big )^{\frac{1}{2}} \Big \vert L_{r}(I) \Big \Vert  \\
& \qquad \lesssim \Big \Vert \sum_{p \geq 0}  \Big ( \sum_{i = 0}^{\infty}  2^{2i(s+1)}   \Big ( \int_{h \in V^{N}(x,2^{-i})} \vert ( \Delta^{N}_{h}f_{p} )(x) \vert dh \Big )^{2}   \Big )^{\frac{1}{2}} \Big \vert L_{r}(I) \Big \Vert  .
\end{align*}
Now let $ \delta > 1  $ be as before.  Applying the Cauchy-Schwarz inequality we obtain
\begin{align*}
& \Big \Vert \Big ( \int_{0}^{1} t^{-2s}   \Big ( t^{-1} \int_{h \in V^{N}(x,t)} \vert ( \Delta^{N}_{h}f )(x) \vert dh \Big )^{2}   \frac{dt}{t} \Big )^{\frac{1}{2}} \Big \vert L_{r}(I) \Big \Vert  \\
& \quad \lesssim \Big \Vert \sum_{p \geq 0} (p+1)^{- \frac{ \delta }{2}} (p+1)^{ \frac{\delta }{2}} \Big ( \sum_{i = 0}^{\infty}  2^{2i(s+1)}   \Big ( \int_{h \in V^{N}(x,2^{-i})} \vert ( \Delta^{N}_{h}f_{p} )(x) \vert dh \Big )^{2}   \Big )^{\frac{1}{2}} \Big \vert L_{r}(I) \Big \Vert  \\
& \quad \lesssim \Big \Vert \Big ( \sum_{p \geq 0} (p+1)^{- \delta} \Big )^{\frac{1}{2}}    \Big ( \sum_{p \geq 0}   \sum_{i = 0}^{\infty} (p+1)^{\delta }  2^{2i(s+1)}   \Big ( \int_{h \in V^{N}(x,2^{-i})} \vert ( \Delta^{N}_{h}f_{p} )(x) \vert dh \Big )^{2}   \Big )^{\frac{1}{2}} \Big \vert L_{r}(I) \Big \Vert  \\
& \quad \lesssim \Big \Vert    \Big ( \sum_{p \geq 0}   \sum_{i = 0}^{\infty} (p+1)^{\delta }  2^{2i(s+1)}   \Big ( \int_{h \in V^{N}(x,2^{-i})} \vert ( \Delta^{N}_{h}f_{p} )(x) \vert dh \Big )^{2}   \Big )^{\frac{1}{2}} \Big \vert L_{r}(I) \Big \Vert  .
\end{align*}
In the last step we used that due to $ \delta > 1   $ the first series is finite.

\textit{Substep 3.2. Estimates concerning differences.} In what follows we have to estimate the term that contains differences from above by some maximal functions. Let us fix $  i \in \mathbb{N}_{0} $ for a moment. We can apply  \eqref{m-1_step1eq0} to get
\begin{align*}
& \Big ( \int_{h \in V^{N}(x,2^{-i})} \vert ( \Delta^{N}_{h}f_{p} )(x) \vert dh \Big )^{2}  \\
& \quad \lesssim \Big ( \int_{h \in V^{N}(x,2^{-i})} \sum_{j \geq j_{0} - 1}  \vert ( \Delta^{N}_{h}f_{p,j} )(x) \vert dh \Big )^{2}  \\
& \quad \lesssim \Big ( \sum_{j = j_{0} - 1}^{i-1}  \int_{h \in V^{N}(x,2^{-i})}  \vert ( \Delta^{N}_{h}f_{p,j} )(x) \vert dh + \sum_{j = \max ( j_{0} - 1,i)}^{\infty}  \int_{h \in V^{N}(x,2^{-i})}  \vert ( \Delta^{N}_{h}f_{p,j} )(x) \vert dh \Big )^{2}  .
\end{align*}
If $ i - 1 < j_{0} - 1 $ the first sum is empty.   To continue we choose $ \varepsilon > 0  $ such that $ 0 < \varepsilon < \min(s,t-s)   $. Then we observe
\begin{align*}
& \Big ( \int_{h \in V^{N}(x,2^{-i})} \vert ( \Delta^{N}_{h}f_{p} )(x) \vert dh \Big )^{2}  \\
& \qquad \qquad \lesssim \Big ( \sum_{j = j_{0} - 1}^{i-1} \frac{2^{-(i-j)(-\varepsilon )}}{2^{-(i-j)(-\varepsilon )}}  \int_{h \in V^{N}(x,2^{-i})}  \vert ( \Delta^{N}_{h}f_{p,j} )(x) \vert dh \\
& \qquad \qquad \qquad \qquad + \sum_{j = \max ( j_{0} - 1,i)}^{\infty} \frac{2^{(j-i)\varepsilon }}{2^{(j-i)\varepsilon }}  \int_{h \in V^{N}(x,2^{-i})}  \vert ( \Delta^{N}_{h}f_{p,j} )(x) \vert dh \Big )^{2}  .
\end{align*}
For each $  j \in \mathbb{N}_{0} \cup \{ -1 \}   $ we define
\[
\alpha_{j}:= \left\{ \begin{array}{lll}
0  & \quad & \mbox{for} \qquad j < j_{0} - 1 ;
\\  
1 & \quad & \mbox{for}\qquad j \geq j_{0} - 1 .
\\
\end{array}
\right.
\]
Then we can also write
\begin{align*}
& \Big ( \int_{h \in V^{N}(x,2^{-i})} \vert ( \Delta^{N}_{h}f_{p} )(x) \vert dh \Big )^{2}  \\
& \qquad \qquad \lesssim \Big ( \sum_{j = -1}^{i-1} \alpha_{j} \frac{2^{-(i-j)(-\varepsilon )}}{2^{-(i-j)(-\varepsilon )}}  \int_{h \in V^{N}(x,2^{-i})}  \vert ( \Delta^{N}_{h}f_{p,j} )(x) \vert dh \\
& \qquad \qquad \qquad \qquad + \sum_{j = i}^{\infty} \alpha_{j} \frac{2^{(j-i)\varepsilon }}{2^{(j-i)\varepsilon }}  \int_{h \in V^{N}(x,2^{-i})}  \vert ( \Delta^{N}_{h}f_{p,j} )(x) \vert dh \Big )^{2}  .
\end{align*}
There exists a constant $ C_{\varepsilon } $ depending on $ \varepsilon $ (and on $i$, which can be neglected later), such that
\begin{align*}
\sum_{j = -1}^{i-1} 2^{(i-j)(-\varepsilon )} + \sum_{j =  i}^{\infty} 2^{(j-i)(-\varepsilon )}  =  \frac{2^{\varepsilon} - 2^{- \varepsilon ( i + 1)} + 1}{2^{\varepsilon} - 1} := C_{\varepsilon } .
\end{align*}
For $   - 1 \leq j \leq i - 1  $ we put $ \lambda_{j} :=  C_{\varepsilon}^{-1} 2^{(i-j)(-\varepsilon )}  $. For $ j \geq  i   $ we write $ \lambda_{j} :=  C_{\varepsilon}^{-1} 2^{(j-i)(-\varepsilon )} $. Using this notation we obtain
\begin{align*}
& \Big ( \int_{h \in V^{N}(x,2^{-i})} \vert ( \Delta^{N}_{h}f_{p} )(x) \vert dh \Big )^{2}  \\
& \qquad \qquad \lesssim C_{\varepsilon}^{2} \Big ( \sum_{j =  - 1}^{i-1} \alpha_{j} \lambda_{j} 2^{-(i-j)(-\varepsilon )}  \int_{h \in V^{N}(x,2^{-i})}  \vert ( \Delta^{N}_{h}f_{p,j} )(x) \vert dh \\
& \qquad \qquad \qquad \qquad + \sum_{j = i}^{\infty} \alpha_{j} \lambda_{j} 2^{(j-i)\varepsilon }  \int_{h \in V^{N}(x,2^{-i})}  \vert ( \Delta^{N}_{h}f_{p,j} )(x) \vert dh \Big )^{2}  .
\end{align*}
Of course the function $ g(x) = x^{2}  $ is convex. Consequently the Jensen inequality yields
\begin{align*}
& \Big ( \int_{h \in V^{N}(x,2^{-i})} \vert ( \Delta^{N}_{h}f_{p} )(x) \vert dh \Big )^{2}  \\
& \qquad \qquad \lesssim C_{\varepsilon}^{2} \Big ( \sum_{j =  - 1}^{i-1} \alpha_{j}^{2}  2^{-(i-j)(-2 \varepsilon )} \Big ( \int_{h \in V^{N}(x,2^{-i})}  \vert ( \Delta^{N}_{h}f_{p,j} )(x) \vert dh \Big )^{2} \\
& \qquad \qquad \qquad \qquad + \sum_{j = i}^{\infty} \alpha_{j}^{2}  2^{(j-i)2 \varepsilon } \Big (  \int_{h \in V^{N}(x,2^{-i})}  \vert ( \Delta^{N}_{h}f_{p,j} )(x) \vert dh \Big )^{2} \Big )  .
\end{align*}
Here we also used $ \lambda_{j} \leq 1  $. Since $   \alpha_{j}^{2} = \alpha_{j} $ and
\begin{align*}
C_{\varepsilon} =  \frac{2^{\varepsilon} - 2^{- \varepsilon ( i + 1)} + 1}{2^{\varepsilon} - 1} \leq \frac{2^{\varepsilon}  + 1}{2^{\varepsilon} - 1}
\end{align*}
this can be rewritten as 
\begin{align*}
& \Big ( \int_{h \in V^{N}(x,2^{-i})} \vert ( \Delta^{N}_{h}f_{p} )(x) \vert dh \Big )^{2}  \\
& \qquad \qquad \lesssim   \sum_{j = j_{0} - 1}^{i-1}   2^{-(i-j)(-2 \varepsilon )} \Big ( \int_{h \in V^{N}(x,2^{-i})}  \vert ( \Delta^{N}_{h}f_{p,j} )(x) \vert dh \Big )^{2} \\
& \qquad \qquad \qquad \qquad + \sum_{j = \max ( j_{0}-1, i)}^{\infty}   2^{(j-i)2 \varepsilon } \Big (  \int_{h \in V^{N}(x,2^{-i})}  \vert ( \Delta^{N}_{h}f_{p,j} )(x) \vert dh \Big )^{2}  .
\end{align*}
Hence we get
\begin{align*}
& J := \sum_{i = 0}^{\infty}   2^{2i(s+1)} \Big ( \int_{h \in V^{N}(x,2^{-i})} \vert ( \Delta^{N}_{h}f_{p} )(x) \vert dh \Big )^{2}  \\
& \qquad \qquad \lesssim \sum_{j = j_{0} - 1}^{\infty} \sum_{i = j+1}^{\infty}  2^{2i(s+1)}  2^{-(i-j)(-2 \varepsilon )} \Big ( \int_{h \in V^{N}(x,2^{-i})}  \vert ( \Delta^{N}_{h}f_{p,j} )(x) \vert dh \Big )^{2} \\
& \qquad \qquad \qquad \qquad + \sum_{j =  j_{0} - 1}^{\infty} \sum_{i = 0}^{j}  2^{2i(s+1)} 2^{(j-i)2 \varepsilon } \Big (  \int_{h \in V^{N}(x,2^{-i})}  \vert ( \Delta^{N}_{h}f_{p,j} )(x) \vert dh \Big )^{2}  .
\end{align*}
Recall that we defined an additional parameter $t$ with $ m - 1 \geq N > t > s  $. It can be used to obtain
\begin{align*}
J &  \lesssim \sum_{j = j_{0} - 1}^{\infty} \sum_{i = j+1}^{\infty}  2^{2i(s+1)}  2^{-(i-j)(-2 \varepsilon )} 2^{-2i(t+1)} 2^{2i} 2^{2it} \Big ( \int_{h \in V^{N}(x,2^{-i})}  \vert ( \Delta^{N}_{h}f_{p,j} )(x) \vert dh \Big )^{2} \\
& \qquad \qquad  + \sum_{j =  j_{0} - 1}^{\infty} \sum_{i = 0}^{j}  2^{2i(s+1)} 2^{(j-i)2 \varepsilon } \Big (  \int_{h \in V^{N}(x,2^{-i})}  \vert ( \Delta^{N}_{h}f_{p,j} )(x) \vert dh \Big )^{2}  \\
&  \lesssim \sum_{j = j_{0} - 1}^{\infty} \sum_{i = j+1}^{\infty}  2^{2is}  2^{-(i-j)(-2 \varepsilon )} 2^{-2it} \Big [ 2^{2i} \sum_{\ell = i}^{\infty}  2^{2 \ell t} \Big ( \int_{h \in V^{N}(x,2^{- \ell})}  \vert ( \Delta^{N}_{h}f_{p,j} )(x) \vert dh \Big )^{2} \Big ] \\
& \qquad \qquad  + \sum_{j =  j_{0} - 1}^{\infty} \sum_{i = 0}^{j}  2^{2i(s+1)} 2^{(j-i)2 \varepsilon } \Big (  \int_{h \in V^{N}(x,2^{-i})}  \vert ( \Delta^{N}_{h}f_{p,j} )(x) \vert dh \Big )^{2}  .
\end{align*}
Let us look at the case $i \leq j$. Here we find
\begin{align*}
& 2^{2i} \Big (  \int_{h \in V^{N}(x,2^{-i})}  \vert ( \Delta^{N}_{h}f_{p,j} )(x) \vert dh \Big )^{2} \\
& \qquad =  2^{2i} \Big (  \int_{h \in V^{N}(x,2^{-i})} \Big  \vert \sum_{n = 0}^{N} (-1)^{N-n} { N \choose n } f_{p,j} ( x + n h ) \Big \vert dh \Big )^{2} \\
& \qquad \lesssim \sum_{n = 0}^{N}  2^{2i} \Big (  \int_{h \in V^{N}(x,2^{-i})}   \vert  f_{p,j} ( x + n h )  \vert dh \Big )^{2} \\
& \qquad \lesssim 2^{2i} \Big (  \int_{h \in V^{N}(x,2^{-i})}   \vert  f_{p,j} ( x  )  \vert dh \Big )^{2}  + \sum_{n = 1}^{N}   \Big ( \frac{1}{2^{-i}}  \int_{h \in V^{N}(x,2^{-i})}   \vert  f_{p,j} ( x + n h )  \vert dh \Big )^{2} \\
& \qquad \lesssim    \vert  f_{p,j} ( x  )  \vert ^{2}  +  \Big ( \frac{1}{2^{-i}}  \int_{h \in V^{N}(x,2^{-i})}   \vert  f_{p,j} ( x + N h )  \vert dh \Big )^{2} .
\end{align*}
To continue by $E f_{p,j}$ we denote the extension of $f_{p,j}$ by zero, namely 
\[
E f_{p,j}(x):= \left\{ \begin{array}{lll}
f_{p,j}(x)   & \quad & \mbox{for} \qquad x \in I ;
\\  
0 & \quad & \mbox{for}\qquad x \in \mathbb{R} \setminus I .
\\
\end{array}
\right.
\]
Moreover, let $ \textbf{M}   $ be the Hardy-Littlewood-Maximal Operator, see Chapter 1.2.3 in \cite{Tr83} for a definition. Then for $x \in I$ of course we have
\begin{align*}
\vert  f_{p,j} ( x  )  \vert = \vert  E f_{p,j} ( x  )  \vert \lesssim (\textbf{M} \vert  E f_{p,j}  \vert )(x) = (\textbf{M} \vert  E f_{p,j}  \vert )(x) \chi_{I} (x) .
\end{align*}
Furthermore, we observe
\begin{align*}
V^{N}(x,2^{-i}) \subseteq \{ h \in \mathbb{R} : |h| < 2^{-i} \; \mbox{and} \; x + N h \in I   \} .
\end{align*}
Consequently for all $x \in I$ and $h \in V^{N}(x,2^{-i}) $ we have $x + N h \in I$ and therefore $f_{p,j} ( x + N h ) = E f_{p,j} ( x + N h )$. We calculate
\begin{align*}
& \Big ( \frac{1}{2^{-i}}  \int_{h \in V^{N}(x,2^{-i})}   \vert  f_{p,j} ( x + N h )  \vert dh \Big )^{2} \\
& \qquad = \Big ( \frac{1}{2^{-i}}  \int_{h \in V^{N}(x,2^{-i})}   \vert  E f_{p,j} ( x + N h )  \vert dh \Big )^{2} \\
& \qquad \leq \Big ( \frac{1}{2^{-i}}  \int_{- 2^{-i}}^{2^{-i}}   \vert  E f_{p,j} ( x + N h )  \vert dh \Big )^{2} \\
& \qquad \leq  ( \textbf{M} \vert  E f_{p,j} \vert )(x)^{2} .
\end{align*}
All in all for the case $i \leq j$ we find
\begin{equation}\label{eq_prof_i<j}
2^{2i} \Big (  \int_{h \in V^{N}(x,2^{-i})}  \vert ( \Delta^{N}_{h}f_{p,j} )(x) \vert dh \Big )^{2} \lesssim ( \textbf{M} \vert  E f_{p,j} \vert )(x)^{2} . 
\end{equation}
Now we have to deal with the case $i > j$. Here it is important to estimate
\begin{equation}\label{eq_prof_i>jeq11}
2^{2i} \sum_{\ell = i}^{\infty}  2^{2 \ell t} \Big ( \int_{h \in V^{N}(x,2^{- \ell})}  \vert ( \Delta^{N}_{h}f_{p,j} )(x) \vert dh \Big )^{2} .
\end{equation}
Recall that for $p \in \mathbb{N}_{0}$ and $j \geq j_{0} - 1$ we have
\begin{align*}
f_{p,j} = \sum_{k \in \nabla_{j, \sigma} }   c_{p,j,k} \psi_{p,j,k}^{\sigma} ,
\end{align*}
see \eqref{def_inner_quark} and Definition \ref{def_bound_quarklet} for the definition of the quarklets. Hence the functions $f_{p,j} $ are linear combinations of Schoenberg B-splines $B^{m}_{j+1,l} $ that are multiplied with certain polynomials. Thereby the B-spline $B^{m}_{j+1,l} $ is defined via the knots $ t^{j+1}_{l}, \ldots , t^{j+1}_{l+m}$, see \eqref{eq_def_SchoenbergBs}. For all involved B-splines we observe that they are at least $m-2$ times continuously differentiable at all inner knots $2^{-j-1}a$ with $a = 1, \ldots , 2^{j+1}-1$. At the boundary knots $0$ and $1$ some of the involved B-splines are not continuously differentiable. Now let $p \in \mathbb{N}_{0}$, $j \geq j_{0} - 1$ and $\ell \geq i$ be fixed. Let $x \in I$ and $h \in V^{N}(x,2^{- \ell})$. Moreover, recall that $ m - 1 \geq N > t > s  $. Then the mean value theorem yields  
\begin{align*}
\vert  \Delta^{N}_{h}f_{p,j} (x) \vert \lesssim | h |^{N-1} \max_{y \in [x-(N-1)h,x+(N-1)h]} \vert  \Delta^{1}_{h}f_{p,j}^{(N-1)} (y) \vert .
\end{align*} 
Since $ m - 2 \geq N - 1  $ the function $f_{p,j}^{(N-1)}$ is continuous on the given interval. Consequently there exists $y^{*} \in [x-(N-1)h,x+(N-1)h]$ for that the maximum is realized. We can write
\begin{align*}
\vert  \Delta^{N}_{h}f_{p,j} (x) \vert & \lesssim | h |^{N-1} \max_{y \in [x-(N-1)h,x+(N-1)h]} \vert  f_{p,j}^{(N-1)} (y+h) - f_{p,j}^{(N-1)} (y) \vert \\
& = | h |^{N-1}  \vert  f_{p,j}^{(N-1)} (y^{*}+h) - f_{p,j}^{(N-1)} (y^{*}) \vert .
\end{align*} 
Now two different cases are possible. For small $h$ it might be possible that the function $ f_{p,j}^{(N-1)} $ is continuously differentiable between $ y^{*} $ and $ y^{*} + h $. Then we can apply the mean value theorem one more time. Otherwise, due to $j < i \leq \ell$ and $|h| < 2^{- \ell}$ there exists at most one critical point $y_{c}$ between $ y^{*} $ and $ y^{*} + h $ where $ f_{p,j}^{(N-1)} $ is not differentiable. We observe
\begin{align*}
& \vert  f_{p,j}^{(N-1)} (y^{*}+h) - f_{p,j}^{(N-1)} (y^{*}) \vert \\
& \qquad = \vert  f_{p,j}^{(N-1)} (y^{*}+h) - f_{p,j}^{(N-1)}(y_{c}) + f_{p,j}^{(N-1)}(y_{c})  - f_{p,j}^{(N-1)} (y^{*}) \vert \\
& \qquad \leq \vert  f_{p,j}^{(N-1)} (y^{*}+h) - f_{p,j}^{(N-1)}(y_{c}) \vert  + \vert f_{p,j}^{(N-1)}(y_{c})  - f_{p,j}^{(N-1)} (y^{*}) \vert \\
& \qquad \lesssim |h| \vert  f_{p,j}^{(N)} (\varrho_{1})  \vert  + |h| \vert f_{p,j}^{(N)}(\varrho_{2})   \vert .
\end{align*}
Here $ \varrho_{1} \in I  $ is between $ y^{*} + h   $ and $   y_{c}$. The number $ \varrho_{2} \in I  $ is between $ y^{*}    $ and $   y_{c}$. Combining this with our previous estimates we obtain
\begin{align*}
\vert  \Delta^{N}_{h}f_{p,j} (x) \vert \lesssim | h |^{N} ( \vert  f_{p,j}^{(N)} (\varrho_{1})  \vert  +  \vert f_{p,j}^{(N)}(\varrho_{2})   \vert  )  .
\end{align*} 
Recall that $ \varrho_{1}, \varrho_{2} \in (0,1)  $ holds because of $x \in I$ and $h \in V^{N}(x,2^{- \ell})$. Since the functions $f_{p,j}$ are based on quarklets, which are linear combinations of Schoenberg B-splines multiplied with polynomials, in a next step we can apply the Markov inequality as given in \cite{Tim}, see Eq. (37) in Chapter 4.9.6. With $ |h| < 2^{- \ell} \leq 2^{-i} < 2^{-j} $ it yields
\begin{align*}
\vert  \Delta^{N}_{h}f_{p,j} (x) \vert \lesssim (p+1)^{2N} 2^{jN} | h |^{N} \max_{y \in [x- N 2^{-j}, x + N 2^{-j}] \cap I}  \vert  f_{p,j} (y)  \vert   .
\end{align*} 
When we use this to estimate \eqref{eq_prof_i>jeq11} we find
\begin{align*}
& \Big ( \int_{h \in V^{N}(x,2^{- \ell})}  \vert ( \Delta^{N}_{h}f_{p,j} )(x) \vert dh \Big )^{2} \\
& \qquad \lesssim (p+1)^{4N} 2^{2jN} \max_{y \in [x- N 2^{-j}, x + N 2^{-j}] \cap I}  \vert  f_{p,j} (y)  \vert^{2}  \Big ( \int_{|h| <  2^{- \ell}}   | h |^{N}   dh \Big )^{2} \\
& \qquad \lesssim (p+1)^{4N} 2^{2jN} 2^{- 2 N \ell} 2^{- 2 \ell} \max_{y \in [x- N 2^{-j}, x + N 2^{-j}] \cap I}  \vert  f_{p,j} (y)  \vert^{2} .
\end{align*}
Therefore we also get
\begin{align*}
& 2^{2i} \sum_{\ell = i}^{\infty}  2^{2 \ell t} \Big ( \int_{h \in V^{N}(x,2^{- \ell})}  \vert ( \Delta^{N}_{h}f_{p,j} )(x) \vert dh \Big )^{2} \\
& \qquad \lesssim 2^{2i} \sum_{\ell = i}^{\infty}  2^{2 \ell t} (p+1)^{4N} 2^{2jN} 2^{- 2 N \ell} 2^{- 2 \ell} \max_{y \in [x- N 2^{-j}, x + N 2^{-j}] \cap I}  \vert  f_{p,j} (y)  \vert^{2} \\
& \qquad =    (p+1)^{4N} 2^{2i} 2^{2jN}  \max_{y \in [x- N 2^{-j}, x + N 2^{-j}] \cap I}  \vert  f_{p,j} (y)  \vert^{2} \sum_{\ell = i}^{\infty} 2^{2 \ell (t-N)}  2^{- 2 \ell} .
\end{align*}
Since $\ell \geq i$ we observe $ 2^{- 2 \ell} \leq 2^{- 2 i} $. Moreover, because of $ t < N  $ we calculate
\begin{align*}
\sum_{\ell= i}^{\infty}  2^{2 \ell (t-N) } \lesssim \frac{2^{2i(t-N)}}{1 - 2^{2(t-N)}} \lesssim 2^{2i(t-N)}.
\end{align*}
Recall $i>j$. Hence we find $ 2^{2i(t-N)} < 2^{2j(t-N)} $. Consequently we have
\begin{align*}
& 2^{2i} \sum_{\ell = i}^{\infty}  2^{2 \ell t} \Big ( \int_{h \in V^{N}(x,2^{- \ell})}  \vert ( \Delta^{N}_{h}f_{p,j} )(x) \vert dh \Big )^{2} \\
& \qquad \lesssim    (p+1)^{4N} 2^{2jt}   \max_{y \in [x- N 2^{-j}, x + N 2^{-j}] \cap I}  \vert  f_{p,j} (y)  \vert^{2}   .
\end{align*}
Now we use an estimate given in \cite{Tim}, see page 236. It yields
\begin{align*}
& 2^{2i} \sum_{\ell = i}^{\infty}  2^{2 \ell t} \Big ( \int_{h \in V^{N}(x,2^{- \ell})}  \vert ( \Delta^{N}_{h}f_{p,j} )(x) \vert dh \Big )^{2} \\
& \qquad \lesssim    (p+1)^{4N} 2^{2jt}   \max_{y \in [x- N 2^{-j}, x + N 2^{-j}] \cap I}  \vert  E f_{p,j} (y)  \vert^{2}   \\
& \qquad \lesssim    (p+1)^{4(N+1)} 2^{2jt}    \Big ( \frac{1}{2^{-j}}  \int_{- N 2^{-j}}^{N 2^{-j}}  \vert  E f_{p,j} (x + h)  \vert dh  \Big )^{2}   \\
& \qquad \lesssim    (p+1)^{4(N+1)} 2^{2jt}    ( \textbf{M} \vert  E f_{p,j}   \vert )(x)^{2}   .
\end{align*}
To continue we combine the previous estimate with \eqref{eq_prof_i<j}. Then we obtain
\begin{align*}
J &   \lesssim \sum_{j = j_{0} - 1}^{\infty} \sum_{i = j+1}^{\infty}  2^{2is}  2^{-(i-j)(-2 \varepsilon )} 2^{-2it} (p+1)^{4(N+1)} 2^{2jt}    ( \textbf{M} \vert  E f_{p,j}   \vert )(x)^{2} \\
& \qquad \qquad  + \sum_{j =  j_{0} - 1}^{\infty} \sum_{i = 0}^{j}  2^{2is} 2^{(j-i)2 \varepsilon }  ( \textbf{M} \vert  E f_{p,j} \vert )(x)^{2} \\
&   = \sum_{j = j_{0} - 1}^{\infty}    (p+1)^{4(N+1)} 2^{2 j s}   ( \textbf{M} \vert  E f_{p,j}   \vert )(x)^{2}  \sum_{i = j+1}^{\infty} 2^{-2(i-j)(t - \varepsilon - s)}  \\
& \qquad \qquad  + \sum_{j =  j_{0} - 1}^{\infty}    2^{2 j s}   ( \textbf{M} \vert  E f_{p,j} \vert )(x)^{2}  \sum_{i = 0}^{j}  2^{2(j-i)(\varepsilon - s)} .
\end{align*}
Recall $ 0 < \varepsilon < \min(s,t-s)   $. Therefore the series converge and we find
\begin{equation}\label{eq_prof_eqJestcomp}
J   \lesssim \sum_{j = j_{0} - 1}^{\infty}    (p+1)^{4(N+1)} 2^{2 j s}   ( \textbf{M} \vert  E f_{p,j}   \vert )(x)^{2}   .
\end{equation}
This is the desired estimate in terms of maximal functions. 

\textit{Substep 3.3. Use the quarklet decomposition of $f$.}

To complete the proof we have to combine the result of Substep 3.1 with \eqref{eq_prof_eqJestcomp}. Then we find
\begin{align*}
& \Big \Vert \Big ( \int_{0}^{1} t^{-2s}   \Big ( t^{-1} \int_{h \in V^{N}(x,t)} \vert ( \Delta^{N}_{h}f )(x) \vert dh \Big )^{2}   \frac{dt}{t} \Big )^{\frac{1}{2}} \Big \vert L_{r}(I) \Big \Vert  \\
& \qquad \lesssim \Big \Vert    \Big ( \sum_{p \geq 0}     \sum_{j = j_{0} - 1}^{\infty}    (p+1)^{4(N+1) + \delta} 2^{2 j s}   ( \textbf{M} \vert  E f_{p,j}   \vert )(x)^{2}  \Big )^{\frac{1}{2}} \Big \vert L_{r}(I) \Big \Vert  \\
& \qquad = \Big \Vert    \Big ( \sum_{p \geq 0}     \sum_{j = j_{0} - 1}^{\infty}       ( ( \textbf{M} \vert (p+1)^{2(N+1) + \frac{ \delta}{2} }  2^{j s}   E f_{p,j}   \vert )(x))^{2}  \Big )^{\frac{1}{2}} \Big \vert L_{r}(I) \Big \Vert  .
\end{align*}
Since $1 < r < \infty$ we can apply the Fefferman-Stein maximal inequality, see Chapter 1.2.3 in \cite{Tr83}. It yields
\begin{align*}
& \Big \Vert \Big ( \int_{0}^{1} t^{-2s}   \Big ( t^{-1} \int_{h \in V^{N}(x,t)} \vert ( \Delta^{N}_{h}f )(x) \vert dh \Big )^{2}   \frac{dt}{t} \Big )^{\frac{1}{2}} \Big \vert L_{r}(I) \Big \Vert  \\
& \qquad \lesssim \Big \Vert    \Big ( \sum_{p \geq 0}     \sum_{j = j_{0} - 1}^{\infty}       ( ( \textbf{M} \vert (p+1)^{2(N+1) + \frac{ \delta}{2} }  2^{j s}   E f_{p,j}   \vert )(x))^{2}  \Big )^{\frac{1}{2}} \Big \vert L_{r}(\mathbb{R}) \Big \Vert  \\
& \qquad \lesssim \Big \Vert    \Big ( \sum_{p \geq 0}     \sum_{j = j_{0} - 1}^{\infty}         (p+1)^{4(N+1) +  \delta }  2^{2 j s} \vert  E f_{p,j} (x)   \vert^{2}   \Big )^{\frac{1}{2}} \Big \vert L_{r}(\mathbb{R}) \Big \Vert  \\
& \qquad = \Big \Vert  \chi_{I}(x)  \Big ( \sum_{p \geq 0}     \sum_{j = j_{0} - 1}^{\infty}         (p+1)^{4(N+1) +  \delta }  2^{2 j s} \vert   f_{p,j} (x)   \vert^{2}    \Big )^{\frac{1}{2}} \Big \vert L_{r}(\mathbb{R}) \Big \Vert  \\
& \qquad = \Big \Vert    \Big ( \sum_{p \geq 0}     \sum_{j = j_{0} - 1}^{\infty}         (p+1)^{4(N+1) +  \delta }  2^{2 j s} \vert   f_{p,j} (x)   \vert^{2}    \Big )^{\frac{1}{2}} \Big \vert L_{r}(I) \Big \Vert  .
\end{align*}
Now we use the special structure of $f_{p,j}$, see \eqref{m-1_step1eq0}. Then we obtain
\begin{align*}
& \Big \Vert \Big ( \int_{0}^{1} t^{-2s}   \Big ( t^{-1} \int_{h \in V^{N}(x,t)} \vert ( \Delta^{N}_{h}f )(x) \vert dh \Big )^{2}   \frac{dt}{t} \Big )^{\frac{1}{2}} \Big \vert L_{r}(I) \Big \Vert  \\
& \qquad \lesssim  \Big \Vert    \Big ( \sum_{p \geq 0}     \sum_{j = j_{0} - 1}^{\infty}         (p+1)^{4(N+1) +  \delta }  2^{2 j s} \Big \vert     \sum_{k \in \nabla_{j, \sigma} }   c_{p,j,k} \psi_{p,j,k}^{\sigma}(x) \Big \vert^{2}    \Big )^{\frac{1}{2}} \Big \vert L_{r}(I) \Big \Vert  \\
& \qquad \lesssim  \Big \Vert    \Big ( \sum_{p \geq 0}     \sum_{j = j_{0} - 1}^{\infty}         (p+1)^{4(N+1) +  \delta }  2^{2 j s} \Big ( \sum_{k \in \nabla_{j, \sigma} }   | c_{p,j,k}|   | \psi_{p,j,k}^{\sigma}(x)|  \Big )^{2}    \Big )^{\frac{1}{2}} \Big \vert L_{r}(I) \Big \Vert  \\
& \qquad \lesssim  \Big \Vert    \Big ( \sum_{p \geq 0}     \sum_{j = j_{0} - 1}^{\infty}         (p+1)^{4(N+1) +  \delta }  2^{2 j s} 2^{j} \Big ( \sum_{k \in \nabla_{j, \sigma} }   | c_{p,j,k}|  | \tilde{\chi}_{j,k}(x) |  \Big )^{2}    \Big )^{\frac{1}{2}} \Big \vert L_{r}(I) \Big \Vert .
\end{align*}
Next we utilize the Cauchy-Schwarz inequality in combination with $ |\tilde{\chi}_{j,k}(x)| = |\tilde{\chi}_{j,k}(x)|^{\frac{1}{2}}  |\tilde{\chi}_{j,k}(x)|^{\frac{1}{2}}    $. We get
\begin{align*}
& \Big \Vert \Big ( \int_{0}^{1} t^{-2s}   \Big ( t^{-1} \int_{h \in V^{N}(x,t)} \vert ( \Delta^{N}_{h}f )(x) \vert dh \Big )^{2}   \frac{dt}{t} \Big )^{\frac{1}{2}} \Big \vert L_{r}(I) \Big \Vert  \\
& \quad \lesssim  \Big \Vert    \Big ( \sum_{p \geq 0}     \sum_{j = j_{0} - 1}^{\infty}        \sum_{k \in \nabla_{j, \sigma} }  (p+1)^{4(N+1) +  \delta }  2^{2 j s} 2^{j}    | c_{p,j,k}|^{2}  | \tilde{\chi}_{j,k}(x) |  \Big ( \sum_{k' \in \nabla_{j, \sigma} } | \tilde{\chi}_{j,k'}(x) |   \Big )     \Big )^{\frac{1}{2}} \Big \vert L_{r}(I) \Big \Vert \\
& \quad \lesssim  \Big \Vert    \Big ( \sum_{p \geq 0}     \sum_{j = j_{0} - 1}^{\infty}        \sum_{k \in \nabla_{j, \sigma} }  (p+1)^{4m + 2   \delta }  2^{2 j s} 2^{j}    | c_{p,j,k}|^{2}  | \tilde{\chi}_{j,k}(x) |      \Big )^{\frac{1}{2}} \Big \vert L_{r}(I) \Big \Vert .
\end{align*}
Here in the last step we used $N + 1 \leq m$. Since the representation given in \eqref{m-1_step1eq0} was arbitrary, taking the infimum provides the claim. The proof is complete. 
\end{proof}

\begin{rem}
The proof of Proposition \ref{res_Fspq_m-1} employs some ideas given in \cite{HoDa}, see Proposition 12. However, in \cite{HoDa} neither function spaces on domains nor boundary adapted quarklets show up. Therefore in the proof above several technical modifications became necessary. It seems to be possible to extend Proposition \ref{res_Fspq_m-1} to the case of more general Triebel-Lizorkin spaces $ F^{s}_{r,q}(I)  $ using a similar proof technique, as long as suitable conditions concerning the parameters are fulfilled.
\end{rem}

\subsection{Lower Estimates}

Now let us prove the converse estimate. For that purpose we have to deal with the underlying Primbs wavelets and their duals. We will see that they can be interpreted as kernels of local means.

\begin{prop}\label{Hsr_finallow}
Let $  1 < r < \infty $ and $ m \in \mathbb{N}  $ with $ m \geq 2  $. Let $ 0 < s  < m - 1 $ and $\delta > 1$. Let arbitrary boundary conditions $ \sigma $ as defined in \eqref{def_bound_cond_sigma1} be given. Let $ f \in (H^{s}_{r}(I))_{\sigma}  $.  Then there exists a sequence $ \{ c_{p,j,k}  \}_{(p,j,k) \in \nabla_{\sigma}} \subset \mathbb{C}  $, such that $ f $ can be represented as
\begin{equation}\label{rep_qua}
f = \sum_{(p,j,k) \in \nabla_{\sigma}}   c_{p,j,k} \psi_{p,j,k}^{\sigma}
\end{equation}
with convergence in $ \mathcal{S}'(I)  $, whereby we have
\begin{align*}
\Big\| \Big[   \sum_{(p,j,k) \in \nabla_{\sigma}} (p+1)^{4m + 2   \delta }  2^{2 j s} 2^{j}    | c_{p,j,k}|^{2}  | \tilde{\chi}_{j,k}(x) |   \Big]^{\frac{1}{2}}\Big|  L_{r} (I) \Big\| < \infty .
\end{align*}
Moreover, there exists a constant $ C > 0  $ independent of $ f \in (H^{s}_{r}(I))_{\sigma}  $, such that
\begin{align*}
\inf_{\eqref{rep_qua}} \Big\| \Big[   \sum_{(p,j,k) \in \nabla_{\sigma}} (p+1)^{4m + 2   \delta }  2^{2 j s} 2^{j}    | c_{p,j,k}|^{2}  | \tilde{\chi}_{j,k}(x) |   \Big]^{\frac{1}{2}}\Big|  L_{r} (I) \Big\|  \leq C  \|\, f \, |(H^s_{r}(I))_{\sigma}\| .
\end{align*}
Here the infimum is taken over all sequences $ \{ c_{p,j,k}  \}_{(p,j,k) \in \nabla_{\sigma}} \subset \mathbb{C}  $ such that \eqref{rep_qua} is fulfilled.
\end{prop}

\begin{proof}
The main idea for the proof is to work with the quarklets of polynomial degree $p=0$.

\textit{Step 1. Preparations. Quarklets with $p=0$.} 

At first let us collect some properties of the quarklets $\psi_{p,j,k}^{\sigma}$ with polynomial degree $p = 0$. They will be important for the proof later on. Let us start with the inner quarklets, see \eqref{def_inner_quark}. That means we deal with $ j \geq j_{0}  $, $ k \in  \nabla_{j, \sigma}  $ and $ m - 1 \leq k \leq 2^{j} - m    $. For $p = 0$ we observe
\begin{align*}
\psi_{0,j,k}^{\sigma} = \sum_{l \in \Delta_{j+1, \sigma}}  b^{0, j , \sigma}_{k,l} \varphi_{0,j+1,l} = \sum_{l \in \Delta_{j+1, \sigma}}  b^{j , \sigma}_{k,l} \varphi_{j+1,l} .
\end{align*} 
Hence we recover the inner wavelets constructed by Primbs, see \eqref{wav_prim_bjk}. On the other hand we have to deal with the boundary quarklets, see Definition \ref{def_bound_quarklet}. It already has been observed in \cite{SiDiss}, see Theorem 4.23, that for $p=0$ they coincide with the boundary wavelets constructed in \cite{Pri1}. Consequently the quarklets addressed by the index set $\nabla_{\sigma}^{P} $ coincide with the Primbs wavelets, see \cite{Pri1}. It is well-known that the Primbs wavelet system
\begin{equation}\label{eq_Primbs_wavsy}
\Psi_{\sigma}^{P} := \{  \psi^{\sigma}_{0,j,k} : (0,j,k) \in   \nabla_{\sigma}^{P}    \}
\end{equation}
provides a Riesz basis for $L_{2}(I)$, see Theorem 2.12 in \cite{SiDiss}. Moreover, there exists a second wavelet system 
\begin{equation}\label{eq_Primbs_wavsy_dual}
\tilde{\Psi}_{\sigma}^{P} := \{  \tilde{\psi}^{\sigma}_{j,k} : (j,k) \in   \nabla_{\sigma}^{P}    \}
\end{equation}
which is biorthogonal to \eqref{eq_Primbs_wavsy} and also provides a Riesz basis for $L_{2}(I)$, see Theorem 2.10 in \cite{SiDiss}. Details concerning the construction of the dual wavelet basis can be found in \cite{Pri1}, see Chapter 4. The dual wavelets have the following properties that are important for us later on:
\begin{itemize}
\item[(i)] There exists a constant $ C > 0 $ such that $\supp \tilde{\psi}^{\sigma}_{j,k} \subset ( C I_{j,k}) \cap I \subset C I_{j,k}  $ holds for all $(j,k) \in \nabla_{\sigma}^{P}$. This follows from the construction summarized in Chapter 2.2  in \cite{SiDiss} in combination with Theorem 2.10 from there. Much more explanations and explicit constants $C$ depending on $m$ and $\tilde{m}$ can be found in \cite{Pri1}, see Chapter 4.4, Remark 5.13 and Theorem 5.16. 

\item[(ii)] For all $\tilde{\psi}^{\sigma}_{j,k}$ with $(j,k) \in \nabla_{\sigma}^{P}$ it holds $ \tilde{\psi}^{\sigma}_{j,k} \in L_{\infty}(I) $ and $\Vert \tilde{\psi}^{\sigma}_{j,k} \vert L_{\infty}(I) \Vert \lesssim 2^{\frac{j}{2}} $. For the inner functions this already has been observed in \cite{CoDau}, see Section 6.A. For the boundary functions this is a consequence of the construction given in \cite{Pri1}, see Chapters 4.1 - 4.4 and Chapter 5.1. 

\item[(iii)] The functions $\tilde{\psi}^{\sigma}_{j,k} $ with $(j,k) \in  \{  (j,k) : j \in \mathbb{N}, j \geq j_{0}, k \in \nabla_{j, \sigma}     \}$ have $m$ vanishing moments, see Theorem 2.12 in \cite{SiDiss}.
\end{itemize}
Notice that the dual wavelets satisfy all conditions given in Definition 1.9 in \cite{Tr10}, whereby we put $A = 0$ and $B = m$. This will be important for us later.

\textit{Step 2. Interpret the dual Primbs wavelets as kernels of local means.} 

To continue for $(j,k) \in \{ (j,k) : j \in \mathbb{N}_{0} , k \in \mathbb{Z}  \} \setminus \nabla_{\sigma}^{P}   $ we put $\tilde{\psi}^{\sigma}_{j,k} = 0$. Of course the zero function fulfills the properties (i)-(iii) listed in Step 1. Consequently the system $\{ 2^{\frac{j}{2}} \tilde{\psi}^{\sigma}_{j,k} : j \in \mathbb{N}_{0} , k \in \mathbb{Z}   \}$ can be seen as a system of kernels of local means consistend with Definition 1.9 in \cite{Tr10}. Now for $ f \in (H^{s}_{r}(I))_{\sigma}  $ we set
\begin{equation}\label{eq_localmean_zeroext}
\lambda_{j,k} := \int_{I}  \tilde{\psi}^{\sigma}_{j,k}(x) f(x) dx = \int_{- \infty}^{\infty}  \tilde{\psi}^{\sigma}_{j,k}(x) f(x) dx .
\end{equation}
Since $s > 0$ and $1 < r < \infty$ we observe that $ f \in (H^{s}_{r}(I))_{\sigma}  $ yields $f \in L_{1}(I)$. Therefore the numbers $ 2^{\frac{j}{2}} \lambda_{j,k}$ are well-defined. Moreover, since $s < m-1 < m$ they can be interpreted as local means according to Definition 1.13 and formula (1.52) in \cite{Tr10}. Hence we can apply Theorem 1.15 (ii) in \cite{Tr10}. For that purpose let $g \in H^{s}_{r}(\mathbb{R})$ with $f = g$ in $I$, such that $\Vert g \vert H^{s}_{r}(\mathbb{R}) \Vert \lesssim 2 \Vert f \vert (H^{s}_{r}(I))_{\sigma} \Vert $. Since $ f \in (H^{s}_{r}(I))_{\sigma}  $ we have $ \Vert f \vert (H^{s}_{r}(I))_{\sigma} \Vert =   \Vert f \vert H^{s}_{r}(I) \Vert $. We find
\begin{align*}
\Big\| \Big (   \sum_{j \in \mathbb{N}_{0}, k \in \mathbb{Z}}   2^{2 j s} 2^{j}    | \lambda_{j,k}|^{2}  | \chi_{j,k}(x) |   \Big )^{\frac{1}{2}}\Big|  L_{r} (\mathbb{R}) \Big\| & \lesssim \Vert g \vert H^{s}_{r}(\mathbb{R}) \Vert  \lesssim 2 \Vert f \vert (H^{s}_{r}(I))_{\sigma} \Vert .
\end{align*}
Due to \eqref{eq_localmean_zeroext} we observe that for $(j,k) \in \{ (j,k) : j \in \mathbb{N}_{0} , k \in \mathbb{Z}  \} \setminus \nabla_{\sigma}^{P}   $ we have $\lambda_{j,k} =0$. Consequently we also can write 
\begin{align*}
\Big\| \Big (   \sum_{(j,k) \in \nabla_{\sigma}^{P}  }   2^{2 j s} 2^{j}    | \lambda_{j,k}|^{2}  | \chi_{j,k}(x) |   \Big )^{\frac{1}{2}}\Big|  L_{r} (\mathbb{R}) \Big\| & \lesssim 2 \Vert f \vert (H^{s}_{r}(I))_{\sigma} \Vert .
\end{align*}
Due to the definition of the norm $ \Vert \cdot \vert  L_{r}(I) \Vert    $ via an infimum for $ f \in (H^{s}_{r}(I))_{\sigma}  $ we obtain
\begin{equation}\label{eq_quadual_locmean_eq1}
\Big\| \Big (   \sum_{(j,k) \in \nabla_{\sigma}^{P}  }   2^{2 j s} 2^{j}    | \lambda_{j,k}|^{2}  | \tilde{\chi}_{j,k}(x) |   \Big )^{\frac{1}{2}}\Big|  L_{r} (I) \Big\|  \lesssim 2 \Vert f \vert (H^{s}_{r}(I))_{\sigma} \Vert .
\end{equation}

\textit{Step 3. The function $f$ can be represented via Primbs wavelets.} 

Now we prove that $ f \in (H^{s}_{r}(I))_{\sigma}  $ can be represented as
\begin{equation}\label{rep_qua_Primbsw}
f = \sum_{(j,k) \in \nabla_{\sigma}^{P}}   \lambda_{j,k} \psi_{0,j,k}^{\sigma}
\end{equation}
with convergence in $ \mathcal{S}'(I)  $. That means we are interested in the function
\begin{align*}
h = \sum_{(j,k) \in \nabla_{\sigma}^{P}}   \lambda_{j,k} \psi_{0,j,k}^{\sigma} .
\end{align*}
A combination of Proposition \ref{res_Fspq_m-1}, formula \eqref{eq_quadual_locmean_eq1} and $f \in (H^{s}_{r}(I))_{\sigma}$ yields $h \in (H^{s}_{r}(I))_{\sigma}$. We have to show that $f = h$ in the sense of $\mathcal{S}'(I) $. For that purpose we prove that for all $ \eta \in \mathcal{S}(I)  $ we have $ \langle f - h ,  \eta \rangle_{L_{2}(I)} = 0 $. We use the biorthogonality of $\Psi_{\sigma}^{P} $ and $\tilde{\Psi}_{\sigma}^{P} $, see Theorem 5.7 in \cite{Pri1}. It yields that  $ \eta \in \mathcal{S}(I)  \subset L_{2}(I)   $ can be written in the form
\begin{align*}
\eta =   \sum_{(j,k) \in \nabla_{\sigma}^{P}} \left\langle \eta , \psi_{j,k}^{\sigma}     \right\rangle_{L_{2}(I)} \tilde{\psi}_{j,k}^{\sigma}
\end{align*}
with convergence in $ L_{2}(I)  $, see Lemma 2.6 in \cite{Pri1}. Convergence in $ L_{2}(I)  $ yields convergence pointwise almost everywhere for an appropriate subsequence. Now let $ (j' , k') \in \nabla_{\sigma}^{P}  $ be fixed. Then we find
\begin{align*}
& \int_{I} (f(x) - h(x)) \langle \eta , \psi_{j',k'}^{\sigma}     \rangle_{L_{2}(I)} \tilde{\psi}_{j',k'}^{\sigma}(x) dx \\
& \qquad =  \int_{I} \Big (f(x) -  \sum_{(j,k) \in \nabla_{\sigma}^{P}} \lambda_{j,k} \psi_{j,k}^{\sigma}(x) \Big) \langle \eta , \psi_{j',k'}^{\sigma}     \rangle_{L_{2}(I)} \tilde{\psi}_{j',k'}^{\sigma}(x) dx \\
& \qquad = \langle \eta , \psi_{j',k'}^{\sigma}     \rangle_{L_{2}(I)}    \langle f , \tilde{\psi}_{j',k'}^{\sigma} \rangle_{L_{2}(I)}   - \langle \eta , \psi_{j',k'}^{\sigma}     \rangle_{L_{2}(I)}    \sum_{(j,k) \in \nabla_{\sigma}^{P}} \lambda_{j,k}  \int_{I} \psi_{j,k}^{\sigma}(x) \tilde{\psi}_{j',k'}^{\sigma}(x)    dx .
\end{align*}
We apply the definition of the numbers $ \lambda_{j,k}  $, see \eqref{eq_localmean_zeroext}, and biorthogonality. We get
\begin{align*}
& \int_{I} (f(x) - h(x)) \langle \eta , \psi_{j',k'}^{\sigma}     \rangle_{L_{2}(I)} \tilde{\psi}_{j',k'}^{\sigma}(x) dx \\
& \qquad = \langle \eta , \psi_{j',k'}^{\sigma}     \rangle_{L_{2}(I)}    \langle f , \tilde{\psi}_{j',k'}^{\sigma} \rangle_{L_{2}(I)}   - \langle \eta , \psi_{j',k'}^{\sigma}     \rangle_{L_{2}(I)} \langle f , \tilde{\psi}_{j',k'}^{\sigma} \rangle_{L_{2}(I)} = 0 .  
\end{align*}
A similar calculation can be done for every linear combination of $  \langle \eta , \psi_{j,k}^{\sigma}     \rangle_{L_{2}(I)} \tilde{\psi}_{j,k}^{\sigma} $. Consequently we obtain $ f = h   $ in the sense of $ \mathcal{S}'(I)  $. This step of the proof is complete. 

\textit{Step 4. Complete the proof.}

To complete the proof let $ f \in (H^{s}_{r}(I))_{\sigma}  $. Since $\Psi_{\sigma}^{P} \subset \Psi_{\sigma} $ Step 3 yields that there exists a sequence $ \{ c_{p,j,k}  \}_{(p,j,k) \in \nabla_{\sigma}} \subset \mathbb{C}  $, such that $ f $ can be represented as
\begin{align*}
f = \sum_{(p,j,k) \in \nabla_{\sigma}}   c_{p,j,k} \psi_{p,j,k}^{\sigma}
\end{align*}
with convergence in $ \mathcal{S}'(I)  $. Moreover, we observe
\begin{align*}
& \inf_{\eqref{rep_qua}} \Big\| \Big[   \sum_{(p,j,k) \in \nabla_{\sigma}} (p+1)^{4m + 2   \delta }  2^{2 j s} 2^{j}    | c_{p,j,k}|^{2}  | \tilde{\chi}_{j,k}(x) |   \Big]^{\frac{1}{2}}\Big|  L_{r} (I) \Big\| \\
& \qquad \leq \inf_{\eqref{rep_qua_p=0}} \Big\| \Big[   \sum_{(j,k) \in \nabla_{\sigma}^{P}}   2^{2 j s} 2^{j}    | \tilde{c}_{j,k}|^{2}  | \tilde{\chi}_{j,k}(x) |   \Big]^{\frac{1}{2}}\Big|  L_{r} (I) \Big\| , 
\end{align*}
at which at the right hand side the infimum is taken over all representations
\begin{equation}\label{rep_qua_p=0}
f = \sum_{(j,k) \in \nabla_{\sigma}^{P}}   \tilde{c}_{j,k} \psi_{0,j,k}^{\sigma} .
\end{equation} 
Here we used that for $p=0$ we have $1^{4m + 2   \delta } = 1$. To estimate the infimum from above we can use the representation obtained in \eqref{rep_qua_Primbsw}. Then we find
\begin{align*}
& \inf_{\eqref{rep_qua}} \Big\| \Big[   \sum_{(p,j,k) \in \nabla_{\sigma}} (p+1)^{4m + 2   \delta }  2^{2 j s} 2^{j}    | c_{p,j,k}|^{2}  | \tilde{\chi}_{j,k}(x) |   \Big]^{\frac{1}{2}}\Big|  L_{r} (I) \Big\| \\
& \qquad \leq  \Big\| \Big[   \sum_{(j,k) \in \nabla_{\sigma}^{P}}   2^{2 j s} 2^{j}    | \lambda_{j,k}|^{2}  | \tilde{\chi}_{j,k}(x) |   \Big]^{\frac{1}{2}}\Big|  L_{r} (I) \Big\| .
\end{align*}
An application of \eqref{eq_quadual_locmean_eq1} yields
\begin{align*}
\inf_{\eqref{rep_qua}} \Big\| \Big[   \sum_{(p,j,k) \in \nabla_{\sigma}} (p+1)^{4m + 2   \delta }  2^{2 j s} 2^{j}    | c_{p,j,k}|^{2}  | \tilde{\chi}_{j,k}(x) |   \Big]^{\frac{1}{2}}\Big|  L_{r} (I) \Big\| \lesssim  \Vert f \vert (H^{s}_{r}(I))_{\sigma} \Vert .
\end{align*}
The proof is complete.
\end{proof}

Now we can combine our previous findings to deduce a characterization of the univariate Bessel-Potential spaces $ (H^{s}_{r}(I))_{\sigma} $ in terms of quarklets. It is the first main result of this paper.  

\begin{satz}\label{mainresult1_Hsr_d=1}
Let $  1 < r < \infty $ and $ m \in \mathbb{N}  $ with $ m \geq 2  $. Let $ 0 < s  < m - 1 $ and $\delta > 1$. Let arbitrary boundary conditions $ \sigma $ as defined in \eqref{def_bound_cond_sigma1} be given. Let $ f \in L_{r}(I)   $. Then we have   $ f \in (H^{s}_{r}(I))_{\sigma}  $ if and only if $ f $ can be represented as
\begin{equation}\label{rep_main_Hsr_d=1}
f = \sum_{(p,j,k) \in \nabla_{\sigma}}   c_{p,j,k} \psi_{p,j,k}^{\sigma}
\end{equation}
with convergence in $ \mathcal{S}'(I)  $, where we have that
\begin{equation}\label{rep_main_seqnorm_Hsr_d=1}
\Big\| \Big[   \sum_{(p,j,k) \in \nabla_{\sigma}} (p+1)^{4m + 2   \delta }  2^{2 j s} 2^{j}    | c_{p,j,k}|^{2}  | \tilde{\chi}_{j,k}(x) |   \Big]^{\frac{1}{2}}\Big|  L_{r} (I) \Big\|
\end{equation}
is finite. Moreover the norms $ \Vert f \vert (H^{s}_{r}(I))_{\sigma}  \Vert    $ and
\begin{align*}
\Vert f \vert (H^{s}_{r}(I))_{\sigma}  \Vert^{(\delta, \nabla_{\sigma})}  := \inf_{\eqref{rep_main_Hsr_d=1}} \Big\| \Big[   \sum_{(p,j,k) \in \nabla_{\sigma}} (p+1)^{4m + 2   \delta }  2^{2 j s} 2^{j}    | c_{p,j,k}|^{2}  | \tilde{\chi}_{j,k}(x) |   \Big]^{\frac{1}{2}}\Big|  L_{r} (I) \Big\|
\end{align*}
are equivalent. Here the infimum is taken over all sequences $ \{ c_{p,j,k}  \}_{(p,j,k) \in \nabla_{\sigma}}   $ such that \eqref{rep_main_Hsr_d=1} is fulfilled.
\end{satz}

\begin{proof}
This result is a direct consequence of Proposition \ref{res_Fspq_m-1} and Proposition \ref{Hsr_finallow}.
\end{proof}

\begin{rem}
We observe that in Theorem \ref{mainresult1_Hsr_d=1} the condition $ s  < m - 1    $ shows up. On the other hand Lemma \ref{spline_inF} implies that the assertion of Theorem \ref{mainresult1_Hsr_d=1} becomes wrong for $  s \geq m - 1 + \frac{1}{r}   $. Most likely it is possible to prove a counterpart of Theorem \ref{mainresult1_Hsr_d=1} for $ m - 1 \leq s <  m - 1 + \frac{1}{r}   $. However, for that purpose another method for the proof must be applied. Some attempts to close this gap in the shift-invariant setting have been made in \cite{HoDa} using complex interpolation. Though a final answer is not given there. Further discussions concerning the necessity of conditions on the parameters also can be found in \cite{GSU}, \cite{SU17} and \cite{Tr78Ser}.
\end{rem}

Consistent with Theorem \ref{mainresult1_Hsr_d=1} we also can prove a quarklet characterization for the Lebesgue spaces $ L_{r}(I)  $.

\begin{cor}\label{mainresult1_Lr_d=1}
Let $  1 < r < \infty $ and $ m \in \mathbb{N}  $ with $ m \geq 2  $. Let $ \tilde{m} \in \mathbb{N}_{0}   $ with $ m + \tilde{m} \in 2 \mathbb{N}   $ and $  \tilde{m} > 5 m + 12 $. Let $\delta > 1$ and $\sigma = (0,0)$. Let $ f \in \mathcal{S}'(I)   $. Then we have   $ f \in L_{r}(I)  $ if and only if $ f $ can be represented via \eqref{rep_main_Hsr_d=1} with convergence in $ \mathcal{S}'(I)  $, where we have that
\begin{equation}\label{rep_main_seqnorm_Lr_d=1}
\Big\| \Big[   \sum_{(p,j,k) \in \nabla_{\sigma}} (p+1)^{2 \delta }   2^{j}    | c_{p,j,k}|^{2}  | \tilde{\chi}_{j,k}(x) |   \Big]^{\frac{1}{2}}\Big|  L_{r} (I) \Big\|
\end{equation}
is finite. Moreover, the norms $ \Vert f \vert L_{r}(I)  \Vert    $ and
\begin{align*}
\Vert f \vert L_{r}(I)  \Vert^{(\delta, \nabla_{\sigma})}  := \inf_{\eqref{rep_main_Hsr_d=1}} \Big\| \Big[   \sum_{(p,j,k) \in \nabla_{\sigma}} (p+1)^{ 2   \delta }  2^{j}    | c_{p,j,k}|^{2}  | \tilde{\chi}_{j,k}(x) |   \Big]^{\frac{1}{2}}\Big|  L_{r} (I) \Big\|
\end{align*}
are equivalent. Here the infimum is taken over all sequences $ \{ c_{p,j,k}  \}_{(p,j,k) \in \nabla_{\sigma}}   $ such that \eqref{rep_main_Hsr_d=1} is fulfilled.
\end{cor}

\begin{proof}
\textit{Step 1.} At first we prove a counterpart of Proposition \ref{res_Fspq_m-1}. For that purpose let $f$ be a function as defined in \eqref{m-1_step1eq0}. We follow the proof of Proposition \ref{res_Fspq_m-1}, see Step 2. Using the Cauchy-Schwarz inequality and $\delta > 1$ we find
\begin{align*}
\|\, f \, |L_r (I)\| &   \lesssim  \Big \|\, \sum_{p \geq 0}^{} \sum_{j \geq j_{0}-1} \sum_{k \in \nabla_{j, \sigma}} (p+1)^{\delta} (p+1)^{- \delta}   2^{\frac{j}{2}} |c_{p,j,k}|  | \tilde{\chi}_{j,k}(\cdot) | \, \Big |L_r(I) \Big \| \\
&  \lesssim  \Big \|\,  \Big ( \sum_{p \geq 0}^{} \sum_{j \geq j_{0}-1} (p+1)^{ 2 \delta}  2^{j} \Big [   \sum_{k \in \nabla_{j, \sigma}}  |c_{p,j,k}|  | \tilde{\chi}_{j,k}(\cdot) | \Big ]^{2} \Big )^{\frac{1}{2}}   \, \Big |L_r (I) \Big \| \\
&  \lesssim  \Big \|\,  \Big (  \sum_{p \geq 0} \sum_{j \geq j_{0}-1} \sum_{k \in \nabla_{j, \sigma}} (p+1)^{ 2 \delta}   2^{j}   |c_{p,j,k}|^{2}  | \tilde{\chi}_{j,k}(\cdot) |  \Big (   \sum_{k' \in \nabla_{j, \sigma}}     | \tilde{\chi}_{j,k'}(\cdot) | \Big )  \Big )^{\frac{1}{2}} \, \Big |L_r (I) \Big \| \\
&  \lesssim   \Big \|\, \Big (  \sum_{p \geq 0} \sum_{j \geq j_{0}-1} \sum_{k \in \nabla_{j, \sigma}} (p+1)^{2 \delta}   2^{j}   |c_{p,j,k}|^{2}  | \tilde{\chi}_{j,k}(\cdot) |  \Big )^{\frac{1}{2}} \, \Big |L_r (I) \Big \| .
\end{align*}
So Step 1 of the proof is complete.

\textit{Step 2.} Now we prove a counterpart of Proposition \ref{Hsr_finallow}. As before it will be enough to work with the underlying Primbs wavelet basis $ \Psi_{\sigma}^{P}  $. We are interested in a counterpart of formula \eqref{eq_quadual_locmean_eq1}. For that purpose again we apply Theorem 1.15 (ii) in \cite{Tr10}. This is possible since for $1 < r < \infty$ we have $L_{r}(\mathbb{R}) = F^{0}_{r,2}(\mathbb{R})$. However, when we use Theorem 1.15 (ii) with $s = 0$ the additional condition $A > 0$ shows up, see also Definition 1.9 in \cite{Tr10}. That means we require some additional smoothness for the dual wavelets $\tilde{\psi}^{\sigma}_{j,k}$. Put $A = 1$. We already have seen in the proof of Proposition \ref{Hsr_finallow}, see Step 1, that  for all $\tilde{\psi}^{\sigma}_{j,k}$ with $(j,k) \in \nabla_{\sigma}^{P}$ it holds $ \tilde{\psi}^{\sigma}_{j,k} \in L_{\infty}(I) $ and $\Vert \tilde{\psi}^{\sigma}_{j,k} \vert L_{\infty}(I) \Vert \lesssim 2^{\frac{j}{2}} $. It follows by the construction of Primbs, that the dual boundary functions have the same Sobolev regularity as the dual inner scaling functions given in \cite{CoDau}, see Remark 4.10 and Remark 4.18 in \cite{Pri1}. It has been observed in \cite{CoDau}, see formula (6.17), that the dual inner scaling function $\tilde{\varphi}$ belongs to $C^{1}(\mathbb{R})$ (and therefore is also Lipschitz continuous), if
\begin{align*}
\tilde{m} > 5 m + 12 .
\end{align*}
Consequently in view of Step 1 of the proof of Proposition \ref{Hsr_finallow} we are now able to apply Theorem 1.15 (ii) in \cite{Tr10} with $s = 0$. Using some ideas of Step 2 of the proof of Proposition \ref{Hsr_finallow} for $f \in L_{r}(I)$ we obtain
\begin{align*}
\Big\| \Big (   \sum_{(j,k) \in \nabla_{\sigma}^{P}  }   2^{j}    | \lambda_{j,k}|^{2}  | \tilde{\chi}_{j,k}(x) |   \Big )^{\frac{1}{2}}\Big|  L_{r} (I) \Big\|  \lesssim  \Vert f \vert L_{r}(I) \Vert .
\end{align*}
Here the numbers $ \lambda_{j,k} $ are the same as in \eqref{eq_localmean_zeroext}. They are well-defined since $f \in L_{r}(I) \subset L_{1}(I)$. Now following Step 3 of the proof of Proposition \ref{Hsr_finallow} we can show that each $ f \in L_{r}(I)  $ can be represented via \eqref{rep_qua_Primbsw} with convergence in $ \mathcal{S}'(I)  $. Finally to complete the proof we can proceed as described in Step 4 of the proof of Proposition \ref{Hsr_finallow}.
\end{proof}

\begin{rem}
Theorem \ref{mainresult1_Hsr_d=1} and Corollary \ref{mainresult1_Lr_d=1} can be interpreted as extensions of Theorem 2.7 and Theorem 2.9 in \cite{DaFKRaa}. Here corresponding results for the Hilbert space case $r = 2$ have been formulated. However, in   \cite{DaFKRaa} the proofs are not given in detail. Hence our calculations above also can be seen as a subsequent verification of some results given in \cite{DaFKRaa}.
\end{rem}

\section{Bivariate Bessel-Potential Spaces via Tensor Products}\label{sec_tensor}

It is the main goal of this paper to characterize bivariate Bessel-Potential spaces defined on squares via bivariate quarklets that are constructed out of univariate quarklets using tensor product methods. In \cite{DaFKRaa} multivariate quarklets already have been used to describe the Sobolev spaces $ H^{s}((0,1)^d) = H^{s}_{2}((0,1)^d) $ with $0 \leq s < m - \frac{1}{2} $ and $s \not \in \mathbb{N}_{0} + \frac{1}{2}$, see Theorem 3.11 in \cite{DaFKRaa}. Its proof is based on the fact that the spaces $ H^{s}((0,1)^d)   $ can be written as an intersection of function spaces that are tensor products of univariate Sobolev spaces $H^{s}(I)$ and univariate Lebesgue spaces $L_{2}(I)$. For results concerning this topic we refer to Section 3 in \cite{DauSte} and to Section 3 in \cite{ChDaFSte}. However, in \cite{DaFKRaa} only the Hilbert space case $ r = 2 $ has been investigated. Fortunately a recent result of Hansen and Sickel, see \cite{HaSi2022}, yields that also the more general Bessel-Potential spaces $   H^{s}_{r}((0,1)^2)$ with $ 1 < r < \infty   $ can be written as an intersection of function spaces which have a tensor product structure. In connection with that we have to deal with a special crossnorm, the so-called \textit{r-nuclear norm} $g_{r}$. In order to give a precise definition and to explicitly state the result of Hansen and Sickel, in what follows we have to recall some basic facts concerning tensor products of Banach spaces. For that purpose we follow \cite{Rya}, \cite{LightChe} and \cite{DeFlo}. Moreover, let us refer to Appendix B in \cite{SiUl}, where a brief summary regarding tensor products can be found.

\begin{defi}\label{def_tensor_prod}
\textbf{Tensor Products of Banach Spaces.}
\hspace{0,1 cm}

\begin{itemize}

\item[(i)]
Let $X$ and $Y$ be Banach spaces. Let $X^{'}$ be the dual space of $X$. To define the algebraic tensor product $ X \otimes Y$ we consider the set of all expressions $\sum_{j = 1}^{a} f_{j} \otimes g_{j} $ with $a \in \mathbb{N}$, $f_{j} \in X$ and $g_{j} \in Y$. We introduce an equivalence relation 
\begin{align*}
\sum_{j = 1}^{a} f_{j} \times g_{j} \sim \sum_{j = 1}^{b} u_{j} \otimes v_{j}
\end{align*}
with $b \in \mathbb{N}$, $u_{j} \in X$ and $v_{j} \in Y$, if both expressions generate the same operator $A : X^{'} \rightarrow Y$, that is
\begin{align*}
\sum_{j = 1}^{a} \varphi(f_{j})  g_{j} = \sum_{j = 1}^{b} \varphi(u_{j})  v_{j} \qquad \mbox{for all $\varphi \in X^{'}$.}
\end{align*}
Then the algebraic tensor product $X \otimes Y$ of $X$ and $Y$ is defined to be the set of all such equivalence classes.

\item[(ii)] 

The completion of $X \otimes Y$ with respect to a tensor norm $\alpha$ will be denoted by $X \otimes_{\alpha} Y$.

\end{itemize}
\end{defi}

In order to write bivariate Bessel-Potential spaces in terms of univariate function spaces by using tensor products we have to work with the so-called $r$-nuclear norm $g_{r}$. There is the following definition, see also Definition 1.45 in \cite{LightChe}.

\begin{defi}\label{def_r_nuclear_norm}
Let $X$ and $Y$ be Banach spaces. Let $Y^{'}$ be the dual space of $Y$. Let $1 \leq r \leq \infty $ and let $ 1/r + 1/r^{'} = 1$. Let $h \in X \otimes Y$ be given by 
\begin{equation}\label{eq_h_tensorrep}
h = \sum_{j = 1}^{a} f_{j} \otimes g_{j} , \qquad f_{j} \in X, \; g_{j} \in Y 
\end{equation}
with $a \in \mathbb{N}$. Then the $r$-nuclear tensor norm $g_{r} := g_{r}(\cdot , X, Y)$ is given by
\begin{align*}
g_{r}(h,X,Y) := \inf_{\eqref{eq_h_tensorrep}} \Big \{ \Big ( \sum_{j=1}^{a}   \Vert f_{j} \vert X \Vert^{r} \Big )^{\frac{1}{r}}  \times \sup \Big \{ \Big ( \sum_{j=1}^{a}  \vert \psi(g_{j}) \vert^{r^{'}} \Big )^{\frac{1}{r^{'}}}   : \psi \in Y^{'}, \Vert \psi \vert Y^{'} \Vert \leq 1 \Big \} \Big \} .
\end{align*}
Here the infimum is taken over all representations of $h$, see \eqref{eq_h_tensorrep}.
\end{defi}

\begin{rem}\label{rem_Nuclear_norm_alternative}
In the definition of the $r$-nuclear norm $g_{r}$, see Definition \ref{def_r_nuclear_norm}, it is possible to replace 
\begin{align*}
\sup \Big \{ \Big ( \sum_{j=1}^{a}  \vert \psi(g_{j}) \vert^{r^{'}} \Big )^{\frac{1}{r^{'}}}   : \psi \in Y^{'}, \Vert \psi \vert Y^{'} \Vert \leq 1 \Big \}
\end{align*}
by
\begin{align*}
\sup \Big \{  \Big \Vert \sum_{j=1}^{a} \lambda_{j} g_{j} \Big \vert Y \Big \Vert : \Big ( \sum_{j=1}^{a} \vert \lambda_{j} \vert^{r} \Big )^{\frac{1}{r}} \leq 1  \Big \} .
\end{align*}
For this observation we refer to \cite{LightChe}, see Lemma 1.44.  
\end{rem}
In Definition \ref{def_r_nuclear_norm} formally the $r$-nuclear tensor norm $  g_{r}(h,X,Y)  $ is only defined for functions $h$ of the form \eqref{eq_h_tensorrep} that can be written as a sum of tensor products $   f_{j} \otimes g_{j}   $. In order to define a norm on function spaces $ X \otimes_{g_{r}} Y    $ which are defined by completion, we have to extend Definition \ref{def_r_nuclear_norm} in the following way.

\begin{defi}\label{def_r_nuclear_norm_value}
Let $X$ and $Y$ be Banach spaces. Let $ 1 \leq r \leq \infty $ and $ f \in   X \otimes_{g_{r}} Y $. Then we define
\begin{align*}
\Vert f \vert  X \otimes_{g_{r}} Y  \Vert := \lim_{n \rightarrow \infty} \inf_{\eqref{eq_def_tensornorm_value11}} g_{r}(h_{n},X,Y) .
\end{align*} 
Here the infimum is taken over all sequences $ \{ h_{n} \}_{n \in \mathbb{N}}  $ of the form 
\begin{equation}\label{eq_def_tensornorm_value11}
h_{n} = \sum_{j = 1}^{a(n)} (f^{n})_{j} \otimes (g^{n})_{j} , \qquad (f^{n})_{j} \in X, \; (g^{n})_{j}  \in Y ,
\end{equation}
with $a(n) \in \mathbb{N}$, whereby we have
\begin{align*}
\lim_{n \rightarrow \infty} h_{n} = f \qquad \qquad \mbox{with convergence in $ X \otimes_{g_{r}} Y  $.}
\end{align*}
\end{defi}

The $r$-nuclear norm can be used to write multivariate Lebesgue spaces as a tensor product of univariate Lebesgue spaces. For that purpose let us recall the following result, see Corollary 1.52 in \cite{LightChe}.

\begin{lem}\label{lem_Lp_pnuc_tens}
Let $ 1 \leq r < \infty $ and $d \in \mathbb{N}$. Then we have
\begin{align*}
L_{r}((0,1)^d) = L_{r}((0,1)) \otimes_{g_{r}} L_{r}((0,1)) \otimes_{g_{r}} \ldots \otimes_{g_{r}} L_{r}((0,1)) .
\end{align*}
\end{lem}

In order to write bivariate Bessel-Potential spaces in terms of tensor products we have to deal with intersections of certain function spaces. In connection with that we will use the definition below.

\begin{defi}\label{defi_norm_intersec}
Let $n \in \mathbb{N}$. Let $X_{j}$ with $j = 1, 2, \ldots , n$ be Banach spaces with associated norms $\Vert \cdot \vert X_{j} \Vert $. Then we equip the intersection space $\bigcap_{j=1}^{n} X_{j}$ with the quasi-norm 
\begin{align*}
\Big \Vert \cdot \Big \vert \bigcap_{j=1}^{n} X_{j} \Big \Vert := \sum_{j = 1}^{n} \Vert \cdot \vert X_{j} \Vert .
\end{align*}
\end{defi}

Now we are well-prepared to formulate a result of Hansen and Sickel, see Theorem 2.16 in \cite{HaSi2022}. It will be very important for us later on.

\begin{prop}\label{prop_tensor_Hsr_d2}
Let $1 < r < \infty$ and $s > 0$. Then it holds
\begin{align*}
H^{s}_{r}((0,1)^{2}) = ( H^{s}_{r}((0,1)) \otimes_{g_{r}} L_{r}((0,1)) ) \cap ( L_{r}((0,1)) \otimes_{g_{r}}  H^{s}_{r}((0,1)) )
\end{align*}
in the sense of equivalent norms. 
\end{prop}
Let us remark that Proposition \ref{prop_tensor_Hsr_d2} can be seen as a continuation of \cite{GriOs}, see Example 2 in Section 2. There the special case $r = 2$ and $s = 1$ has been investigated.

\section{Bivariate Quarklets via Tensor Product Methods}\label{sec_bessel_multi}

In this section we describe bivariate Bessel-Potential spaces $H^{s}_{r}((0,1)^2)$ defined on squares via quarklets. Thereto we use bivariate quarklets that are constructed out of univariate quarklets via tensor product methods. To this end let us introduce some additional notation. We follow \cite{DaFKRaa}, see Section 3.3. In what follows $i \in \{ 1, 2  \}$ always refers to the $i-$th Cartesian direction. The boundary conditions in direction $i$ are described by $  \sigma_{i} \in \{ 0 , \ldots  , \lfloor   s + 1 - \frac{2}{r} - \varepsilon  \rfloor   \}^{2}   $ with $\varepsilon > 0$ arbitrary small. The collection of all boundary conditions is given by
\begin{equation}
\boldsymbol{\sigma} := ( \sigma_{1}, \sigma_{2} ) = ( (\sigma_{1}^{l} , \sigma_{1}^{r}) , (\sigma_{2}^{l} , \sigma_{2}^{r}) ) \in ( \{ 0 , \ldots  , \lfloor   s + 1 - \frac{2}{r} - \varepsilon  \rfloor   \}^{2} )^{2} . 
\end{equation}
The collection $ \boldsymbol{\sigma}   $ can be used to define bivariate Bessel-Potential spaces with certain boundary conditions. So $  (H^{s}_{r}((0,1)^2))_{\boldsymbol{\sigma}}  $ is defined to be the space of all $  H^{s}_{r}((0,1)^2)$-functions whose normal derivatives of up to orders $ \sigma_{i}^{l}  $ and $  \sigma_{i}^{r}  $ vanish at the facets $ I^{i-1} \times \{ 0 \}  \times I^{2-i}    $ and $ I^{i-1} \times \{ 1 \}  \times I^{2-i}    $, respectively, for $ i \in \{ 1 , 2 \}   $. Here the boundary conditions are meaningful due to Theorem 3.3.1 in \cite{SiTri}. The spaces $  (H^{s}_{r}((0,1)^2))_{\boldsymbol{\sigma}}  $ are equipped with the usual norm. That means 
\begin{equation}\label{def_norm_boundcond1_d=2}
\Vert f \vert  (H^{s}_{r}((0,1)^{2}))_{\boldsymbol{\sigma}}   \Vert := \inf \Big\{\|  g \vert H^{s}_{r}( \mathbb{R}^{2} ) \|: \; g_{|_{(0,1)^{2}}} =f \; , \mbox{$g$ fulfills conditions $\boldsymbol{\sigma}$ at $\partial (0,1)^2$ }  \Big\} \, .
\end{equation}
There is the following counterpart of Proposition \ref{prop_tensor_Hsr_d2} involving boundary conditions.

\begin{prop}\label{prop_tensor_Hsr_d2_boundco}
Let $1 < r < \infty$ and $s > 0$. Let $ \boldsymbol{\sigma} = ( \sigma_{1}, \sigma_{2} )  $ be a collection of boundary conditions. Then it holds
\begin{align*}
(H^{s}_{r}((0,1)^{2}))_{\boldsymbol{\sigma}} = ( (H^{s}_{r}((0,1)))_{\sigma_{1}} \otimes_{g_{r}} L_{r}((0,1)) ) \cap ( L_{r}((0,1)) \otimes_{g_{r}}  (H^{s}_{r}((0,1)))_{\sigma_{2}} )
\end{align*}
in the sense of equivalent norms. 
\end{prop}

\begin{proof}
This result can be proved with similar methods as Proposition \ref{prop_tensor_Hsr_d2}, see Theorem 2.16 in \cite{HaSi2022}. In the proof given there one has to replace $  H^{s}_{r}((0,1)^{2})  $ by $  (H^{s}_{r}((0,1)^{2}))_{\boldsymbol{\sigma}}  $  and $ H^{s}_{r}((0,1))   $ either by $ (H^{s}_{r}((0,1)))_{\sigma_{1}}   $ or $  (H^{s}_{r}((0,1)))_{\sigma_{2}}  $ in an appropriate way. Those modifications are possible since the proof of Theorem 2.16 in \cite{HaSi2022} is completely based on norm estimates, for that the particular function values on $\partial (0,1)^2$ do not play any role. 
\end{proof}

\begin{rem}
A counterpart of Proposition \ref{prop_tensor_Hsr_d2_boundco} for Sobolev spaces with boundary conditions can be found in \cite{DauSte}, see Section 3, and in \cite{ChDaFSte}, see Section 3.
\end{rem}

Recall that each univariate quarklet $\psi_{p,j,k}^{\sigma}$ can be identified via a triple $\lambda = (p,j,k)$ with $p \in \mathbb{N}_{0}$, $j \geq j_{0}-1$ and $ k \in \nabla_{j, \sigma}$. Sometimes $\lambda$ also is called a quarklet index. To obtain bivariate quarklets we require quarklet indices for each Cartesian direction $i$. They are denoted by $\lambda_{i} = (p_{i},j_{i},k_{i})$. We put $\boldsymbol{\lambda} := ( \lambda_{1}, \lambda_{2} )$. The index set for the whole quarklet system concerning direction $i$ is given by $\nabla_{\sigma_{i}}$, see also \eqref{index_quarkl_full1}. Recall that it also depends on the boundary conditions assigned with direction $i$. Now for quarklet indices $\lambda_{i}$ and boundary conditions $\sigma_{i}$ we define bivariate quarklets via tensor products of univariate quarklets. That means we put
\begin{equation}\label{def_tensor_quarklet}
\boldsymbol{\psi^{\sigma}_{\lambda}} :=  \psi_{\lambda_{1}}^{\sigma_{1}} \otimes \psi_{\lambda_{2}}^{\sigma_{2}} .
\end{equation}
To address the bivariate quarklets we define the index set
\begin{equation}\label{def_tensorqua_index}
\boldsymbol{\nabla_{\sigma}} :=  \nabla_{\sigma_{1}} \times \nabla_{\sigma_{2}} .
\end{equation}
Then the collection of all bivariate quarklets that can be obtained via tensor products as described above is given by
\begin{equation}\label{def_collect_multqua}
\boldsymbol{\Psi_{\sigma}} :=  \Psi_{\sigma_{1}} \otimes \Psi_{\sigma_{2}} = \Big \{ \boldsymbol{\psi^{\sigma}_{\lambda}} : \boldsymbol{\lambda}  \in  \boldsymbol{\nabla_{\sigma}} \Big \} .
\end{equation}
Here $\Psi_{\sigma_{i}}$ refers to the whole system of univariate quarklets concerning direction $i$, see \eqref{def_quarklsy_ful1}.

\subsection{Quarklet Characterizations for Bivariate Bessel-Potential Spaces $H^{s}_{r}((0,1)^2)$}

In order to describe bivariate Bessel-Potential spaces in terms of quarklets we require some additional notation. It is connected with the characterization of univariate Bessel-Potential spaces via quarklets given in Theorem \ref{mainresult1_Hsr_d=1}.

\begin{defi}\label{def_seq_sp_mitp}
Let $  1 < r < \infty $ and $ m \in \mathbb{N}  $ with $ m \geq 2  $. Let $ 0 \leq s  < m - 1 $ and $\delta > 1$. Let $\nabla_{\sigma}$ be the index set defined in \eqref{index_quarkl_full1}. Let $ \{ c_{p,j,k}  \}_{(p,j,k) \in \nabla_{\sigma}} \subset \mathbb{C}   $ be a sequence. Then we define
\begin{align*}
& \Vert  \{ c_{p,j,k}  \}_{(p,j,k) \in \nabla_{\sigma}}  \vert  h^{s,m}_{r,\delta}(\nabla_{\sigma})  \Vert \\
& \qquad \qquad :=  \Big\| \Big (   \sum_{(p,j,k) \in \nabla_{\sigma}} (p+1)^{\sgn(s) 4m + 2   \delta }  2^{2 j s} 2^{j}    | c_{p,j,k}|^{2}  | \tilde{\chi}_{j,k}(x) |   \Big )^{\frac{1}{2}}\Big|  L_{r} (I) \Big\| .
\end{align*}
Here $\sgn$ refers to the sign function. Recall $I = (0,1) $.
\end{defi}

Motivated by Propositions \ref{prop_tensor_Hsr_d2} and \ref{prop_tensor_Hsr_d2_boundco} below we investigate the spaces $ (H^{s}_{r}(I))_{\sigma_{1}} \otimes_{g_{r}} L_{r}(I)  $ and $ L_{r}(I) \otimes_{g_{r}}  (H^{s}_{r}(I))_{\sigma_{2}} $. It will be an important intermediate step for us to describe these function spaces in terms of bivariate quarklets. Let us start with $ (H^{s}_{r}(I))_{\sigma_{1}} \otimes_{g_{r}} L_{r}(I) $.

\begin{prop}\label{prop_char_Hs_Lr}
Let $d = 2$. Let $  1 < r < \infty $ and $ m \in \mathbb{N}  $ with $ m \geq 2  $. Let $ \tilde{m} \in \mathbb{N}_{0}   $ with $ m + \tilde{m} \in 2 \mathbb{N}   $ and $  \tilde{m} > 5 m + 12 $. Let $ 0 < s  < m - 1 $ and $\delta_{1}, \delta_{2} > 1$. Let $ \sigma_{1}  $ be a boundary condition and $  \sigma_{2} = ( 0 , 0 )$. Put $  \boldsymbol{{\sigma}^{1}} := ( \sigma_{1}, \sigma_{2} ) $. Let $f \in L_{r}((0,1)^{2})$. Then we have $ f \in  (H^{s}_{r}(I))_{\sigma_{1}} \otimes_{g_{r}} L_{r}(I)  $ if and only if there exists a sequence of functions $(h_{n})_{n \in \mathbb{N}}$ of the form 
\begin{equation}\label{qua_rep_d=2_1_a}
h_{n} = \sum_{( \lambda_{1}, \lambda_{2}) \in \boldsymbol{\nabla_{{\sigma}^{1}}} }  \Big (  \sum_{\ell = 1}^{a(n)}    (c^{n}_{1, \ell})_{\lambda_{1} }    (c^{n}_{2, \ell})_{\lambda_{2}}   \Big )  (  \psi_{\lambda_{1}}^{\sigma_{1}}  \otimes  \psi_{\lambda_{2}}^{\sigma_{2}} )
\end{equation}
with $a(n) \in \mathbb{N}$, such that
\begin{equation}\label{qua_rep_d=2_1_b}
\Vert  \boldsymbol{c}^{1}  \vert  h^{s,m}_{r, \boldsymbol{\delta}}(\boldsymbol{\nabla_{{\sigma}^{1}}} , 1)  \Vert :=  \Big ( \sum_{\ell =1}^{a(n)}  \Vert  (c^{n}_{1, \ell})_{\lambda_{1}} \vert  h^{s,m}_{r,\delta_{1}}(\nabla_{\sigma_{1}})  \Vert   \Big )  \times  \Big (  \sum_{\ell =1}^{a(n)}    \Vert   (c^{n}_{2, \ell})_{\lambda_{2}}  \vert  h^{0,m}_{r,\delta_{2}}(\nabla_{\sigma_{2}})  \Vert \Big ) < \infty
\end{equation}
and
\begin{equation}\label{qua_rep_d=2_1_d}
\lim_{n \rightarrow \infty} \Vert f - h_{n} \vert (H^{s}_{r}(I))_{\sigma_{1}} \otimes_{g_{r}} L_{r}(I) \Vert = 0.
\end{equation}
Moreover, the norms $ \Vert f  \vert (H^{s}_{r}(I))_{\sigma_{1}} \otimes_{g_{r}} L_{r}(I) \Vert$ and 
\begin{equation}\label{qua_rep_d=2_1_e}
\begin{split}
& \Vert f  \vert (H^{s}_{r}(I))_{\sigma_{1}} \otimes_{g_{r}} L_{r}(I) \Vert^{( \boldsymbol{\delta}, \boldsymbol{\nabla_{{\sigma}^{1}}})}  \\
& \qquad \qquad := \lim_{n \rightarrow \infty} \inf_{\eqref{qua_rep_d=2_1_f}} \Big \{  \Big ( \sum_{\ell =1}^{a(n)}  \inf_{\eqref{qua_rep_d=2_1_g}} \Vert  (c^{n}_{1, \ell})_{\lambda_{1}} \vert  h^{s,m}_{r,\delta_{1}}(\nabla_{\sigma_{1}})  \Vert^{r}   \Big )^{\frac{1}{r}} \\
& \qquad \qquad \qquad \qquad \times \sup \Big \{ \inf_{\eqref{qua_rep_d=2_1_h}} \Vert  (c^{n}_{2})_{\lambda_{2}} \vert  h^{0,m}_{r,\delta_{2}}(\nabla_{\sigma_{2}})  \Vert  : \Big ( \sum_{\ell =1}^{a(n)} \vert \mu_{\ell} \vert^{r} \Big )^{\frac{1}{r}} \leq 1  \Big \} \Big \} 
\end{split}
\end{equation}
are equivalent. Here the infimum \eqref{qua_rep_d=2_1_f} is taken over all sequences of functions $(h_{n})_{n \in \mathbb{N}}$ of the form 
\begin{equation}\label{qua_rep_d=2_1_f}
h_{n} = \sum_{\ell = 1}^{a(n)} (f^{n})_{\ell} \otimes (g^{n})_{\ell} , \qquad (f^{n})_{\ell} \in (H^{s}_{r}(I))_{\sigma_{1}}, \; (g^{n})_{\ell} \in L_{r}(I) 
\end{equation}
with 
\begin{align*}
\lim_{n \rightarrow \infty} \Vert f - h_{n} \vert (H^{s}_{r}(I))_{\sigma_{1}} \otimes_{g_{r}} L_{r}(I) \Vert = 0.
\end{align*}
The infimum \eqref{qua_rep_d=2_1_g} is taken over all sequences $ \{ (c^{n}_{1, \ell})_{\lambda_{1}}  \}_{\lambda_{1} \in \nabla_{\sigma_{1}}}   $ such that
\begin{equation}\label{qua_rep_d=2_1_g}
(f^{n})_{\ell} = \sum_{\lambda_{1} \in \nabla_{\sigma_{1}}}   (c^{n}_{1, \ell})_{\lambda_{1}} \psi_{\lambda_{1}}^{\sigma_{1}} .
\end{equation}
The infimum \eqref{qua_rep_d=2_1_h} is taken over all sequences $ \{ (c^{n}_{2})_{\lambda_{2}}  \}_{\lambda_{2} \in \nabla_{\sigma_{2}}}   $ such that
\begin{equation}\label{qua_rep_d=2_1_h}
\sum_{\ell =1}^{a(n)} \mu_{\ell} (g^{n})_{\ell} = \sum_{\lambda_{2} \in \nabla_{\sigma_{2}}}   (c^{n}_{2})_{\lambda_{2}} \psi_{\lambda_{2}}^{\sigma_{2}} .
\end{equation}
\end{prop}

\begin{rem}
In the notation $  \Vert  \boldsymbol{c}^{1}  \vert  h^{s,m}_{r, \boldsymbol{\delta}}(\boldsymbol{\nabla_{{\sigma}^{1}}} , 1)  \Vert    $ introduced in \eqref{qua_rep_d=2_1_b} the symbol $\boldsymbol{c}^{1}$ collects all sequences showing up there. Moreover, we use the abbreviation $\boldsymbol{\delta} := ( \delta_{1}, \delta_{2} )$.
\end{rem}

\begin{proof}
\textit{Step 1.}
Let $f \in (H^{s}_{r}(I))_{\sigma_{1}} \otimes_{g_{r}} L_{r}(I) $. We want to prove that there exists a sequence of functions $(h_{n})_{n \in \mathbb{N}}$ of the form \eqref{qua_rep_d=2_1_a}, such that \eqref{qua_rep_d=2_1_b} and \eqref{qua_rep_d=2_1_d} are fulfilled. For that purpose at first we can use Definitions \ref{def_tensor_prod}, \ref{def_r_nuclear_norm} and \ref{def_r_nuclear_norm_value} to find that there exists a sequence of functions $(h_{n})_{n \in \mathbb{N}}$, such that
\begin{equation}\label{eq__h_tensorrep_Hs_Lr_conveq}
\lim_{n \rightarrow \infty} \Vert f - h_{n} \vert (H^{s}_{r}(I))_{\sigma_{1}} \otimes_{g_{r}} L_{r}(I) \Vert = 0 ,
\end{equation}
at which each function $ h_{n}  $ has a representation 
\begin{equation}\label{eq_h_tensorrep_Hs_Lr}
h_{n} = \sum_{\ell = 1}^{a(n)} (f^{n})_{\ell} \otimes (g^{n})_{\ell} , \qquad (f^{n})_{\ell} \in (H^{s}_{r}(I))_{\sigma_{1}}, \; (g^{n})_{\ell} \in L_{r}(I) ,
\end{equation}
whereat $a(n) \in \mathbb{N}$ depends on $n$. For each $n \in \mathbb{N}$ we have $ g_{r}(h_{n}, (H^{s}_{r}(I))_{\sigma_{1}}, L_{r}(I)) < \infty $. A combination of \eqref{eq_h_tensorrep_Hs_Lr} and Theorem \ref{mainresult1_Hsr_d=1} yields that $h_{n}$ can be written as
\begin{align*}
h_{n}  = \sum_{\ell = 1}^{a(n)} (f^{n})_{\ell} \otimes (g^{n})_{\ell}  = \sum_{\ell = 1}^{a(n)} \Big ( \sum_{\lambda_{1} \in \nabla_{\sigma_{1}}}   (c^{n}_{1, \ell})_{\lambda_{1}} \psi_{\lambda_{1}}^{\sigma_{1}} \Big ) \otimes (g^{n})_{\ell} 
\end{align*}
with
\begin{equation}\label{eq_prof_auxeq11}
 \sum_{\ell =1}^{a(n)}  \Vert  (c^{n}_{1, \ell})_{\lambda_{1}} \vert  h^{s,m}_{r,\delta_{1}}(\nabla_{\sigma_{1}})  \Vert < \infty .
\end{equation}
Recall that each $(g^{n})_{\ell} \in L_{r}(I)$ has a representation in terms of univariate quarklets, see Corollary \ref{mainresult1_Lr_d=1}. Therefore we also can find a representation
\begin{equation}\label{eq_f_zwei_tensor}
h_{n}  = \sum_{\ell = 1}^{a(n)} \Big ( \sum_{\lambda_{1} \in \nabla_{\sigma_{1}}}   (c^{n}_{1, \ell})_{\lambda_{1}} \psi_{\lambda_{1}}^{\sigma_{1}} \Big ) \otimes  \Big ( \sum_{\lambda_{2} \in \nabla_{\sigma_{2}}}   (c^{n}_{2, \ell})_{\lambda_{2}} \psi_{\lambda_{2}}^{\sigma_{2}} \Big ) 
\end{equation}
with
\begin{align*}
\sum_{\ell =1}^{a(n)}    \Vert   (c^{n}_{2, \ell})_{\lambda_{2}}  \vert  h^{0,m}_{r,\delta_{2}}(\nabla_{\sigma_{2}})  \Vert < \infty .
\end{align*}
In combination with \eqref{eq_prof_auxeq11} this yields  \eqref{qua_rep_d=2_1_b}. In order to obtain \eqref{qua_rep_d=2_1_a} in what follows we rearrange \eqref{eq_f_zwei_tensor}. We observe that since $\nabla_{\sigma_{i}} \subset \mathbb{N}_{0} \times \mathbb{N}_{0} \times \mathbb{Z} $ there exists a bijection $\pi_{i} : \nabla_{\sigma_{i}} \rightarrow \mathbb{N}_{0}$. Consequently we can use the Cauchy product to find
\begin{align*}
h_{n}  = \sum_{\ell = 1}^{a(n)} \Big ( \sum_{\pi_{1}(\lambda_{1}) = 0}^{\infty} \sum_{\pi_{2}(\lambda_{2}) = 0}^{\pi_{1}(\lambda_{1})}  ( (c^{n}_{1, \ell})_{\pi_{1}^{-1}( \pi_{1}(\lambda_{1})-\pi_{2}(\lambda_{2})) } \psi_{\pi_{1}^{-1}( \pi_{1}(\lambda_{1})-\pi_{2}(\lambda_{2}))}^{\sigma_{1}} ) \otimes (   (c^{n}_{2, \ell})_{\lambda_{2}} \psi_{\lambda_{2}}^{\sigma_{2}} ) \Big ) .
\end{align*}
The convergence is ensured by Theorem \ref{mainresult1_Hsr_d=1} and Corollary \ref{mainresult1_Lr_d=1}. It yields that we also can write
\begin{align*}
h_{n} & =   \sum_{\pi_{1}(\lambda_{1}) = 0}^{\infty} \sum_{\pi_{2}(\lambda_{2}) = 0}^{\pi_{1}(\lambda_{1})}  \sum_{\ell = 1}^{a(n)} \Big (  ( (c^{n}_{1, \ell})_{\pi_{1}^{-1}( \pi_{1}(\lambda_{1})-\pi_{2}(\lambda_{2})) } \psi_{\pi_{1}^{-1}( \pi_{1}(\lambda_{1})-\pi_{2}(\lambda_{2}))}^{\sigma_{1}} ) \otimes (   (c^{n}_{2, \ell})_{\lambda_{2}} \psi_{\lambda_{2}}^{\sigma_{2}} ) \Big ) \\
& =   \sum_{\pi_{1}(\lambda_{1}) = 0}^{\infty} \sum_{\pi_{2}(\lambda_{2}) = 0}^{\pi_{1}(\lambda_{1})} \Big [ \Big (  \psi_{\pi_{1}^{-1}( \pi_{1}(\lambda_{1})-\pi_{2}(\lambda_{2}))}^{\sigma_{1}}  \otimes \psi_{\lambda_{2}}^{\sigma_{2}} \Big )  \Big (  \sum_{\ell = 1}^{a(n)}   ( (c^{n}_{1, \ell})_{\pi_{1}^{-1}( \pi_{1}(\lambda_{1})-\pi_{2}(\lambda_{2})) } (   (c^{n}_{2, \ell})_{\lambda_{2}}  ) \Big ) \Big ] .
\end{align*}
This representation has the desired form, see \eqref{qua_rep_d=2_1_a}, whereby here we have a particular arrangement of the coefficients. In combination with \eqref{eq__h_tensorrep_Hs_Lr_conveq}, which yields \eqref{qua_rep_d=2_1_d}, this completes Step 1 of the proof.

\textit{Step 2.} 
In order to obtain the equivalent norm again let $f \in (H^{s}_{r}(I))_{\sigma_{1}} \otimes_{g_{r}} L_{r}(I) $. As before the sequence of functions $(h_{n})_{n \in \mathbb{N}}$ is given by \eqref{eq_h_tensorrep_Hs_Lr}, whereby \eqref{eq__h_tensorrep_Hs_Lr_conveq} is fulfilled. Then we can use Definition \ref{def_r_nuclear_norm}  and Remark \ref{rem_Nuclear_norm_alternative} to find
\begin{align*}
& g_{r}(h_{n}, (H^{s}_{r}(I))_{\sigma_{1}}, L_{r}(I))  \\
& \quad = \inf_{\eqref{eq_h_tensorrep_Hs_Lr}} \Big \{ \Big ( \sum_{\ell =1}^{a(n)}   \Vert (f^{n})_{\ell} \vert (H^{s}_{r}(I))_{\sigma_{1}} \Vert^{r} \Big )^{\frac{1}{r}}  \times \sup \Big \{  \Big \Vert \sum_{\ell =1}^{a(n)} \mu_{\ell} (g^{n})_{\ell} \Big \vert L_{r}(I) \Big \Vert : \Big ( \sum_{\ell =1}^{a(n)} \vert \mu_{\ell} \vert^{r} \Big )^{\frac{1}{r}} \leq 1  \Big \} \Big \} .
\end{align*}
Here the infimum is taken over all representations of $h_{n}$ of the form \eqref{eq_h_tensorrep_Hs_Lr}. To investigate the first term we can apply Theorem \ref{mainresult1_Hsr_d=1}. It yields
\begin{align*}
& \Big ( \sum_{\ell=1}^{a(n)}   \Vert (f^{n})_{\ell} \vert (H^{s}_{r}(I))_{\sigma_{1}} \Vert^{r} \Big )^{\frac{1}{r}} \\
& \qquad \approx \Big ( \sum_{\ell =1}^{a(n)}  \inf_{\eqref{rep_main_Hsr_d=1_tensor_d2_H1}} \Big\| \Big[   \sum_{(p,j,k) = \lambda_{1} \in \nabla_{\sigma_{1}}} (p+1)^{4m + 2   \delta_{1} }  2^{2 j s} 2^{j}    | (c^{n}_{1, \ell})_{\lambda_{1}}|^{2}  | \tilde{\chi}_{j,k}(x) |   \Big]^{\frac{1}{2}}\Big|  L_{r} (I) \Big\|^{r} \Big )^{\frac{1}{r}} .
\end{align*}
Here the infimum is taken over all sequences $ \{ (c^{n}_{1, \ell})_{\lambda_{1}}  \}_{\lambda_{1} \in \nabla_{\sigma_{1}}}   $ such that
\begin{equation}\label{rep_main_Hsr_d=1_tensor_d2_H1}
(f^{n})_{\ell} = \sum_{\lambda_{1} \in \nabla_{\sigma_{1}}}   (c^{n}_{1, \ell})_{\lambda_{1}} \psi_{\lambda_{1}}^{\sigma_{1}}
\end{equation}
with convergence in $ \mathcal{S}'(I)  $. To deal with the second term we use Corollary \ref{mainresult1_Lr_d=1}. We obtain
\begin{align*}
& \sup \Big \{  \Big \Vert \sum_{\ell =1}^{a(n)} \mu_{\ell} (g^{n})_{\ell} \Big \vert L_{r}(I) \Big \Vert : \Big ( \sum_{\ell =1}^{a(n)} \vert \mu_{\ell} \vert^{r} \Big )^{\frac{1}{r}} \leq 1  \Big \} \\
& \quad \approx \sup \Big \{ \inf_{\eqref{rep_main_Hsr_d=1_d2_L2}} \Big\| \Big[   \sum_{(p,j,k) = \lambda_{2} \in \nabla_{\sigma_{2}}} (p+1)^{ 2   \delta_{2} }  2^{j}    | (c^{n}_{2})_{\lambda_{2}}|^{2}  | \tilde{\chi}_{j,k}(x) |   \Big]^{\frac{1}{2}}\Big|  L_{r} (I) \Big\|  : \Big ( \sum_{\ell =1}^{a(n)} \vert \mu_{\ell} \vert^{r} \Big )^{\frac{1}{r}} \leq 1  \Big \} .
\end{align*}
Here the infimum is taken over all sequences $ \{ (c^{n}_{2})_{\lambda_{2}}  \}_{\lambda_{2} \in \nabla_{\sigma_{2}}}   $ such that
\begin{equation}\label{rep_main_Hsr_d=1_d2_L2}
\sum_{\ell =1}^{a(n)} \mu_{\ell} (g^{n})_{\ell} = \sum_{\lambda_{2} \in \nabla_{\sigma_{2}}}   (c^{n}_{2})_{\lambda_{2}} \psi_{\lambda_{2}}^{\sigma_{2}}
\end{equation}
with convergence in $ \mathcal{S}'(I)  $. Consequently using Definition \ref{def_seq_sp_mitp} we also find
\begin{align*}
& g_{r}(h_{n}, (H^{s}_{r}(I))_{\sigma_{1}}, L_{r}(I))  \\
& \qquad  \approx \inf_{\eqref{eq_h_tensorrep_Hs_Lr}} \Big \{  \Big ( \sum_{\ell =1}^{a(n)}  \inf_{\eqref{rep_main_Hsr_d=1_tensor_d2_H1}} \Vert  (c^{n}_{1, \ell})_{\lambda_{1}} \vert  h^{s,m}_{r,\delta_{1}}(\nabla_{\sigma_{1}})  \Vert^{r}   \Big )^{\frac{1}{r}} \\
& \qquad \qquad \qquad  \times \sup \Big \{ \inf_{\eqref{rep_main_Hsr_d=1_d2_L2}} \Vert  (c^{n}_{2})_{\lambda_{2}} \vert  h^{0,m}_{r,\delta_{2}}(\nabla_{\sigma_{2}})  \Vert  : \Big ( \sum_{\ell =1}^{a(n)} \vert \mu_{\ell} \vert^{r} \Big )^{\frac{1}{r}} \leq 1  \Big \} \Big \} .
\end{align*}
The above equivalence in combination with \eqref{eq__h_tensorrep_Hs_Lr_conveq} shows
\begin{align*}
& \lim_{n \rightarrow \infty} \inf_{\eqref{eq_h_tensorrep_Hs_Lr}} \Big \{  \Big ( \sum_{\ell =1}^{a(n)}  \inf_{\eqref{rep_main_Hsr_d=1_tensor_d2_H1}} \Vert  (c^{n}_{1, \ell})_{\lambda_{1}} \vert  h^{s,m}_{r,\delta_{1}}(\nabla_{\sigma_{1}})  \Vert^{r}   \Big )^{\frac{1}{r}} \\
& \qquad \qquad   \times \sup \Big \{ \inf_{\eqref{rep_main_Hsr_d=1_d2_L2}} \Vert  (c^{n}_{2})_{\lambda_{2}} \vert  h^{0,m}_{r,\delta_{2}}(\nabla_{\sigma_{2}})  \Vert  : \Big ( \sum_{\ell =1}^{a(n)} \vert \mu_{\ell} \vert^{r} \Big )^{\frac{1}{r}} \leq 1  \Big \} \Big \}  \\
& \qquad \qquad \qquad \qquad \lesssim  \lim_{n \rightarrow \infty} g_{r}(h_{n}, (H^{s}_{r}(I))_{\sigma_{1}}, L_{r}(I))   \\
& \qquad \qquad \qquad \qquad =  \lim_{n \rightarrow \infty} \Vert  h_{n} \vert (H^{s}_{r}(I))_{\sigma_{1}} \otimes_{g_{r}} L_{r}(I) \Vert    \\
& \qquad \qquad \qquad \qquad =  \lim_{n \rightarrow \infty} \Vert  h_{n} - f + f \vert (H^{s}_{r}(I))_{\sigma_{1}} \otimes_{g_{r}} L_{r}(I) \Vert    \\
& \qquad \qquad \qquad \qquad \leq  \lim_{n \rightarrow \infty} \Vert  h_{n} - f  \vert (H^{s}_{r}(I))_{\sigma_{1}} \otimes_{g_{r}} L_{r}(I) \Vert +  \Vert  f  \vert (H^{s}_{r}(I))_{\sigma_{1}} \otimes_{g_{r}} L_{r}(I) \Vert   \\
& \qquad \qquad \qquad \qquad =   \Vert  f  \vert (H^{s}_{r}(I))_{\sigma_{1}} \otimes_{g_{r}} L_{r}(I) \Vert   .
\end{align*}
This is the desired estimate. Hence Step 2 of the proof is complete.

\textit{Step 3.}
Now let $f \in L_{r}((0,1)^{2})$ such that there exists a sequence of functions $(h_{n})_{n \in \mathbb{N}}$ of the form \eqref{qua_rep_d=2_1_a}, whereat \eqref{qua_rep_d=2_1_b} and \eqref{qua_rep_d=2_1_d} are fulfilled. Then we observe
\begin{align*}
& \Vert f  \vert (H^{s}_{r}(I))_{\sigma_{1}} \otimes_{g_{r}} L_{r}(I) \Vert \\ 
& \qquad  = \lim_{n \rightarrow \infty} \Vert f - h_{n} + h_{n}  \vert (H^{s}_{r}(I))_{\sigma_{1}} \otimes_{g_{r}} L_{r}(I) \Vert \\
& \qquad \leq \lim_{n \rightarrow \infty}  \Vert f - h_{n}   \vert (H^{s}_{r}(I))_{\sigma_{1}} \otimes_{g_{r}} L_{r}(I) \Vert + \lim_{n \rightarrow \infty} \Vert  h_{n}  \vert (H^{s}_{r}(I))_{\sigma_{1}} \otimes_{g_{r}} L_{r}(I) \Vert \\
& \qquad = \lim_{n \rightarrow \infty} \Vert f - h_{n}   \vert (H^{s}_{r}(I))_{\sigma_{1}} \otimes_{g_{r}} L_{r}(I) \Vert + \lim_{n \rightarrow \infty} g_{r}(h_{n}, (H^{s}_{r}(I))_{\sigma_{1}}, L_{r}(I)) .
\end{align*}
We can use \eqref{qua_rep_d=2_1_d} to get
\begin{align*}
\Vert f  \vert (H^{s}_{r}(I))_{\sigma_{1}} \otimes_{g_{r}} L_{r}(I) \Vert \leq  \lim_{n \rightarrow \infty} g_{r}(h_{n}, (H^{s}_{r}(I))_{\sigma_{1}}, L_{r}(I)) .
\end{align*}
To continue with similar methods as in Step 1 of the proof using Theorem \ref{mainresult1_Hsr_d=1} and Corollary \ref{mainresult1_Lr_d=1} we find
\begin{equation}\label{eq_Prop6_Step3_est11}
\begin{split}
& \Vert f  \vert (H^{s}_{r}(I))_{\sigma_{1}} \otimes_{g_{r}} L_{r}(I) \Vert \\
&  \lesssim \lim_{n \rightarrow \infty}  \Big \{  \Big ( \sum_{\ell =1}^{a(n)}   \Vert  (c^{n}_{1, \ell})_{\lambda_{1}} \vert  h^{s,m}_{r,\delta_{1}}(\nabla_{\sigma_{1}})  \Vert^{r}   \Big )^{\frac{1}{r}}  \times \sup \Big \{  \Vert  (c^{n}_{2})_{\lambda_{2}} \vert  h^{0,m}_{r,\delta_{2}}(\nabla_{\sigma_{2}})  \Vert  : \Big ( \sum_{\ell =1}^{a(n)} \vert \mu_{\ell} \vert^{r} \Big )^{\frac{1}{r}} \leq 1  \Big \} \Big \} .
\end{split}
\end{equation}
Notice that we do not require the infima here since we already have found a suitable representation, see \eqref{qua_rep_d=2_1_a}. Now we utilize the relation
\begin{align*}
\sum_{\lambda_{2} \in \nabla_{\sigma_{2}}}   \sum_{\ell =1}^{a(n)} \mu_{\ell} (c^{n}_{2, \ell})_{\lambda_{2}} \psi_{\lambda_{2}}^{\sigma_{2}}  = \sum_{\lambda_{2} \in \nabla_{\sigma_{2}}}   (c^{n}_{2})_{\lambda_{2}} \psi_{\lambda_{2}}^{\sigma_{2}} .
\end{align*}
It yields
\begin{align*}
& \sup \Big \{  \Vert  (c^{n}_{2})_{\lambda_{2}} \vert  h^{0,m}_{r,\delta_{2}}(\nabla_{\sigma_{2}})  \Vert  : \Big ( \sum_{\ell =1}^{a(n)} \vert \mu_{\ell} \vert^{r} \Big )^{\frac{1}{r}} \leq 1  \Big \} \\
& \qquad = \sup \Big \{ \Big \Vert  \sum_{\ell =1}^{a(n)} \mu_{\ell} (c^{n}_{2, \ell})_{\lambda_{2}} \Big \vert  h^{0,m}_{r,\delta_{2}}(\nabla_{\sigma_{2}}) \Big \Vert  : \Big ( \sum_{\ell =1}^{a(n)} \vert \mu_{\ell} \vert^{r} \Big )^{\frac{1}{r}} \leq 1  \Big \} \\
& \qquad \leq \sup \Big \{ \sum_{\ell =1}^{a(n)} \vert \mu_{\ell} \vert  \Big \Vert   (c^{n}_{2, \ell})_{\lambda_{2}} \Big \vert  h^{0,m}_{r,\delta_{2}}(\nabla_{\sigma_{2}}) \Big \Vert  : \Big ( \sum_{\ell =1}^{a(n)} \vert \mu_{\ell} \vert^{r} \Big )^{\frac{1}{r}} \leq 1  \Big \} .
\end{align*}
For each $ l \in \{ 1, 2, \ldots , a(n) \}   $ we observe
\begin{align*}
\vert \mu_{l} \vert = ( \vert \mu_{l} \vert^{r} )^{\frac{1}{r}} \leq \Big ( \sum_{\ell =1}^{a(n)} \vert \mu_{\ell} \vert^{r} \Big )^{\frac{1}{r}} \leq 1 .
\end{align*}
Consequently we also find
\begin{align*}
& \sup \Big \{  \Vert  (c^{n}_{2})_{\lambda_{2}} \vert  h^{0,m}_{r,\delta_{2}}(\nabla_{\sigma_{2}})  \Vert  : \Big ( \sum_{\ell =1}^{a(n)} \vert \mu_{\ell} \vert^{r} \Big )^{\frac{1}{r}} \leq 1  \Big \} \leq   \sum_{\ell =1}^{a(n)}   \Big \Vert   (c^{n}_{2, \ell})_{\lambda_{2}} \Big \vert  h^{0,m}_{r,\delta_{2}}(\nabla_{\sigma_{2}}) \Big \Vert   .
\end{align*}
Applying \eqref{qua_rep_d=2_1_b} this yields
\begin{align*}
& \Vert f  \vert (H^{s}_{r}(I))_{\sigma_{1}} \otimes_{g_{r}} L_{r}(I) \Vert \\
& \qquad \lesssim \lim_{n \rightarrow \infty}  \Big \{  \Big ( \sum_{\ell =1}^{a(n)}   \Vert  (c^{n}_{1, \ell})_{\lambda_{1}} \vert  h^{s,m}_{r,\delta_{1}}(\nabla_{\sigma_{1}})  \Vert^{r}   \Big )^{\frac{1}{r}}  \times \Big (  \sum_{\ell =1}^{a(n)}   \Big \Vert   (c^{n}_{2, \ell})_{\lambda_{2}} \Big \vert  h^{0,m}_{r,\delta_{2}}(\nabla_{\sigma_{2}}) \Big \Vert \Big ) \Big \} \\
& \qquad < \infty .
\end{align*}
Thus we have $ f \in (H^{s}_{r}(I))_{\sigma_{1}} \otimes_{g_{r}} L_{r}(I)$. To obtain the equivalent norm we can use similar methods. This time in \eqref{eq_Prop6_Step3_est11} we have to work with the infimum over all possible sequences of functions $(h_{n})_{n \in \mathbb{N}}$ of the form \eqref{qua_rep_d=2_1_f}, at which we apply decompositions \eqref{qua_rep_d=2_1_g} and \eqref{qua_rep_d=2_1_h}. To this end once again we can apply Theorem \ref{mainresult1_Hsr_d=1} and Corollary \ref{mainresult1_Lr_d=1}. Putting all together the proof is complete. 
\end{proof}

The counterpart of Proposition \ref{prop_char_Hs_Lr} for the function space $ L_{r}(I)    \otimes_{g_{r}} (H^{s}_{r}(I))_{\sigma_{2}}  $ reads as follows.

\begin{prop}\label{prop_char_Lr_Hs}
Let $d = 2$. Let $  1 < r < \infty $ and $ m \in \mathbb{N}  $ with $ m \geq 2  $. Let $ \tilde{m} \in \mathbb{N}_{0}   $ with $ m + \tilde{m} \in 2 \mathbb{N}   $ and $  \tilde{m} > 5 m + 12 $. Let $ 0 < s  < m - 1 $ and $\delta_{1}, \delta_{2} > 1$. Let $ \sigma_{1} = ( 0 , 0 )  $ and $ \sigma_{2}  $ be a boundary condition. Put $  \boldsymbol{{\sigma}^{2}} := ( \sigma_{1}, \sigma_{2} )  $. Let $f \in L_{r}((0,1)^{2})$. Then we have $ f \in L_{r}(I)    \otimes_{g_{r}} (H^{s}_{r}(I))_{\sigma_{2}} $ if and only if there exists a sequence of functions $(h_{n})_{n \in \mathbb{N}}$ of the form 
\begin{equation}\label{qua_rep_d=2_2_a}
h_{n} = \sum_{( \lambda_{1}, \lambda_{2}) \in \boldsymbol{\nabla_{{\sigma}^2}} }  \Big (  \sum_{\ell = 1}^{a(n)}    (c^{n}_{1, \ell})_{\lambda_{1} }    (c^{n}_{2, \ell})_{\lambda_{2}}   \Big )  (  \psi_{\lambda_{1}}^{\sigma_{1}}  \otimes  \psi_{\lambda_{2}}^{\sigma_{2}} )
\end{equation}
with $a(n) \in \mathbb{N}$, such that
\begin{equation}\label{qua_rep_d=2_2_b}
\Vert  \boldsymbol{c}^{2}  \vert  h^{s,m}_{r, \boldsymbol{\delta}}(\boldsymbol{\nabla_{{\sigma}^{2}}} , 2)  \Vert :=  \Big ( \sum_{\ell =1}^{a(n)}  \Vert  (c^{n}_{1, \ell})_{\lambda_{1}} \vert  h^{0,m}_{r,\delta_{1}}(\nabla_{\sigma_{1}})  \Vert   \Big ) \times \Big ( \sum_{\ell =1}^{a(n)}  \Vert  (c^{n}_{2, \ell})_{\lambda_{2}} \vert  h^{s,m}_{r,\delta_{2}}(\nabla_{\sigma_{2}})  \Vert   \Big )   < \infty
\end{equation}
and
\begin{equation}\label{qua_rep_d=2_2_d}
\lim_{n \rightarrow \infty} \Vert f - h_{n} \vert   L_{r}(I)  \otimes_{g_{r}} (H^{s}_{r}(I))_{\sigma_{2}} \Vert = 0.
\end{equation}
Moreover, the norms $ \Vert f  \vert  L_{r}(I) \otimes_{g_{r}} (H^{s}_{r}(I))_{\sigma_{2}} \Vert$ and 
\begin{equation}\label{qua_rep_d=2_2_e}
\begin{split}
&\Vert f  \vert  L_{r}(I) \otimes_{g_{r}} (H^{s}_{r}(I))_{\sigma_{2}}  \Vert^{( \boldsymbol{\delta}, \boldsymbol{\nabla_{{\sigma}^{2}}})} \\
& \qquad \qquad := \lim_{n \rightarrow \infty} \inf_{\eqref{qua_rep_d=2_2_f}} \Big \{  \Big ( \sum_{\ell =1}^{a(n)}  \inf_{\eqref{qua_rep_d=2_2_g}} \Vert  (c^{n}_{1, \ell})_{\lambda_{1}} \vert  h^{0,m}_{r,\delta_{1}}(\nabla_{\sigma_{1}})  \Vert^{r}   \Big )^{\frac{1}{r}} \\
& \qquad \qquad \qquad \qquad \times \sup \Big \{ \inf_{\eqref{qua_rep_d=2_2_h}} \Vert  (c^{n}_{2})_{\lambda_{2}} \vert  h^{s,m}_{r,\delta_{2}}(\nabla_{\sigma_{2}})  \Vert  : \Big ( \sum_{\ell =1}^{a(n)} \vert \mu_{\ell} \vert^{r} \Big )^{\frac{1}{r}} \leq 1  \Big \} \Big \} 
\end{split}
\end{equation}
are equivalent. Here the infimum \eqref{qua_rep_d=2_2_f} is taken over all sequences of functions $(h_{n})_{n \in \mathbb{N}}$ of the form 
\begin{equation}\label{qua_rep_d=2_2_f}
h_{n} = \sum_{\ell = 1}^{a(n)} (f^{n})_{\ell} \otimes (g^{n})_{\ell} , \qquad (f^{n})_{\ell} \in L_{r}(I), \; (g^{n})_{\ell} \in (H^{s}_{r}(I))_{\sigma_{2}} 
\end{equation}
with 
\begin{align*}
\lim_{n \rightarrow \infty} \Vert f - h_{n} \vert L_{r}(I) \otimes_{g_{r}} (H^{s}_{r}(I))_{\sigma_{2}} \Vert = 0.
\end{align*}
The infimum \eqref{qua_rep_d=2_2_g} is taken over all sequences $ \{ (c^{n}_{1, \ell})_{\lambda_{1}}  \}_{\lambda_{1} \in \nabla_{\sigma_{1}}}   $ such that
\begin{equation}\label{qua_rep_d=2_2_g}
(f^{n})_{\ell} = \sum_{\lambda_{1} \in \nabla_{\sigma_{1}}}   (c^{n}_{1, \ell})_{\lambda_{1}} \psi_{\lambda_{1}}^{\sigma_{1}} .
\end{equation}
The infimum \eqref{qua_rep_d=2_2_h} is taken over all sequences $ \{ (c^{n}_{2})_{\lambda_{2}}  \}_{\lambda_{2} \in \nabla_{\sigma_{2}}}   $ such that
\begin{equation}\label{qua_rep_d=2_2_h}
\sum_{\ell =1}^{a(n)} \mu_{\ell} (g^{n})_{\ell} = \sum_{\lambda_{2} \in \nabla_{\sigma_{2}}}   (c^{n}_{2})_{\lambda_{2}} \psi_{\lambda_{2}}^{\sigma_{2}} .
\end{equation}
\end{prop}

\begin{proof}
This result can be proved with similar methods as used for the proof of Proposition \ref{prop_char_Hs_Lr}. Again we can apply Theorem \ref{mainresult1_Hsr_d=1} and Corollary \ref{mainresult1_Lr_d=1}. This time the roles of the spaces $ (H^{s}_{r}(I))_{\sigma_{2}}  $ and $ L_{r}(I)  $ must be interchanged. No further modifications are necessary. 
\end{proof}

Now we are well-prepared to describe bivariate Bessel-Potential spaces $ (H^{s}_{r}((0,1)^{2}))_{\boldsymbol{\sigma}} $ in terms of bivariate quarklets. For that purpose we apply Proposition \ref{prop_char_Hs_Lr} and Proposition \ref{prop_char_Lr_Hs} and combine them by using Proposition \ref{prop_tensor_Hsr_d2} and Proposition \ref{prop_tensor_Hsr_d2_boundco}. We obtain the following theorem. It is one of the main results of this paper.

\begin{satz}\label{THM_d2_Hsr2_main1}
Let $d = 2$. Let $  1 < r < \infty $ and $ m \in \mathbb{N}  $ with $ m \geq 2  $. Let $ \tilde{m} \in \mathbb{N}_{0}   $ with $ m + \tilde{m} \in 2 \mathbb{N}   $ and $  \tilde{m} > 5 m + 12 $. Let $ 0 < s  < m - 1 $ and $\delta_{1}, \delta_{2} > 1$. Let $ \boldsymbol{\sigma} = ( \sigma_{1}, \sigma_{2} )  $ be a collection of boundary conditions. Let $f \in L_{r}((0,1)^{2})$. Then we have $ f \in (H^{s}_{r}((0,1)^{2}))_{\boldsymbol{\sigma} } $ if and only if there exists a sequence of functions $(h_{n})_{n \in \mathbb{N}}$ of the form 
\begin{equation}\label{qua_rep_d=2_THM1_a}
h_{n} := h_{n}^{1} + h_{n}^{2} := \sum_{( \lambda_{1}, \lambda_{2}) \in \boldsymbol{\nabla_{\sigma}} }  \Big (  \sum_{\ell = 1}^{a(n)}    (c^{n}_{1, \ell})_{\lambda_{1} }    (c^{n}_{2, \ell})_{\lambda_{2}}  + \sum_{\ell = 1}^{b(n)}    (d^{n}_{1, \ell})_{\lambda_{1} }    (d^{n}_{2, \ell})_{\lambda_{2}}  \Big )  (  \psi_{\lambda_{1}}^{\sigma_{1}}  \otimes  \psi_{\lambda_{2}}^{\sigma_{2}} )
\end{equation}
with $a(n), b(n) \in \mathbb{N}$, such that 
\begin{equation}\label{qua_rep_d=2_THM1_b}
\begin{split}
 \Vert  \boldsymbol{c} \vert  h^{s,m}_{r, \boldsymbol{\delta}}(\boldsymbol{\nabla_{\sigma}})  \Vert & := \Vert  \boldsymbol{c}^{1}  \vert  h^{s,m}_{r, \boldsymbol{\delta}}(\boldsymbol{\nabla_{\sigma}} , 1)  \Vert +  \Vert  \boldsymbol{c}^{2}  \vert  h^{s,m}_{r, \boldsymbol{\delta}}(\boldsymbol{\nabla_{\sigma}} , 2)  \Vert   \\
&  = \Big ( \sum_{\ell =1}^{a(n)}  \Vert  (c^{n}_{1, \ell})_{\lambda_{1}} \vert  h^{s,m}_{r,\delta_{1}}(\nabla_{\sigma_{1}})  \Vert   \Big ) \times \Big ( \sum_{\ell =1}^{a(n)}  \Vert  (c^{n}_{2, \ell})_{\lambda_{2}} \vert  h^{0,m}_{r,\delta_{2}}(\nabla_{\sigma_{2}})  \Vert   \Big )  \\
& \qquad  + \Big ( \sum_{\ell =1}^{b(n)}  \Vert  (d^{n}_{1, \ell})_{\lambda_{1}} \vert  h^{0,m}_{r,\delta_{1}}(\nabla_{\sigma_{1}})  \Vert   \Big )  \times \Big ( \sum_{\ell =1}^{b(n)}  \Vert  (d^{n}_{2, \ell})_{\lambda_{2}} \vert  h^{s,m}_{r,\delta_{2}}(\nabla_{\sigma_{2}})  \Vert   \Big )  < \infty
\end{split}
\end{equation}
and
\begin{equation}\label{qua_rep_d=2_THM1_d}
\lim_{n \rightarrow \infty} \Big \Vert \frac{1}{2} f - h_{n}^{1} \Big \vert (H^{s}_{r}(I))_{\sigma_{1}} \otimes_{g_{r}} L_{r}(I) \Big \Vert + \lim_{n \rightarrow \infty} \Big \Vert \frac{1}{2} f - h_{n}^{2} \Big \vert   L_{r}(I)  \otimes_{g_{r}} (H^{s}_{r}(I))_{\sigma_{2}} \Big \Vert = 0 .
\end{equation}
Moreover, the norms $ \Vert f \vert  (H^{s}_{r}((0,1)^{2}))_{\boldsymbol{\sigma}}  \Vert  $ and
\begin{equation}\label{eq_new_norm_d2}
\begin{split}
& \Vert f \vert  (H^{s}_{r}((0,1)^{2}))_{\boldsymbol{\sigma}}  \Vert^{( \boldsymbol{\delta}, \boldsymbol{\nabla_{\sigma}})}  
\\
& \qquad \qquad :=  \Vert f  \vert (H^{s}_{r}(I))_{\sigma_{1}} \otimes_{g_{r}} L_{r}(I) \Vert^{( \boldsymbol{\delta}, \boldsymbol{\nabla_{{\sigma}^{1}}})} + \Vert f  \vert  L_{r}(I) \otimes_{g_{r}} (H^{s}_{r}(I))_{\sigma_{2}}  \Vert^{( \boldsymbol{\delta}, \boldsymbol{\nabla_{{\sigma}^{2}}})}   \\
& \qquad \qquad  = \lim_{n \rightarrow \infty} \inf_{\eqref{qua_rep_d=2_1_f}} \Big \{  \Big ( \sum_{\ell =1}^{a(n)}  \inf_{\eqref{qua_rep_d=2_1_g}} \Vert  (c^{n}_{1, \ell})_{\lambda_{1}} \vert  h^{s,m}_{r,\delta_{1}}(\nabla_{\sigma_{1}})  \Vert^{r}   \Big )^{\frac{1}{r}} \\
& \qquad \qquad \qquad \qquad  \times \sup \Big \{ \inf_{\eqref{qua_rep_d=2_1_h}} \Vert  (c^{n}_{2})_{\lambda_{2}} \vert  h^{0,m}_{r,\delta_{2}}(\nabla_{(0,0)})  \Vert  : \Big ( \sum_{\ell =1}^{a(n)} \vert \mu_{\ell} \vert^{r} \Big )^{\frac{1}{r}} \leq 1  \Big \} \Big \}  \\
& \qquad \qquad   + \lim_{n \rightarrow \infty} \inf_{\eqref{qua_rep_d=2_2_f}} \Big \{  \Big ( \sum_{\ell =1}^{b(n)}  \inf_{\eqref{qua_rep_d=2_2_g}} \Vert  (d^{n}_{1, \ell})_{\lambda_{1}} \vert  h^{0,m}_{r,\delta_{1}}(\nabla_{(0,0)})  \Vert^{r}   \Big )^{\frac{1}{r}} \\
& \qquad \qquad \qquad \qquad  \times \sup \Big \{ \inf_{\eqref{qua_rep_d=2_2_h}} \Vert  (d^{n}_{2})_{\lambda_{2}} \vert  h^{s,m}_{r,\delta_{2}}(\nabla_{\sigma_{2}})  \Vert  : \Big ( \sum_{\ell =1}^{b(n)} \vert \mu_{\ell} \vert^{r} \Big )^{\frac{1}{r}} \leq 1  \Big \} \Big \} 
\end{split}
\end{equation} 
are equivalent. Here the infimum \eqref{qua_rep_d=2_1_f} is taken over all sequences of functions $(h^{1}_{n})_{n \in \mathbb{N}}$ of the form \eqref{qua_rep_d=2_1_f}
with 
\begin{align*}
\lim_{n \rightarrow \infty} \Vert f - h^{1}_{n} \vert (H^{s}_{r}(I))_{\sigma_{1}} \otimes_{g_{r}} L_{r}(I) \Vert = 0.
\end{align*}
The infimum \eqref{qua_rep_d=2_1_g} is taken over all sequences $ \{ (c^{n}_{1, \ell})_{\lambda_{1}}  \}_{\lambda_{1} \in \nabla_{\sigma_{1}}}   $ such that \eqref{qua_rep_d=2_1_g} holds. The infimum \eqref{qua_rep_d=2_1_h} is taken over all sequences $ \{ (c^{n}_{2})_{\lambda_{2}}  \}_{\lambda_{2} \in \nabla_{(0,0)}}   $ such that \eqref{qua_rep_d=2_1_h} is fulfilled. For the second term the infimum \eqref{qua_rep_d=2_2_f} is taken over all sequences of functions $(h^{2}_{n})_{n \in \mathbb{N}}$ of the form \eqref{qua_rep_d=2_2_f} with 
\begin{align*}
\lim_{n \rightarrow \infty} \Vert f - h^{2}_{n} \vert L_{r}(I) \otimes_{g_{r}} (H^{s}_{r}(I))_{\sigma_{2}} \Vert = 0,
\end{align*}
at which $a(n)$ is replaced by $b(n) \in \mathbb{N}$. The infimum \eqref{qua_rep_d=2_2_g} is taken over all sequences $ \{ (d^{n}_{1, \ell})_{\lambda_{1}}  \}_{\lambda_{1} \in \nabla_{(0,0)}}   $ such that \eqref{qua_rep_d=2_2_g} holds with $(c^{n}_{1, \ell})$ replaced by $(d^{n}_{1, \ell})$. The infimum \eqref{qua_rep_d=2_2_h} is taken over all sequences $ \{ (d^{n}_{2})_{\lambda_{2}}  \}_{\lambda_{2} \in \nabla_{\sigma_{2}}}   $ such that \eqref{qua_rep_d=2_2_h} is fulfilled with $(c^{n}_{2})$ replaced by $(d^{n}_{2})$.
\end{satz}

\begin{proof}
\textit{Step 1.} Let $ f \in (H^{s}_{r}((0,1)^2))_{\boldsymbol{\sigma}}  $. We want to prove that there exists a sequence of functions $(h_{n})_{n \in \mathbb{N}}$ of the form \eqref{qua_rep_d=2_THM1_a}, such that \eqref{qua_rep_d=2_THM1_b} and \eqref{qua_rep_d=2_THM1_d} are fulfilled. For that purpose we can use Proposition \ref{prop_tensor_Hsr_d2} and Proposition \ref{prop_tensor_Hsr_d2_boundco} to find $   f \in ( (H^{s}_{r}(I))_{\sigma_{1}} \otimes_{g_{r}} L_{r}(I) ) \cap ( L_{r}(I) \otimes_{g_{r}}  (H^{s}_{r}(I))_{\sigma_{2}} ) $. Consequently we can apply versions of Proposition \ref{prop_char_Hs_Lr} and Proposition \ref{prop_char_Lr_Hs} that are adapted to the given boundary conditions to describe the function $  f = \frac{1}{2} f + \frac{1}{2} f   $. We obtain sequences of functions $ ( h_{n} )_{n \in \mathbb{N}}  $, $ ( h_{n}^{1} )_{n \in \mathbb{N}}  $ and $ ( h_{n}^{2} )_{n \in \mathbb{N}}  $ of the form
\begin{align*}
h_{n}  = h_{n}^{1} + h_{n}^{2} & = \sum_{( \lambda_{1}, \lambda_{2}) \in \boldsymbol{\nabla_{\sigma}} }  \Big (  \sum_{\ell = 1}^{a(n)}    (c^{n}_{1, \ell})_{\lambda_{1} }    (c^{n}_{2, \ell})_{\lambda_{2}}   \Big )  (  \psi_{\lambda_{1}}^{\sigma_{1}}  \otimes  \psi_{\lambda_{2}}^{\sigma_{2}} ) \\
& \qquad + \sum_{( \lambda_{1}, \lambda_{2}) \in \boldsymbol{\nabla_{\sigma}} }  \Big (  \sum_{\ell = 1}^{b(n)}    (d^{n}_{1, \ell})_{\lambda_{1} }    (d^{n}_{2, \ell})_{\lambda_{2}}   \Big )  (  \psi_{\lambda_{1}}^{\sigma_{1}}  \otimes  \psi_{\lambda_{2}}^{\sigma_{2}} )
\end{align*}
with $ a(n), b(n) \in \mathbb{N}   $, such that \eqref{qua_rep_d=2_THM1_b} holds and we have
\begin{align*}
\lim_{n \rightarrow \infty} \Big \Vert \frac{1}{2} f - h_{n}^{1} \Big \vert (H^{s}_{r}(I))_{\sigma_{1}} \otimes_{g_{r}} L_{r}(I) \Big \Vert + \lim_{n \rightarrow \infty} \Big \Vert \frac{1}{2} f - h_{n}^{2} \Big \vert   L_{r}(I)  \otimes_{g_{r}} (H^{s}_{r}(I))_{\sigma_{2}} \Big \Vert = 0 .
\end{align*}
This completes Step 1 of the proof.

\textit{Step 2.} Now let $f \in L_{r}((0,1)^{2})$. Let $ ( h_{n} )_{n \in \mathbb{N}}  $, $ ( h_{n}^{1} )_{n \in \mathbb{N}}  $ and $ ( h_{n}^{2} )_{n \in \mathbb{N}}  $ be sequences of functions, such that \eqref{qua_rep_d=2_THM1_a}, \eqref{qua_rep_d=2_THM1_b} and \eqref{qua_rep_d=2_THM1_d} are fulfilled. Then Proposition \ref{prop_char_Hs_Lr} applied for $ ( h_{n}^{1} )_{n \in \mathbb{N}}  $ yields $ f \in  (H^{s}_{r}(I))_{\sigma_{1}} \otimes_{g_{r}} L_{r}(I)  $. Moreover, we can use Proposition \ref{prop_char_Lr_Hs} for $ ( h_{n}^{2} )_{n \in \mathbb{N}}  $ to find $ f \in L_{r}(I)    \otimes_{g_{r}} (H^{s}_{r}(I))_{\sigma_{2}} $. Consequently we also get $ f \in ((H^{s}_{r}(I))_{\sigma_{1}} \otimes_{g_{r}} L_{r}(I)) \cap ( L_{r}(I)    \otimes_{g_{r}} (H^{s}_{r}(I))_{\sigma_{2}}   )$. Then Proposition \ref{prop_tensor_Hsr_d2} and Proposition \ref{prop_tensor_Hsr_d2_boundco} imply $ f \in (H^{s}_{r}((0,1)^{2}))_{\boldsymbol{\sigma}} $.

\textit{Step 3.}
Now we prove the equivalence of the norms. For that purpose Proposition \ref{prop_tensor_Hsr_d2} and Proposition \ref{prop_tensor_Hsr_d2_boundco} yield that we have 
\begin{align*}
(H^{s}_{r}((0,1)^{2}))_{\boldsymbol{\sigma}} = ( (H^{s}_{r}(I))_{\sigma_{1}} \otimes_{g_{r}} L_{r}(I) ) \cap ( L_{r}(I) \otimes_{g_{r}}  (H^{s}_{r}(I))_{\sigma_{2}} )
\end{align*}
with equivalent norms. Consequently using Definition \ref{defi_norm_intersec} we find 
\begin{align*}
\Vert f \vert  (H^{s}_{r}((0,1)^{2}))_{\boldsymbol{\sigma}} \Vert &  \approx \Vert f \vert  ( (H^{s}_{r}(I))_{\sigma_{1}} \otimes_{g_{r}} L_{r}(I) ) \cap ( L_{r}(I) \otimes_{g_{r}}  (H^{s}_{r}(I))_{\sigma_{2}} )  \Vert \\ 
&  = \Vert f \vert   (H^{s}_{r}(I))_{\sigma_{1}} \otimes_{g_{r}} L_{r}(I)   \Vert + \Vert f \vert    L_{r}(I) \otimes_{g_{r}}  (H^{s}_{r}(I))_{\sigma_{2}}   \Vert .
\end{align*}
Now for the first term we can apply Proposition \ref{prop_char_Hs_Lr}. Then we obtain
\begin{align*}
& \Vert f \vert   (H^{s}_{r}(I))_{\sigma_{1}} \otimes_{g_{r}} L_{r}(I)   \Vert \\
& \qquad \approx \lim_{n \rightarrow \infty} \inf_{\eqref{qua_rep_d=2_1_f}} \Big \{  \Big ( \sum_{\ell =1}^{a(n)}  \inf_{\eqref{qua_rep_d=2_1_g}} \Vert  (c^{n}_{1, \ell})_{\lambda_{1}} \vert  h^{s,m}_{r,\delta_{1}}(\nabla_{\sigma_{1}})  \Vert^{r}   \Big )^{\frac{1}{r}} \\
& \qquad \qquad  \times \sup \Big \{ \inf_{\eqref{qua_rep_d=2_1_h}} \Vert  (c^{n}_{2})_{\lambda_{2}} \vert  h^{0,m}_{r,\delta_{2}}(\nabla_{(0,0)})  \Vert  : \Big ( \sum_{\ell =1}^{a(n)} \vert \mu_{\ell} \vert^{r} \Big )^{\frac{1}{r}} \leq 1  \Big \} \Big \}  .
\end{align*}
Here the infimum \eqref{qua_rep_d=2_1_f} is taken over all sequences of functions $(h^{1}_{n})_{n \in \mathbb{N}}$ of the form \eqref{qua_rep_d=2_1_f}
with 
\begin{align*}
\lim_{n \rightarrow \infty} \Vert f - h^{1}_{n} \vert (H^{s}_{r}(I))_{\sigma_{1}} \otimes_{g_{r}} L_{r}(I) \Vert = 0.
\end{align*}
The infimum \eqref{qua_rep_d=2_1_g} is taken over all sequences $ \{ (c^{n}_{1, \ell})_{\lambda_{1}}  \}_{\lambda_{1} \in \nabla_{\sigma_{1}}}   $ such that \eqref{qua_rep_d=2_1_g} holds. The infimum \eqref{qua_rep_d=2_1_h} is taken over all sequences $ \{ (c^{n}_{2})_{\lambda_{2}}  \}_{\lambda_{2} \in \nabla_{(0,0)}}   $ such that \eqref{qua_rep_d=2_1_h} is fulfilled. For the second term we can use Proposition \ref{prop_char_Lr_Hs}. It yields
\begin{align*}
& \Vert f \vert    L_{r}(I) \otimes_{g_{r}}  (H^{s}_{r}(I))_{\sigma_{2}}   \Vert \\
& \qquad  \approx \lim_{n \rightarrow \infty} \inf_{\eqref{qua_rep_d=2_2_f}} \Big \{  \Big ( \sum_{\ell =1}^{b(n)}  \inf_{\eqref{qua_rep_d=2_2_g}} \Vert  (d^{n}_{1, \ell})_{\lambda_{1}} \vert  h^{0,m}_{r,\delta_{1}}(\nabla_{(0,0)})  \Vert^{r}   \Big )^{\frac{1}{r}} \\
& \qquad \qquad  \times \sup \Big \{ \inf_{\eqref{qua_rep_d=2_2_h}} \Vert  (d^{n}_{2})_{\lambda_{2}} \vert  h^{s,m}_{r,\delta_{2}}(\nabla_{\sigma_{2}})  \Vert  : \Big ( \sum_{\ell =1}^{b(n)} \vert \mu_{\ell} \vert^{r} \Big )^{\frac{1}{r}} \leq 1  \Big \} \Big \} . 
\end{align*}
Here the infimum \eqref{qua_rep_d=2_2_f} is taken over all sequences of functions $(h^{2}_{n})_{n \in \mathbb{N}}$ of the form \eqref{qua_rep_d=2_2_f} with 
\begin{align*}
\lim_{n \rightarrow \infty} \Vert f - h^{2}_{n} \vert L_{r}(I) \otimes_{g_{r}} (H^{s}_{r}(I))_{\sigma_{2}} \Vert = 0,
\end{align*}
at which $a(n)$ is replaced by $b(n) \in \mathbb{N}$. The infimum \eqref{qua_rep_d=2_2_g} is taken over all sequences $ \{ (d^{n}_{1, \ell})_{\lambda_{1}}  \}_{\lambda_{1} \in \nabla_{(0,0)}}   $ such that \eqref{qua_rep_d=2_2_g} holds with $(c^{n}_{1, \ell})$ replaced by $(d^{n}_{1, \ell})$. The infimum \eqref{qua_rep_d=2_2_h} is taken over all sequences $ \{ (d^{n}_{2})_{\lambda_{2}}  \}_{\lambda_{2} \in \nabla_{\sigma_{2}}}   $ such that \eqref{qua_rep_d=2_2_h} is fulfilled with $(c^{n}_{2})$ replaced by $(d^{n}_{2})$. Combining both equivalences we get
\begin{align*}
\Vert f \vert  (H^{s}_{r}((0,1)^{2}))_{\boldsymbol{\sigma}}  \Vert   \approx   \Vert f \vert  (H^{s}_{r}((0,1)^{2}))_{\boldsymbol{\sigma}}  \Vert^{( \boldsymbol{\delta}, \boldsymbol{\nabla_{\sigma}})} ,
\end{align*}
see \eqref{eq_new_norm_d2} for an explanation concerning the notation. The proof is complete.
\end{proof}

\begin{rem}
Theorem \ref{THM_d2_Hsr2_main1} also provides a characterization in terms of bivariate wavelets for the spaces $  (H^{s}_{r}((0,1)^{2}))_{\boldsymbol{\sigma}} $. These wavelets are constructed out of the univariate boundary adapted Primbs wavelets using the tensor product methods described above. To write down the corresponding result in the formulation of Theorem \ref{THM_d2_Hsr2_main1} one has to replace $ \nabla_{\sigma_{i}}  $ by $ \nabla_{\sigma_{i}}^{P} $, see \eqref{index_Primbs1_full1}. That means we only deal with the Primbs basis index set. The proof can be carried out in a similar way as that for Theorem \ref{THM_d2_Hsr2_main1}, whereat in each step we only have to work with the polynomial degree $ p_{i} = 0  $. For that purpose we use the wavelet counterparts of Theorem \ref{mainresult1_Hsr_d=1} and Corollary \ref{mainresult1_Lr_d=1}.   
\end{rem}

\begin{rem}
Theorem \ref{THM_d2_Hsr2_main1} also can be used to characterize functions $  f \in (H^{s}_{r}((0,1)^{2}))_{\boldsymbol{\sigma}} $ that can be described via a sequence $(h_{n})_{n \in \mathbb{N}}$ such that $ a(n) = b(n) = 1  $ for all $n \in \mathbb{N}$. In that case the expressions \eqref{qua_rep_d=2_THM1_a} and \eqref{qua_rep_d=2_THM1_b} simplify. Moreover, then similarities between \eqref{qua_rep_d=2_THM1_b} and the weights given in Theorem 3.11 in \cite{DaFKRaa} become apparent. However, in \cite{DaFKRaa} only the Hilbert space case $r = 2$ has been investigated, which explains the easier structure of the results proved there. 
\end{rem}

\begin{rem}
Let us remark that the bivariate tensor quarklets constructed in Theorem \ref{THM_d2_Hsr2_main1} have a highly anisotropic structure. Indeed, our tensor product quarklets can be interpreted as wavelet version of sparse grids enriched with polynomials. Consequently in the long run most likely they can be employed to design numerical approximation methods with dimension-independent convergence rates. The quarklet system obtained in Theorem \ref{THM_d2_Hsr2_main1} is profoundly redundant. This will be essential later on when it comes to the development of adaptive approximation techniques using the tensor quarklets. For example they can be utilized to design adaptive bivariate quarklet tree approximation techniques, whereby in each step of the algorithm either a space refinement or a polynomial enrichment can be carried out. Some first numerical experiments using our bivariate tensor quarklets can be found in \cite{DaFKRaa}, see Section 5.
\end{rem}

\section{Summary and Further Issues}\label{sec_summary}

In this paper we have seen that it is possible to describe bivariate Bessel-Potential spaces in terms of quarklets which are constructed out of univariate functions via tensor product methods. For that purpose we used a recent result of Hansen and Sickel, see \cite{HaSi2022}, which implies that bivariate Bessel-Potential spaces defined on squares can be written as an intersection of function spaces that are tensor products of univariate spaces. Looking at our main result Theorem \ref{THM_d2_Hsr2_main1} several new continuative questions arise. Some of them which also might be subject of future research are collected in the following list.

\begin{itemize}
\item[(i)] In this paper we only obtain quarklet characterizations for the bivariate Bessel-Potential spaces $ H^{s}_{r}((0,1)^2)  $ with $ d = 2  $. Consequently the natural question arises how multivariate counterparts of Theorem \ref{THM_d2_Hsr2_main1} for the spaces $ H^{s}_{r}((0,1)^d)  $ with $ d \in \mathbb{N}  $ and $ d \geq 3  $ look like. In connection with that one requires a multivariate version of Proposition \ref{prop_tensor_Hsr_d2}, see \cite{HaSi2022} for some first attempts. Having this at hand one can iterate the arguments used in the proofs of Proposition \ref{prop_char_Hs_Lr} and Theorem \ref{THM_d2_Hsr2_main1} to obtain multivariate quarklet characterizations via tensor products. 

\item[(ii)] For $ 1 < r < \infty    $ the Bessel-Potential spaces can be interpreted as Triebel-Lizorkin spaces via $  H^{s}_{r}((0,1)^d) = F^{s}_{r,2}((0,1)^d)   $, see \cite{Tr83}, \cite{Tr92} and \cite{Tr06} for definitions and explanations. So it would be interesting to have characterizations in terms of multivariate quarklets using tensor product methods for the spaces $    F^{s}_{r,q}((0,1)^d)$ with $ 0 < r < \infty   $ and $ 0 < q \leq \infty   $. To this end a counterpart of Proposition \ref{prop_tensor_Hsr_d2} for Triebel-Lizorkin spaces is required. However, as far as the author knows, the shape of such decompositions for the case $  q \not = 2  $ is still unknown.  

\item[(iii)] In this paper we obtained characterizations in terms of quarklets for Bessel-Potential spaces defined on $ (0,1)^2   $. However, it also would be interesting to have such characterizations for function spaces that are defined on more complicated domains. For example, in \cite{DaFKRaa} quarklets showed up in the context of Sobolev spaces $ H^{s}(\Omega) = H^{s}_{2}(\Omega)  $, where $ \Omega \subset \mathbb{R}^2   $ is a domain that can be decomposed into squares, whereby special emphasis is given to the case that $ \Omega  $ is an L-shaped domain. Inspired by that it will be the subject of future research to deduce quarklet characterizations for the spaces $ H^{s}_{r}(\Omega)  $ with $ 1 < r < \infty   $ for L-shaped domains $\Omega$.

\end{itemize}

\vspace{0,3 cm}

\textbf{Acknowledgment.} This work has been supported by Deutsche Forschungsgemeinschaft (DFG), grant HO 7444/1-1 with project number 528343051. The author thanks Stephan Dahlke and Dorian Vogel for several valuable discussions.

\end{document}